\documentclass[12pt]{amsart}

\usepackage{geometry}
\usepackage{amsfonts}
\usepackage{amsmath}
\usepackage{amssymb}
\usepackage{amsthm}
\usepackage{mathtools}
\usepackage{mathrsfs}
\usepackage{bbm}
\usepackage{graphicx}
\usepackage{epstopdf}
\usepackage{xcolor}
\usepackage{leftindex}
\usepackage[normalem]{ulem}
\usepackage[lite]{amsrefs}
\usepackage[colorlinks]{hyperref}

\newtheorem{theorem}{Theorem}[section]
\newtheorem{corollary}[theorem]{Corollary}
\newtheorem{lemma}[theorem]{Lemma}
\newtheorem{proposition}[theorem]{Proposition}
\newtheorem*{thm}{Theorem}
\newtheorem{notation}[theorem]{Notation}

\theoremstyle{definition}
\newtheorem{definition}[theorem]{Definition}
\newtheorem{example}[theorem]{Example}

\theoremstyle{remark}
\newtheorem{remark}[theorem]{Remark}

\numberwithin{equation}{section}
\numberwithin{figure}{section}

\usepackage{tikz-cd}
\usepackage{tikz}
\usepackage{tikz-3dplot}

\usepgfmodule{decorations}
\usetikzlibrary{
    arrows.meta,
    calc,
    decorations.markings,
    hobby,
    knots
}

\newenvironment{tric}
    {\begin{tikzpicture}[scale= 0.5, semithick,draw=black,double distance=0.25mm,
        baseline={([yshift=-.8ex]current bounding box.center)}] }
    {\end{tikzpicture}}

\colorlet{darkblue}{blue!70!black}
\colorlet{darkred}{red!70!black}
\colorlet{darkgreen}{green!70!black}

\newcommand{\GL}{\mathfrak{gl}}
\newcommand{\Z}{\mathbb{Z}}
\newcommand{\C}{\mathbb{C}}
\newcommand{\R}{\mathbb{R}}

\newcommand{\webcat}{\mathbf{Web}}
\newcommand{\Web}{\mathbf{W}_N}
\newcommand{\Hom}{\mathrm{Hom}}
\newcommand{\Fund}{\mathbf{Fund}}
\newcommand{\Kar}{\operatorname{Kar}}
\newcommand{\Li}[2]{\overline{L^{#1}_{#2}}}
\newcommand{\Ri}[2]{\overline{R^{#1}_{#2}}}
\newcommand{\Wsix}{\mathcal{W}_6}
\newcommand{\Wzero}{\mathcal{W}_0}
\newcommand{\Wone}{\mathcal{W}_1}

\newcommand{\emptylist}{\varepsilon_{\emptyset}}
\newcommand{\emptyweb}{\Gamma_{\emptyset}}

\newcommand{\Seq}{\mathsf{Seq}}
\newcommand{\Symf}{\mathrm{Sym}_N}
\newcommand{\Foam}{\mathbf{F}_N}

\newcommand{\foamdeg}{\deg_N}
\newcommand{\boundarydeg}{\deg_N}
\newcommand{\smoothing}{\mathcal{S}}
\newcommand{\osum}{\bigoplus}

\DeclareMathOperator{\idfoam}{\mathit{id}}

\title{Two categorifications of the Kekulé relations}

\begin{document}

\author[J.~Grlj]{Jernej Grlj}
\address{Department of Mathematics\\
 University of Southern California \\
  Los Angeles, California 90089, USA}
  \email{grlj@usc.edu}
  
\author[Z.~Lei]{Ziyi Lei}
\address{School of Mathematical Sciences \\
Beijing Normal University \\
Beijing 100875, China}
  \email{202431130003@mail.bnu.edu.cn}

\author[E.~Wagner]{Emmanuel Wagner}
\address{Université Paris Cité, IMJ-PRG \\
Université Paris Sorbonne, UMR 7586 CNRS \\
F-75013 Paris, France}
  \email{emmanuel.wagner@imj-prg.fr}

\begin{abstract}
We give new relations in the kernel of the foam evaluation. As a consequence we obtain categorifications of all the Kelul\'e relations introduced by Morrison. Along the way, we discover an equivalent set of relations for which we provide categorifications as well.
\end{abstract}

\maketitle

\section{Introduction}
In a landmark paper, Cautis, Kamnitzer and Morrison \cite{CKM2014} provided a generators and relations description of the subcategory of representations of quantum $\GL_N$ monoidally generated by exterior powers of the fundamental representation. This is diagrammatically described using so called $\GL_N$-webs. Moreover it is well known that Karoubi completing recovers the full representation category of finite-dimensional representation of quantum $\GL_N$. This question was previously studied by Morrison \cite{Mor2007} in his thesis, where he discovered a set of relations called ``Kekulé relations" which ultimately turned out to be non-necessary to give the presentation.

More recently, the Kekulé relations played an important role in constructing a rotationally invariant basis for the web category in small ranks \cite{Kup1996,SL4Basis}. As mentioned, these relations are redundant in the generators and relations definition of the $\GL_N$ web category. However there is no known diagrammatic derivation for the Kekulé relations outside of the simplest case, benzene relation (for $\GL_4$ webs).

On the other hand, the $\GL_N$-web category was extensively used in the realm of categorification for the Khovanov-Rozansky $\GL_N$-link homology \cite{KR2008, MSV2009}. In fact, a categorification of the $\GL_N$-web category is intrinsically provided by $\GL_N$-foams through the foam evaluation formula of Robert and the third author \cite{RW2020}. This $\GL_N$-foam 2-category is directly linked to many developments in higher representation theory (see \cite{goertz2026monoidal2categoriesfoamevaluation} and the references therein), and has a direct relationship with categorification of quantum groups  through categorified Howe duality \cite{QR2016} and categorification of Hecke algebras \cite{EliasThick} (and more generally Schur algebras \cite{MR2998837}) notably through Soergel bimodules \cite{Elias_2010} and KLR algebras \cite{Khovanov_2009, rouquier20082kacmoody}.

The original aim of this paper was to provide a categorification of the Kekulé relations with two main goals in mind: the first was to have a better understanding of these relations from a categorical perspective in the hope that it would point toward an explicit derivation of them as consequences of the relations from \cite{CKM2014}  
% Cautis, Kamnitzer and Morrison 
which were already categorified in \cite{RW2020}. The second goal was to gain a deeper understanding of the kernel of the foam evaluation.

The second goal is fully achieved -- Proposition \ref{prop:exact-word-projectors} and Theorem \ref{thm:upper-factor-through} contain the most sophisticated relations in the kernel we derived, but Theorem \ref{Null-homotopyTh} and Theorem \ref{thm:composition of slice morphisms} are other prominent examples (and there are plenty of others). All these relations are derived from the foam evaluation formula.

We did not succeed in deriving the Kekulé relations from the classical relations. Instead, we found another set of relations which turned out to be simpler than the Kekulé relations (see Definitions \ref{def:true kekule}, \ref{def:simplified kekule}), while remaining equivalent, see Lemma \ref{lem:true kekule} and Theorem \ref{thm:hocat categorificaiton}.

As the title suggests this paper contains categorifications of the Kekulé relations and fulfills the original aim \textit{twice}. The first categorification of the aforementioned cleaner versions of Kekulé relations is provided by Theorem \ref{thm:hocat categorificaiton}. The categorification therein is given by a collection of long exact sequences in the homotopy category of $\GL_N$-foams. The categorification in the homotopy category was a key step in obtaining (what we believe to be our main result) a categorification in the Karoubi completion of the $\GL_N$-foam category. The latter is achieved by explicitly constructing idempotents in the foam category, recovering all of the Kekulé relations. For the particular choice of webs in Figure \ref{fig:Gammaim}, we strengthen our results further, as stated in the following theorem.

\begin{thm} 
    Let $m, N\geq3$ and let $1\leq i \leq N-2$. Then the web $\Gamma_i^m$ (depicted in Figure \ref{fig:Gammaim}) is an indecomposable object in the $\GL_N$-foam category  if and only if $m\geq N $. Moreover, if $m<N$, then for all $1\leq i\leq N-1$, there exist $\sum_{\max \{0, i-m-1\} \leq c \leq \min \{N - 1 - m, i-1\}}\begin{pmatrix}
        N - m \\ c
    \end{pmatrix}$ mutually orthogonal idempotents  
    of $\Gamma_i^m$ such that the corresponding objects in the Karoubi completion are indecomposable and all are isomorphic up to grading shifts. 
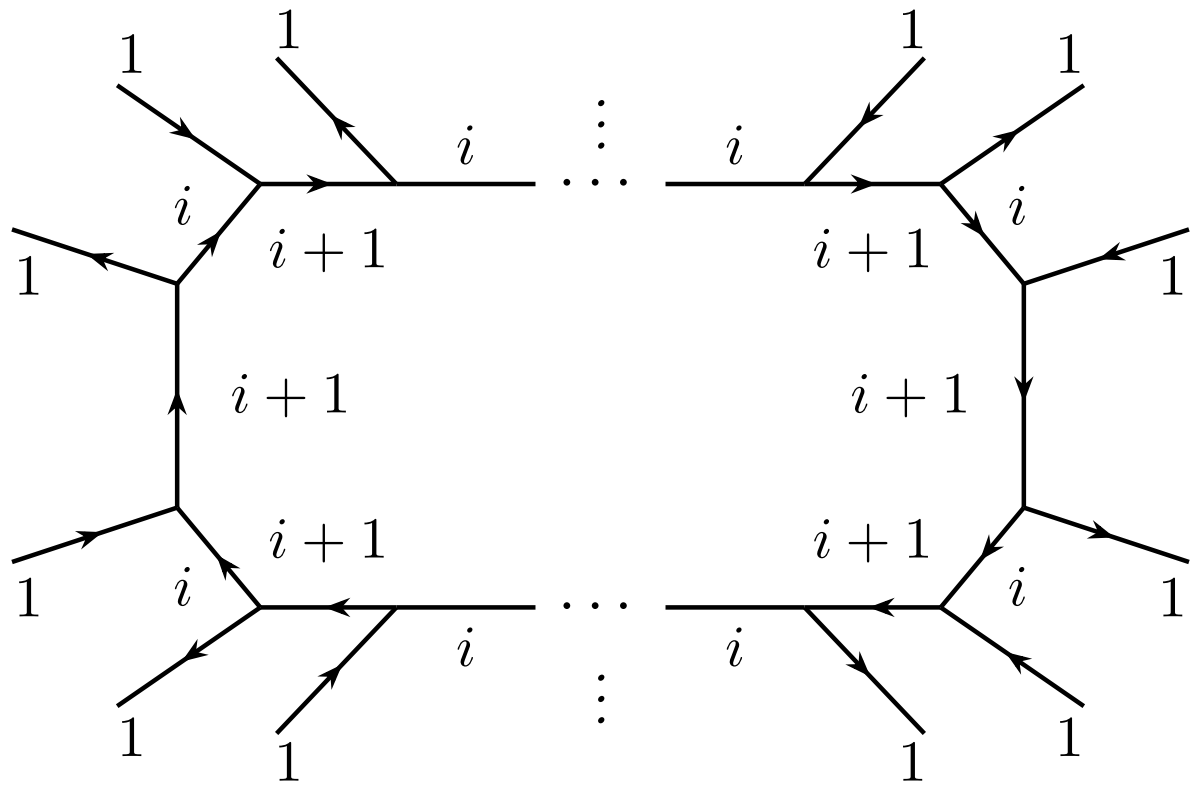
\begin{figure}[htbp!]
    \centering
\tikzset{
  inner cycle/.style={
    draw=black,
    line width=.85pt,
    postaction={decorate},
    decoration={markings,mark=at position .53 with
      {\arrow{Stealth[length=1.7mm,width=1.25mm]}}}
  },
  radial edge/.style={
    draw=black,
    line width=.85pt,
    postaction={decorate},
    decoration={markings,mark=at position .55 with
      {\arrow{Stealth[length=1.7mm,width=1.25mm]}}}
  },
  cycle label/.style={font=\small,inner sep=1pt,fill=white}
}
\begin{tikzpicture}[x=1cm,y=1cm]
  % Normalize the original outward vectors so every free edge
  % has length 1.15 cm, preserving its angle and arrow direction.
  \def\radialLength{1.15}
  \foreach \k/\x/\y/\dx/\dy in {
     0/-1.35/ 1.4/-.95/ 1, 1/-2.25/ 1.4/-1.45/ 1,
     2/-2.8/ .74/-1.7/ .56, 3/-2.8/-.74/-1.7/-.56,
     4/-2.25/-1.4/-1.45/-1, 5/-1.35/-1.4/-.95/-1,
     6/ 1.35/-1.4/ .95/-1, 7/ 2.25/-1.4/ 1.45/-1,
     8/ 2.8/-.74/ 1.7/-.56, 9/ 2.8/ .74/ 1.7/ .56,
    10/ 2.25/ 1.4/ 1.45/ 1,11/ 1.35/ 1.4/ .95/ 1}{
    \coordinate (I\k) at (\x,\y);
    \pgfmathsetmacro{\directionLength}{sqrt(\dx*\dx+\dy*\dy)}
    \coordinate (O\k) at
      ({\x+\radialLength*\dx/\directionLength},
       {\y+\radialLength*\dy/\directionLength});
  }

  % Preserve the directions on the inner polygon.
  \foreach \k in {0,1,2,3,4,6,7,8,9,10}{
    \pgfmathtruncatemacro{\kp}{\k+1}
    \draw[inner cycle] (I\kp) -- (I\k);
  }

  % Alternating outgoing and incoming 1-labeled edges.
  \foreach \k in {0,2,4,6,8,10}
    \draw[radial edge] (I\k) -- (O\k);
  \foreach \k in {1,3,5,7,9,11}
    \draw[radial edge] (O\k) -- (I\k);
  % Place the four top labels above their edges, near the free endpoints.
  \foreach \k in {0,1,10,11}
    \path (I\k) -- (O\k)
      node[pos=.90,above=3pt,cycle label] {$1$};
  \foreach \k in {2,3,4,5,6,7,8,9}
    \path (I\k) -- (O\k)
      node[pos=.90,below=3pt,cycle label] {$1$};

  % Omitted portions of the polygon.
  \foreach \y/\a/\b in {1.4/1.35/.43,-1.4/1.35/.43}{
    \draw[line width=.85pt] (-\a,\y) -- (-\b,\y);
    \draw[line width=.85pt] (\b,\y) -- (\a,\y);
    \node[fill=white,inner sep=0pt] at (0,\y) {$\cdots$};
  }

  % Omitted families of incoming and outgoing edges.
  \node[fill=white,inner sep=0pt] at (0, 1.9) {$\vdots$};
  \node[fill=white,inner sep=0pt] at (0,-1.9) {$\vdots$};

  % Label every visible polygon edge, alternating i+1 and i.
  % Place the four slanted-edge i labels outside the polygon.
  % Keep the i+1 labels inside, clear of the arrows.
  \foreach \x/\y/\lab in {
    -1.80/ .96/{i+1}, -2.76/ 1.26/{i},
    -2.05/ 0/{i+1}, -2.76/-1.26/{i},
    -1.80/-.96/{i+1}, 1.80/-.96/{i+1},
     2.76/-1.26/{i}, 2.05/0/{i+1},
     2.76/ 1.26/{i}, 1.80/.96/{i+1}}
    \node[cycle label] at (\x,\y) {$\lab$};

  % Label both pieces next to each ellipsis outside the polygon.
  \foreach \x in {-.89,.89}{
    \node[cycle label] at (\x, 1.66) {$i$};
    \node[cycle label] at (\x,-1.66) {$i$};
  }
\end{tikzpicture}
    \caption{A particular family of webs. These are objects of further study in Section \ref{sec:indecomposabilis 2m vertex}.}
    \label{fig:Gammaim}
\end{figure}
\end{thm}

The idempotents are constructed in greater generality, but the indecomposability of the corresponding objects remains to be checked for a generic boundary sequence. Nevertheless, the isomorphisms between some of them remain true in full generality (see Equation \eqref{eq:isomorphism for fixed S+i}) and are enough to provide a categorification of all the Kekulé relations.

{\bf Summary of results:} The second and third sections provide all the background on $\GL_N$ webs, $\GL_N$ foams and their evaluations. The fourth section introduces the Kekulé relations themselves and  the equivalent versions we found. The fifth section contains some new relations in the kernel of the foam evaluation and as an application the categorification in the homotopy category. The sixth section contains the main  results, the constructions of the idempotents and a categorification of the Kekulé relations in the Karoubi completion. The seventh section simplifies the idempotents and makes explicit the relationship with \cite{grlj20266valentvertexmathfrakglnweb} resolving a conjecture therein. The eighth section presents our results on the indecomposability of objects.

\subsection*{Acknowledgments} J.G. was partially supported by NSF grants DMS-2200419 and DMS-2601214, Simons Foundation Collaboration Grant on New Structures in Low-Di\-mensional Topology MPS-LDT-00920756 and SFI-MPS-LDT-00014760-01, and Simons Dissertation Fellowship in Mathematics SFI-MPS-SDF-00026689. Z.L. was partially supported by the China Scholarship Council. E.W. was partially supported by the ANR project Caget (ANR-25-CE40-4162).

We are grateful to our coauthors, M. Khovanov, L.-H. Robert, H. Wu, and M. Zhang for the figures and text from \cite{RW2020, grlj20266valentvertexmathfrakglnweb} that we adapted for our use in Sections \ref{sec:gl_nweb} and \ref{sec:gln foam}.

\subsection*{On the use of AI models}
Large Language Models were used to discuss ideas in Section \ref{sec:kar categorification}, they pointed out that our construction in Section \ref{sec:corner} coincides with the concept of corner algebras.
Large Language Models have been used to assist with writing and proofreading, and to assist in proofs of Corollary \ref{prop:two-m extreme coefficients}, Lemma \ref{lem:two-m triangular}, and Propositions \ref{lem:karoubi categorification}, \ref{prop:two-m powers}. 
%to probe and discuss ideas, 
The proofs were checked and then rewritten by the authors. 
% During the revision process, we used Large Language Models to look for typographical and mathematical errors, which we then fixed.
%All errors and are our own.

\section{$\GL_n$ webs}\label{sec:gl_nweb}
The Jones polynomial (also known as the $\mathfrak{sl}_2$ or the $\GL_2$ polynomial) was first introduced by Jones in \cite{Jon1985,Jon1987}. For a given knot the Jones polynomial can be computed by reducing the number of crossings using the skein relations, due to the work of Kauffman \cite{Kau1987}. The natural generalization of skein relations to higher rank Lie algebras was introduced by Kuperberg in \cite{Kup1996} and Murakami, Ohtsuki, and Yamada in \cite{MOY1998}. The rest of the section recalls the definition and some properties of $\GL_N$ webs following the aforementioned work \cite{MOY1998}, the work of Morrison \cite{Mor2007}, and the work of Cautis, Kamnitzer, and Morrison \cite{CKM2014}. 

\begin{definition}
    Let $n,k \in \Z$ be nonnegative integers. Define \emph{the quantum integers}, \emph{quantum factorials} and \emph{the quantum binomial coefficients} as
    \begin{align*}
    [n]_q &:=\frac{q^n-q^{-n}}{q-q^{-1}}, \qquad [-n]_q := - [n]_q \\
    [n]_q! &:= [n]_q \cdot [n-1]_q \cdots [1]_ q \\
            \begin{bmatrix}
    n\\
    k
  \end{bmatrix}_q
  &:=\frac{[n]_q[n-1]_q\cdots [n-k+1]_q}{[k]_q[k-1]_q\cdots[1]_q},
    \end{align*}
where $\begin{bmatrix}
    n\\
    k
  \end{bmatrix}_q$ is understood to be $0$ for $n < k$ or $k < 0$, and $\begin{bmatrix}
    n\\
    0
  \end{bmatrix}_q = 1$.
\end{definition}
\begin{definition}
For some $0 < N \in \Z$, let $\Symf$ denote the ring of symmetric polynomials in $N$ variables. For any $k \in \Z$, $k\geq 0$, define the \emph{elementary symmetric polynomials} and \emph{complete homogeneous symmetric polynomials} in $l$ variables respectively as 
\begin{align*}
    e_k(x_1,  \ldots, x_l) &=\sum_{1\leq i_1<i_2<\cdots<i_k\leq l}x_{i_1}x_{i_2}\cdots x_{i_k}, \\
h_k(x_1,  \ldots, x_l)  &=\sum_{\alpha_i\geq 0,\sum_i \alpha_i=l}x_1^{\alpha_1}x_2^{\alpha_2}\cdots x_l^{\alpha_l}.
\end{align*}
It is understood that $e_0(x_1,  \ldots, x_l) = h_0(x_1,  \ldots, x_l) = 1,$ $e_k(x_1,  \ldots, x_l) = h_k(x_1,  \ldots, x_l)= 0$ for $k<0$, and $e_k(x_1,  \ldots, x_l) = 0$ for $k>l$. 

If $l = N$, then $e_k(x_1,  \ldots, x_N), h_k(x_1,  \ldots, x_N) \in \Symf$. For a set $X$ we will use $e_k(X)$ and $h_k(X)$ to denote the symmetric polynomials in elements of $X$. When the set of variables is clear from the context (see Definition \ref{def:slice morphism} for example) we will simply write $e_k$ or $h_k$, omitting the set of variables. 
\end{definition}

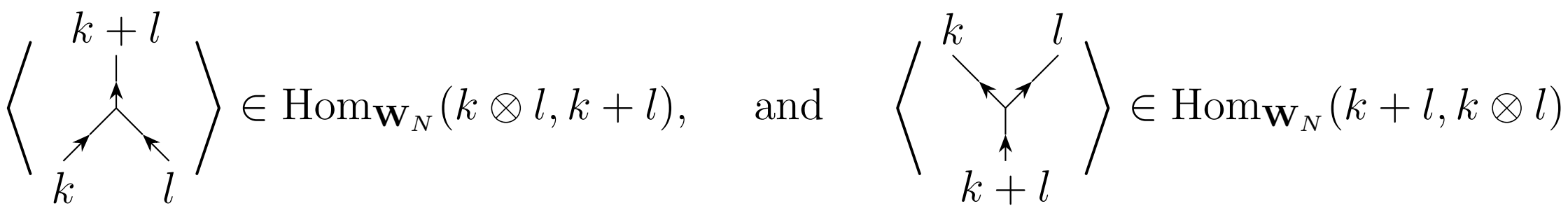
\begin{figure}[htbp]
  \centering
  \begin{equation*}
    \left\langle\begin{tikzpicture}[scale=0.5, baseline=-3]
      \draw[-Stealth] (0,0)--(0,0.5);
      \draw (0,0.5)--(0,1);
      \draw[-Stealth] (-1,-1)--(-0.5,-0.5);
      \draw (-0.5,-0.5)--(0,0);
      \draw[-Stealth] (1,-1)--(0.5,-0.5);
      \draw (0.5,-0.5)--(0,0);
      \node at (0,1.5) {$k+l$};
      \node at (-1,-1.5) {$k$};
      \node at (1,-1.5) {$l$};
    \end{tikzpicture} \right\rangle \in {\Hom}_{\Web}(k\otimes l, k+l),
    \quad \text{ and } \quad
    \left\langle\begin{tikzpicture}[scale=0.5, baseline=-3]
      \draw[-Stealth] (0,-1)--(0,-0.5);
      \draw (0,-0.5)--(0,0);
      \draw[-Stealth] (0,0)--(0.5,0.5);
      \draw (0.5,0.5)--(1,1);
      \draw[-Stealth] (0,0)--(-0.5,0.5);
      \draw (-0.5,0.5)--(-1,1);
      \node at (0,-1.5) {$k+l$};
      \node at (-1,1.5) {$k$};
      \node at (1,1.5) {$l$};
    \end{tikzpicture} \right\rangle \in {\Hom}_{\Web}(k+l, k \otimes l)
  \end{equation*}
  \caption{Trivalent vertices, the defining morphisms of $\Web$.}
  \label{fig:Y}
\end{figure}

\begin{definition} [{\cite[Section 2.2]{CKM2014}}] \label{TypeAWebs}
 Define the $\GL_N$ web category, denoted by 
 $\Web=\webcat_q(\GL_{N})$, to be the pivotal $\mathbb{C}(q)$-linear category monoidally generated by the objects $k^{\pm}$, for $k \in \{0,1,2,...,N\}$,
 where $k^-$ is understood as the dual of $k=k^{+}$. The morphisms of this category are generated by two types of vertices. These are depicted in Figure \ref{fig:Y}, for $l,k \in \mathbb{Z}_{\ge 0}$ and $k+l \leq N$. We impose relations between the morphisms as shown in Figure \ref{relations}.
 \end{definition}

\begin{remark}\label{rmk:web with langle in web category}
    Whenever we will consider a web $W$ as a morphism in $\Web,$ we will write $\langle W \rangle$.
\end{remark}
 
\begin{definition}
    Fix $ k \in \Z$, $k\geq 0$. Let $a_i \in \{1, \ldots, N\}$, and $\epsilon_i = \pm 1$ for $1\leq i \leq k$. Denote by $\varepsilon = a_1^{\epsilon_1}a_2^{\epsilon_2}\cdots a_k^{\epsilon_k}$ \emph{the (boundary) sequence $\varepsilon$} and by $\Seq_{N}$ the set of all such boundary sequences. A boundary sequence is \emph{admissible} when the net flow of labels is balanced, ensuring that the total sum of labels carrying a positive sign equals the total sum of labels carrying a negative sign.
\end{definition}

\begin{notation}
    Given two sequences $\varepsilon_0,\varepsilon_1\in \Seq_{N}$, a $\GL_{N}$ web from $\varepsilon_0$ to $\varepsilon_1$ is a morphism in ${\Hom}_{\Web}(\varepsilon_0,\varepsilon_1)$. Given a web $W_0$ from $\varepsilon_0$ to $\varepsilon_1$ and a web $W_1$ from $\varepsilon_1$ to $\varepsilon_2$, denote by $W_1W_0$ the web from $\varepsilon_0$ to $\varepsilon_2$ by stacking $W_1$ on top of $W_0$. A \emph{closed web} is a web from $\emptylist$ to $\emptylist$, that is, a morphism in ${\Hom}_{\Web}(\emptylist,\emptylist)$. Given a web $W$ from $\varepsilon_0$ to $\varepsilon_1$, denote by $\overline{W}$ the web from $\varepsilon_1$ to $\varepsilon_0$ obtained by taking the reflection of $W$ about a horizontal line and reversing the orientation of each edge.
\end{notation}

\begin{proposition}[{\cite[Theorem 3.3.1]{CKM2014}}]
\label{WebRep}
There is a monoidal functor
\[
\Phi:\webcat_q(\GL_{N}) \rightarrow \Fund (U_q(\GL_{N})),
\]
which is an equivalence of $\mathbb{C}(q)$-linear pivotal categories.
\end{proposition}

$\Fund(U_q(\GL_N))$ is the fundamental representation category of $U_q(\GL_N)$, the full subcategory of the representation category of $U_q(\GL_N)$, monoidally generated by the fundamental representations $\wedge_q^k\mathbb{C}_q^N$ and their duals.

\begin{figure}[htbp]
\newcommand{\RelationsFigureScale}{0.85}
\centering
\scalebox{\RelationsFigureScale}{%
\begin{minipage}{\textwidth}
\begin{equation*}
%
% Circle and associativity relations
%
\left\langle
\begin{tikzpicture}[scale=0.5,baseline=-3]
  \draw (0,0) circle (1);
  \draw[-Stealth] (1,0)--++(0,0.01);
  \node[right] at (1,0) {$k$};
\end{tikzpicture}
\right\rangle
=
\left\langle
\begin{tikzpicture}[scale=0.5,baseline=-3]
  \draw (0,0) circle (1);
  \draw[-Stealth] (1,0)--++(0,-0.01);
  \node[right] at (1,0) {$k$};
\end{tikzpicture}
\right\rangle
=
\begin{bmatrix} n\\ k \end{bmatrix}_q,
\qquad
\left\langle
\begin{tikzpicture}[scale=0.5,baseline=-3]
  \draw (0,-0.5)--(0,-1.5)
        (0,-0.5)--(2,1.5)
        (0,-0.5)--(-2,1.5)
        (1,0.5)--(0,1.5);
  \draw[-Stealth] (0,-1)--++(0,0.01);
  \draw[-Stealth] (0.5,0)--++(0.01,0.01);
  \draw[-Stealth] (-1.5,1)--++(-0.01,0.01);
  \draw[-Stealth] (0.5,1)--++(-0.01,0.01);
  \draw[-Stealth] (1.5,1)--++(0.01,0.01);
  \node at (0,-2) {$i+j+k$};
  \node[right] at (0.5,-0.5) {$j+k$};
  \node at (-2,2) {$i$};
  \node at (0,2) {$j$};
  \node at (2,2) {$k$};
\end{tikzpicture}
\right\rangle
=
\left\langle
\begin{tikzpicture}[scale=0.5,baseline=-3]
  \draw (0,-0.5)--(0,-1.5)
        (0,-0.5)--(2,1.5)
        (0,-0.5)--(-2,1.5)
        (-1,0.5)--(0,1.5);
  \draw[-Stealth] (0,-1)--++(0,0.01);
  \draw[-Stealth] (-0.5,0)--++(-0.01,0.01);
  \draw[-Stealth] (-1.5,1)--++(-0.01,0.01);
  \draw[-Stealth] (-0.5,1)--++(0.01,0.01);
  \draw[-Stealth] (1.5,1)--++(0.01,0.01);
  \node at (0,-2) {$i+j+k$};
  \node[left] at (-0.5,-0.5) {$i+j$};
  \node at (-2,2) {$i$};
  \node at (0,2) {$j$};
  \node at (2,2) {$k$};
\end{tikzpicture}
\right\rangle
\end{equation*}
\begin{gather*}
%
% Digon relations
%
\left\langle
\begin{tikzpicture}[scale=0.5,baseline=-3]
  \draw (0,1)--(0,2) (0,-1)--(0,-2);
  \draw[-Stealth] (0,1.5)--++(0,0.01);
  \draw[-Stealth] (0,-1.5)--++(0,0.01);
  \draw[-Stealth] (-0.5,0)--++(0,0.01);
  \draw[-Stealth] (0.5,0)--++(0,0.01);
  \draw (0,-1) ..controls (-0.5,-0.5) ..(-0.5,0)
        ..controls (-0.5,0.5) ..(0,1);
  \draw (0,-1) ..controls (0.5,-0.5) ..(0.5,0)
        ..controls (0.5,0.5) ..(0,1);
  \node[left] at (-0.5,0) {$i$};
  \node[right] at (0.5,0) {$j$};
  \node at (0,2.5) {$i+j$};
  \node at (0,-2.5) {$i+j$};
\end{tikzpicture}
\right\rangle
=
\begin{bmatrix} i+j\\ i \end{bmatrix}_q
\left\langle
\begin{tikzpicture}[scale=0.5,baseline=-3]
  \draw (0,2)--(0,-2);
  \draw[-Stealth] (0,0)--++(0,0.01);
  \node[right] at (0,0) {$i+j$};
\end{tikzpicture}
\right\rangle,
\qquad
\left\langle
\begin{tikzpicture}[scale=0.5,baseline=-3]
  \draw (0,1)--(0,2) (0,-1)--(0,-2);
  \draw[-Stealth] (0,1.5)--++(0,0.01);
  \draw[-Stealth] (0,-1.5)--++(0,0.01);
  \draw[-Stealth] (-0.5,0)--++(0,0.01);
  \draw[-Stealth] (0.5,0)--++(0,-0.01);
  \draw (0,-1) ..controls (-0.5,-0.5) ..(-0.5,0)
        ..controls (-0.5,0.5) ..(0,1);
  \draw (0,-1) ..controls (0.5,-0.5) ..(0.5,0)
        ..controls (0.5,0.5) ..(0,1);
  \node[left] at (-0.5,0) {$i+j$};
  \node[right] at (0.5,0) {$j$};
  \node at (0,2.5) {$i$};
  \node at (0,-2.5) {$i$};
\end{tikzpicture}
\right\rangle
=
\begin{bmatrix} n-i\\ j \end{bmatrix}_q
\left\langle
\begin{tikzpicture}[scale=0.5,baseline=-3]
  \draw (0,2)--(0,-2);
  \draw[-Stealth] (0,0)--++(0,0.01);
  \node[right] at (0,0) {$i$};
\end{tikzpicture}
\right\rangle
\displaybreak[3]\\[1ex]
%
% First rectangle relation
%
\left\langle
\begin{tikzpicture}[baseline=-3]
  \draw (-1,-1.5)--(-1,1.5)
        (1,-1.5)--(1,1.5)
        (-1,0.5)--(1,0.5)
        (1,-0.5)--(-1,-0.5);
  \draw[-Stealth] (0,0.5)--++(-0.01,0);
  \draw[-Stealth] (0,-0.5)--++(0.01,0);
  \draw[-Stealth] (-1,1)--++(0,0.01);
  \draw[-Stealth] (-1,0)--++(0,-0.01);
  \draw[-Stealth] (-1,-1)--++(0,0.01);
  \draw[-Stealth] (1,1)--++(0,-0.01);
  \draw[-Stealth] (1,0)--++(0,0.01);
  \draw[-Stealth] (1,-1)--++(0,-0.01);
  \node at (0,1) {$k+1$};
  \node at (0,-1) {$k+1$};
  \node[left] at (-1,1) {$1$};
  \node[left] at (-1,0) {$k$};
  \node[left] at (-1,-1) {$1$};
  \node[right] at (1,1) {$k$};
  \node[right] at (1,0) {$1$};
  \node[right] at (1,-1) {$k$};
\end{tikzpicture}
\right\rangle
=
\left\langle
\begin{tikzpicture}[baseline=-3]
  \draw (-1,-1.5)--(-1,1.5) (1,1.5)--(1,-1.5);
  \draw[-Stealth] (-1,0)--++(0,0.01);
  \draw[-Stealth] (1,0)--++(0,-0.01);
  \node[left] at (-1,0) {$1$};
  \node[right] at (1,0) {$k$};
\end{tikzpicture}
\right\rangle
+[n-k-1]_q
\left\langle
\begin{tikzpicture}[baseline=-3]
  \draw (-1,-1.5)--(0,-0.5)
        (1,1.5)--(0,0.5)
        (1,-1.5)--(0,-0.5)
        (-1,1.5)--(0,0.5)
        (0,-0.5)--(0,0.5);
  \draw[-Stealth] (0,0)--++(0,-0.01);
  \draw[-Stealth] (-0.5,1)--++(-0.01,0.01);
  \draw[-Stealth] (-0.5,-1)--++(0.01,0.01);
  \draw[-Stealth] (0.5,1)--++(-0.01,-0.01);
  \draw[-Stealth] (0.5,-1)--++(0.01,-0.01);
  \node[left] at (0,0) {$k-1$};
  \node at (-1,1) {$1$};
  \node at (-1,-1) {$1$};
  \node at (1,1) {$k$};
  \node at (1,-1) {$k$};
\end{tikzpicture}
\right\rangle
\displaybreak[1]\\[1ex]
%
% Second rectangle relation
%
\resizebox{0.98\textwidth}{!}{$\displaystyle
\left\langle
\begin{tikzpicture}[baseline=-3]
  \draw (-1,-1.5)--(-1,1.5)
        (1,-1.5)--(1,1.5)
        (-1,0.5)--(1,0.5)
        (1,-0.5)--(-1,-0.5);
  \draw[-Stealth] (0,0.5)--++(0.01,0);
  \draw[-Stealth] (0,-0.5)--++(-0.01,0);
  \draw[-Stealth] (-1,1)--++(0,0.01);
  \draw[-Stealth] (-1,0)--++(0,0.01);
  \draw[-Stealth] (-1,-1)--++(0,0.01);
  \draw[-Stealth] (1,1)--++(0,0.01);
  \draw[-Stealth] (1,0)--++(0,0.01);
  \draw[-Stealth] (1,-1)--++(0,0.01);
  \node at (0,1) {$j$};
  \node at (0,-1) {$j+k-1$};
  \node[left] at (-1,1) {$k$};
  \node[left] at (-1,0) {$j+k$};
  \node[left] at (-1,-1) {$1$};
  \node[right] at (1,1) {$i$};
  \node[right] at (1,0) {$i-j$};
  \node at (0.7,-1.8) {$i+k-1$};
\end{tikzpicture}
\right\rangle
=
\begin{bmatrix} i-1\\ j \end{bmatrix}_q
\left\langle
\begin{tikzpicture}[baseline=-3]
  \draw (-1,-1.5)--(-1,1.5)
        (1,1.5)--(1,-1.5)
        (-1,0)--(1,0);
  \draw[-Stealth] (-1,1)--++(0,0.01);
  \draw[-Stealth] (-1,-1)--++(0,0.01);
  \draw[-Stealth] (1,1)--++(0,0.01);
  \draw[-Stealth] (1,-1)--++(0,0.01);
  \draw[-Stealth] (0,0)--++(-0.01,0);
  \node at (0,0.5) {$k-1$};
  \node[left] at (-1,1) {$k$};
  \node[left] at (-1,-1) {$1$};
  \node[right] at (1,1) {$i$};
  \node at (0.7,-1.8) {$i+k-1$};
\end{tikzpicture}
\right\rangle
+
\begin{bmatrix} i-1\\ j-1 \end{bmatrix}_q
\left\langle
\begin{tikzpicture}[baseline=-3]
  \draw (-1,-1.5)--(0,-0.5)
        (1,1.5)--(0,0.5)
        (1,-1.5)--(0,-0.5)
        (-1,1.5)--(0,0.5)
        (0,-0.5)--(0,0.5);
  \draw[-Stealth] (0,0)--++(0,0.01);
  \draw[-Stealth] (-0.5,1)--++(-0.01,0.01);
  \draw[-Stealth] (-0.5,-1)--++(0.01,0.01);
  \draw[-Stealth] (0.5,1)--++(0.01,0.01);
  \draw[-Stealth] (0.5,-1)--++(-0.01,0.01);
  \node[left] at (0,0) {$i+k$};
  \node at (-1,1) {$k$};
  \node at (-1,-1) {$1$};
  \node at (1,1) {$i$};
  \node at (0.7,-1.8) {$i+k-1$};
\end{tikzpicture}
\right\rangle
$}
\displaybreak[1]\\[1ex]
%
% Third rectangle relation
%
\resizebox{0.98\textwidth}{!}{$\displaystyle
\left\langle
\begin{tikzpicture}[scale=1.2,baseline=-3]
  \draw (-1,-1.5)--(-1,1.5)
        (1,-1.5)--(1,1.5)
        (-1,0.5)--(1,0.5)
        (1,-0.5)--(-1,-0.5);
  \draw[-Stealth] (0,0.5)--++(0.01,0);
  \draw[-Stealth] (0,-0.5)--++(-0.01,0);
  \draw[-Stealth] (-1,1)--++(0,0.01);
  \draw[-Stealth] (-1,0)--++(0,0.01);
  \draw[-Stealth] (-1,-1)--++(0,0.01);
  \draw[-Stealth] (1,1)--++(0,0.01);
  \draw[-Stealth] (1,0)--++(0,0.01);
  \draw[-Stealth] (1,-1)--++(0,0.01);
  \node at (0,0.8) {$-i+j+k$};
  \node at (0,-0.8) {$k$};
  \node[left] at (-1,1) {$i$};
  \node[left] at (-1,0) {$j+k$};
  \node[left] at (-1,-1) {$j$};
  \node[right] at (1,1) {$j+m$};
  \node[right] at (1,0) {$i-k+m$};
  \node[right] at (1,-1) {$i+m$};
\end{tikzpicture}
\right\rangle
=
\sum_{l=\max(0,i-j)}^i
\begin{bmatrix} m\\ k-l \end{bmatrix}_q
\left\langle
\begin{tikzpicture}[scale=1.2,baseline=-3]
  \draw (-1,-1.5)--(-1,1.5)
        (1,-1.5)--(1,1.5)
        (-1,0.5)--(1,0.5)
        (1,-0.5)--(-1,-0.5);
  \draw[-Stealth] (0,0.5)--++(-0.01,0);
  \draw[-Stealth] (0,-0.5)--++(0.01,0);
  \draw[-Stealth] (-1,1)--++(0,0.01);
  \draw[-Stealth] (-1,0)--++(0,0.01);
  \draw[-Stealth] (-1,-1)--++(0,0.01);
  \draw[-Stealth] (1,1)--++(0,0.01);
  \draw[-Stealth] (1,0)--++(0,0.01);
  \draw[-Stealth] (1,-1)--++(0,0.01);
  \node at (0,0.8) {$l$};
  \node at (0,-0.8) {$-i+j+l$};
  \node[left] at (-1,1) {$i$};
  \node[left] at (-1,0) {$i-l$};
  \node[left] at (-1,-1) {$j$};
  \node[right] at (1,1) {$j+m$};
  \node[right] at (1,0) {$j+m+l$};
  \node[right] at (1,-1) {$i+m$};
\end{tikzpicture}
\right\rangle
$}
\end{gather*}
\end{minipage}}
\caption{Relations between the morphisms in $\Web$.}
\label{relations}
\end{figure}

\section{$\GL_n$ foam category and foam evaluation}\label{sec:gln foam}
Khovanov's categorification of Jones polynomial \cite{Kho2000} pioneered the study of knot homology theories. These include the categorification of the $\GL_3$ knot polynomial in \cite{Kho2004}, where the notion of foams first appeared. The generalization to the $\GL_N$ case is due to Queffelec, and Rose \cite{QR2016}, and Robert and the third author \cite{RW2020}. We will introduce $\GL_N$ foams following the latter. This leads to a $1$-category whose objects are $\GL_N$ webs and whose morphisms are foams. 

\begin{notation}
     We use ``$\equiv$'' to denote congruences modulo $2$ as in \cite{RW2020}.
\end{notation}

\begin{definition}
    Let $P$ be a proposition (in the sense of mathematical logic). Denote
    \begin{align*}
        \mathbbm{1}(P) &:= \begin{cases} 1; \text{ when } P \text{ is true} \\
        0;\text{ when } P \text{ is false.}\end{cases}
    \end{align*}
\end{definition}

\begin{definition}
  A \emph{$\GL_N$ foam} $F$ is a collection of \emph{facets} $\mathcal{F}(F)=(\Sigma_i)_{i\in I}$, that is a finite collection of oriented connected surfaces $(\Sigma_i)_{i\in I}$ with boundary, together with:
  \begin{itemize}
    \item A labeling $l:(\Sigma_i)_{i\in I}\to\{0,\cdots,N\}$,
    \item A gluing recipe of the facets along their boundaries such that when glued together using the recipe there are three possible local models: a neighborhood of a point is homeomorphic to a surface, a tripod times an interval or the cone of the $1$-skeleton of the tetrahedron as depicted by Figure~\ref{fig:localmodels}, which can be seen as the neighborhood on facets, \emph{bindings}, and \emph{singular vertices}. The letter appearing on each facet indicates the label of that facet. Each binding carries:
    \begin{itemize}
    \item an orientation which agrees with the orientations of the facets with
      labels $a$ and $b$ and disagrees with the orientation of the facet with label
      $a+b$.
    \item a cyclic ordering of the three facets around it. When a foam
      is embedded in $\mathbb{R}^3$, we require this cyclic ordering to agree with 
      the left-hand rule with respect to its
      orientation (the dotted circle in the middle indicates that the
      orientation of the binding points to the reader):
      
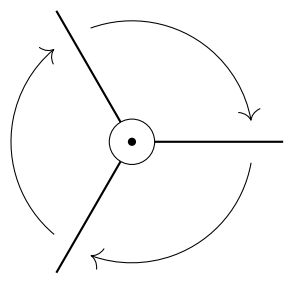
\[
        \begin{tikzpicture}[xscale=1]
          \begin{scope}
            \draw (0,0) -- +(0:1);
            \draw (0,0) -- +(120:1);
            \draw (0,0) -- +(240:1);
            \filldraw[fill= white, draw=black, very thin] (0,0) circle (0.15cm);
            \filldraw[fill = black] (0,0) circle (0.02cm);
            \draw[very thin,->] (-10:0.8cm) arc (-10:-110 :0.8); 
            \draw[very thin,->] (230:0.8cm) arc (230: 130:0.8); 
            \draw[very thin,->] (110:0.8cm) arc (110:10 :0.8); 
           \end{scope}
        \end{tikzpicture}
      \]
    \end{itemize}
    The (left-hand) cyclic orderings of the different bindings adjacent to a
    singular point should be compatible. This means that a
    neighborhood of the singular point is embeddable in $\mathbb{R}^3$ in a
    way that respects the left-hand rule for the four bindings adjacent to this singular vertex.
  \end{itemize}
  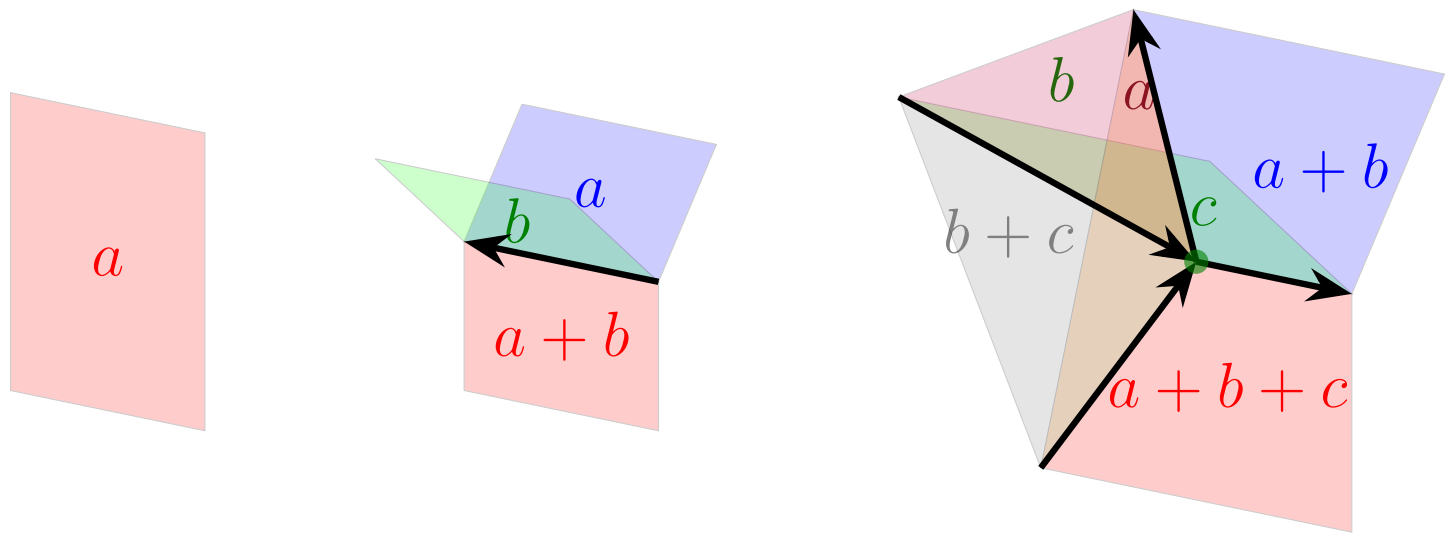
\begin{figure}[htbp]
    \[
      \begin{tikzpicture}
        \begin{scope}
          \tdplotsetmaincoords{80}{140}
            \begin{scope}[tdplot_main_coords]
              \filldraw [very thin, fill=red, opacity = 0.2] (0,-1, -1) -- (0,1,-1) -- (0,1,1) -- (0,-1,1) -- (0,-1,-1);
              \node[sloped, red] at (0,0,0) {$a$};
            \end{scope}
          \begin{scope}[xshift = 3cm, tdplot_main_coords]
            \begin{scope}
              \filldraw [very thin, fill =red, opacity = 0.2] (0,-1, -1) -- (0,1,-1) -- (0,1,0) -- (0,-1,0) -- (0,-1,-1);
                  \node[sloped, red] at (0,0,-0.5) {$a+b$};
              \end{scope}
            \begin{scope}[rotate around y = 150]
              \filldraw [very thin, fill =blue, opacity = 0.2] (0,-1, -1) -- (0,1,-1) -- (0,1,0) -- (0,-1,0) -- (0,-1,-1);
                  \node[sloped, blue] at (0,0,-0.5) {$a$};
            \end{scope}
            \begin{scope}[rotate around y =-130]
              \filldraw [very thin, fill =green, opacity = 0.2] (0,-1, -1) -- (0,1,-1) -- (0,1,0) -- (0,-1,0) -- (0,-1,-1);
                  \node[sloped, green!50!black] at (0,0,-0.5) {$b$};
                  \draw[very thick, -Stealth] (0,1,0) -- (0,-1,0);
            \end{scope}
            \end{scope}
          \begin{scope}[scale = 1.6, xshift = 4.5cm, tdplot_main_coords]
            \begin{scope}
              \filldraw [very thin, fill =red, opacity = 0.2] (0,-1, -1) -- (0,1,-1) -- (0,1,0) -- (0,0,0) -- (0,-1,-1);
              \coordinate (a) at (0, -1, -1);
                  \node[sloped, red] at (0,0.2,-0.5) {$a+b+c$};
              \end{scope}
            \begin{scope}[rotate around y = 150]
              \filldraw [very thin, fill= blue, opacity = 0.2] (0,-1, -1) -- (0,1,-1) -- (0,1,0) -- (0,0,0) -- (0,-1,-1);
              \coordinate (b) at (0, -1, -1);
                  \node[sloped, blue] at (0,0.5,-0.5) {$a+b$};
            \end{scope}
            \begin{scope}[rotate around y =-130]
              \filldraw [very thin, fill= green, opacity = 0.2] (0,-1, -1) -- (0,1,-1) -- (0,1,0) -- (0,0,0) -- (0,-1,-1);
              \coordinate (c) at (0, -1, -1);
                  \node[sloped, green!50!black] at (0,0.5,-0.5) {$c$};
                  \node[sloped, purple!50!black] at (0,1,-1.5) {$a$};
                  \node[sloped, green!50!black] at (0,+0.5,-1.5) {$b$};
                  \node[sloped, gray] at (0,-1.3,0.1) {$b+c$};
            \end{scope}
            \filldraw[very thin, fill= orange, opacity = 0.2] (a) -- (b) -- (0,0,0) -- (a);
            \filldraw[very thin, fill= gray, opacity = 0.2] (a) -- (c) -- (0,0,0) -- (a);
            \filldraw[very thin, fill= purple, opacity = 0.2] (b) -- (c) -- (0,0,0) -- (b);
            \draw[very thick, -Stealth] (a) -- (0,0,0);
            \draw[very thick, Stealth-] (b) -- (0,0,0);
            \draw[very thick, -Stealth] (c) -- (0,0,0);
            \draw[very thick, Stealth-] (0,1, 0) -- (0,0,0);
            \fill[green!50!black, opacity = 0.6] (0,0,0) circle (0.5mm);
          \end{scope}
          \end{scope}          
      \end{tikzpicture}
       \]
       \caption{Local models of foams}
       \label{fig:localmodels}
  \end{figure}
\end{definition}

\begin{definition}
  A \emph{pigment} is an element of $\mathbb{P}= \{1,\dots,N\}$. A \emph{coloring} of a foam $F$ is a map $c: \mathcal{F}(F) \to 2^\mathbb{P}$, such that:
  \begin{itemize}
  \item For each facet $f$, the number of elements $\# c(f)$ of $c(f)$ is equal to the label $l(f)$.
  \item For each binding joining a facet $f_1$ with label $a$, a facet $f_2$ with label $b$, and a facet $f_3$ with label $a+b$, we have $c(f_1) \cup c(f_2) = c(f_3)$. This condition is called the \emph{flow condition}.
  \end{itemize}
A \emph{colored foam} is a foam together with a coloring.
\end{definition}
\begin{definition}
  A \emph{decoration} of a foam $F$ is a map $f\mapsto P_f$ which associates with each facet $f\in F$ a homogeneous symmetric polynomial $P_f$ in $l(f)$ many variables from the set $\{x_1,\cdots,x_N\}$. A foam together with a decoration is called a decorated foam.   Note that every \(x_i\), \(1\le i\le N\), has degree \(2\).
\end{definition}
\begin{definition}
  Let $(F,c)$ be a colored foam.
  \begin{enumerate}
    \item If $i$ is an element of $\mathbb{P}$, the union (with the identification coming from the gluing procedure) of all the facets which contain the pigment $i$ in their colors is a surface. It is called the \emph{monochrome surface of $(F,c)$ associated with $i$} and it is denoted $F_i(c)$. The restriction we imposed on the orientations of facets ensures that $F_i(c)$ is oriented.
    \item  If $i<j$ are two distinct elements of $\mathbb{P}$, the symmetric difference of $F_i(c)$ and $F_j(c)$ is a surface. It is called the \emph{bichrome surface of $(F,c)$ associated with $i,j$} and it is denoted $F_{ij}(c)$ or $F_{ji}(c)$. The restriction we imposed on the orientations of facets ensures that $F_{ij}(c)$ can be oriented by taking the orientation of the facets containing $i$ and the reverse orientation on the facets containing $j$.
    \item If $i<j$ are two distinct elements of $\mathbb{P}$, We consider a binding joining the facets $f_1$, $f_2$ and $f_3$. Suppose that $i$ is in $c(f_1)$, $j$ is in $c(f_2)$ and $\{i,j\}$ is  in $c(f_3)$. We say that the binding is \emph{positive with respect to $(i,j)$} if the cyclic order otherwise. The set $F_{i}(c) \cap F_{j}(c) \cap F_{ij}(c)$ is a collection of disjoint circles. Each of these circles is a union of bindings, for every circle the bindings are either all positive or all negative with respect to $(i,j)$. We denote by $\theta^+_{ij}(F,c)$ (resp. $\theta^-_{ij}(F,c)$) the number of positive (resp. negative) circles with respect to $(i,j)$. Note that $\theta_{ji}, \theta^+_{ji},$ and $\theta^-_{ji}$ are understood to denote the same set and counts as $\theta_{ij}, \theta^+_{ij},$ and $ \theta^-_{ij}$ respectively.
  \end{enumerate}
\end{definition}
\begin{definition}
  If $(F,c)$ is a colored decorated foam, we define:
    \begin{align*}
      s(F,c) &= \sum_{i=1}^N {i\chi(F_i(c))/2}  + \sum_{1\leq i < j \leq N} \theta^+_{ij}(F,c), \\
      P(F,c) &= \prod_{f \textrm{ facet of $F$}} P_f(c(f)),\\%X_{c(f)}), \\
      Q(F,c) &= \prod_{1\leq i < j \leq N} (x_i-x_j)^{\frac{\chi(F_{ij}(c))}{2}} \\
      \langle\,F,c\,\rangle &= (-1)^{s(F,c)} \frac{P(F,c)}{Q(F,c)}.
    \end{align*}
    For a decorated foam $F$, define the \emph{evaluation} of the foam $F$ by:
    \[ \langle\,F\,\rangle := \sum_{c \textrm{ coloring of $F$}} \langle\,F,c\,\rangle. \]
\end{definition}

The foam evaluation of a closed foam results in a symmetric polynomial $\langle F \rangle \in \mathrm{Sym}_N (x_1, \ldots, x_N)$ \cite[Proposition 2.18]{RW2020}. In order to talk about the Karoubi envelope of webs with boundary, we need to introduce the notion of foams with corners. This exposition follows \cite[Section 2.4]{grlj20266valentvertexmathfrakglnweb}.

\subsection{\texorpdfstring{$\GL_N$}{GL\_N} foams with boundary}
\label{foam with boundary}

A closed foam comes with an embedding $\iota : F \longrightarrow
\mathbb{R}^3$. For $r \in \R$, let $T_r$ be the plane
$\{(x,y,z)\mid z=r\}$. We say that $T_r$ intersects $\iota(F)$
\emph{generically} if the following conditions hold:
\begin{itemize}
    \item $\iota(F) \cap T_r$ is a closed $\GL_N$ web $\Gamma$ in
    $T_r$, where the orientations of the edges of $\Gamma$ are induced by
    the orientations of the facets of $\iota(F)$.
    \item There exists $\epsilon>0$ such that
    \[
    \iota(F) \cap \{(x,y,z)\mid x,y\in \mathbb{R},\
      z \in (r-\epsilon,r+\epsilon)\}
    \]
    is homeomorphic to $\Gamma\times(0,1)$.
\end{itemize}

In other words, $T_r$ intersects $\iota(F)$ generically if
$T_r\cap \iota(F)$ does not contain any singular points of $\iota(F)$ and
$T_r$ intersects the binding lines and the facets of $\iota(F)$
transversely.

\begin{definition}\label{def:foam with closed boundary}
Given an embedding of a closed foam
$\iota : F \longrightarrow \mathbb{R}^3$, define a foam with boundary
$(\Gamma_0,\Gamma_1)$ as the intersection of $\iota(F)$ with
$\mathbb{R}^2 \times [0,1]$ such that $T_0$ and $T_1$ intersect
$\iota(F)$ generically and
$T_i \cap \iota(F)=\Gamma_i$ are closed webs, for $i=0,1$.
\end{definition}

We refer to a foam with boundary $(\Gamma_0,\Gamma_1)$ as a
\emph{$(\Gamma_0,\Gamma_1)$-foam} or simply as \emph{foam} when the boundaries are clear from the context.

\begin{remark}
A foam with boundary $(\Gamma_0,\Gamma_1)$ is thought of as a singular
cobordism from $\Gamma_0$ to $\Gamma_1$. We choose the orientations of the
edges of $\Gamma_1$ to agree with the edge orientations induced from the
orientations of the facets of $\iota(F)$, and the orientations of the edges
of $\Gamma_0$ to be opposite to the edge orientations induced by the
orientations of the facets of $\iota(F)$.

Consider a foam $U$ with boundary $(\Gamma_0,\Gamma_1)$ and a foam $V$ with
boundary $(\Gamma_1,\Gamma_2)$. We define the composition of $V$ and $U$,
denoted (as multiplication) by $V\cdot U$ or $V U$, as the foam with boundary $(\Gamma_0,\Gamma_2)$
obtained by stacking $V$ on top of $U$ and rescaling to fit in
$\mathbb{R}^2\times[0,1]$.
\end{remark}

\begin{definition}[{\cite[Definition 2.3]{RW2020}}]
\label{dfn:degreefoam}
The \emph{degree} of a closed $\GL_N$ foam $F$ is defined to be
\[
\foamdeg(F):=
-\sum_{f\,\mathrel{\mathrm{facet}}}d(f)
+\sum_{\substack{i\,\mathrel{\mathrm{interval}}\\
                   \mathrm{binding}}}d(i)
-\sum_{\substack{p\,\mathrel{\mathrm{singular}}\\
                   \mathrm{point}}}d(p)
+\sum_{\substack{s\,\mathrel{\mathrm{labeling}}\\
                   \mathrm{of\ a\ dot}}}\deg(s),
\]
with terms defined as follows.
\begin{itemize}
    \item For each facet $f$ with label $a$, set
    $d(f)=a(N-a)\chi(f)$, where $\chi(f)$ is the Euler characteristic of
    $f$.
    \item For each interval (i.e.\ not a circle) binding $i$ with adjacent
    facets labeled $(a,b,a+b)$, set
    $d(i)=ab+(a+b)(N-a-b)$.
    \item For each singular vertex $p$ with adjacent facets labeled
    $(a,b,c,a+b,b+c,a+b+c)$, set
    $d(p)=ab+bc+cw+wa+ac+bw$, where $w=N-a-b-c$.
    \item Finally, $\deg(s)$ is the degree of the homogeneous symmetric
    polynomial $s$.
\end{itemize}
\end{definition}

\begin{remark}
The same formula is used to define the degree of a $\GL_N$ foam with
boundary; see the discussion in \cite[Section 4]{RW2020}.
\end{remark}

\begin{definition}\label{def:foam between closed webs}
Given a closed web $\Gamma$, the state space of $\Gamma$, denoted by
$\mathcal{F}(\Gamma)$, is a $\Symf$-module generated by foams with boundary
$(\emptyweb,\Gamma)$ modulo the foam evaluation formula.

In other words, $\mathcal{F}(\Gamma)$ is a $\Symf$-module generated by
symbols $\mathcal{F}(U)$ for all foams $U$ with boundary
$(\emptyweb,\Gamma)$, where a relation
$\sum_i a_i\mathcal{F}(U_i)=0$ holds in $\mathcal{F}(\Gamma)$ for
$a_i\in\Symf$ if and only if
\[
\sum_i a_i\langle V\cdot U_i\rangle=0
\]
for any foam $V$ with boundary $(\Gamma,\emptyweb)$. Note that
$\langle V\cdot U_i\rangle\in\Symf$ is the evaluation of the closed foam
$V\cdot U_i$.

The state space $\mathcal{F}(\Gamma)$ is naturally a graded $\Symf$-module.
The degree of an element $\mathcal{F}(U)$ in $\mathcal{F}(\Gamma)$ is given
by the degree of the foam $U$, and for any $a \in\Symf$, 
$\deg(a\mathcal{F}(U))=\deg(a)+\deg(\mathcal{F}(U))$.

There is a direct sum decomposition of the graded abelian group
\[
\mathcal{F}(\Gamma)=\bigoplus_k\mathcal{F}(\Gamma)_k,
\]
where $\mathcal{F}(\Gamma)_k$ is the subspace of elements in
$\mathcal{F}(\Gamma)$ of degree $k$, together with the zero element.

We denote by $q^i\mathcal{F}(\Gamma)$ the graded $\Symf$-module obtained by
raising the grading of $\mathcal{F}(\Gamma)$ by $i$:
\[
q^i\mathcal{F}(\Gamma)=\bigoplus_k(q^i\mathcal{F}(\Gamma))_k,
\qquad
(q^i\mathcal{F}(\Gamma))_k=\mathcal{F}(\Gamma)_{k-i}.
\]
\end{definition}

\begin{proposition}
A foam $V$ with boundary $(\Gamma_0,\Gamma_1)$ induces a
grading-preserving homomorphism
\begin{align*}
\mathcal{F}(V):\mathcal{F}(\Gamma_0)
  &\longrightarrow q^{-\foamdeg(V)}\mathcal{F}(\Gamma_1),\\
\mathcal{F}(U)&\longmapsto\mathcal{F}(V\cdot U).
\end{align*}
\end{proposition}

Now we consider foams between $\GL_N$ webs with a fixed boundary instead of foams between closed $\GL_N$ webs.

\begin{definition}\label{def:foamfromto}
Consider two webs $W_0$ and $W_1$ with the same boundary $\varepsilon$.
Define a foam $F$ from $W_0$ to $W_1$ as a foam $F_\cup$ with boundary
$(\emptyweb,\overline{W_0}W_1)$ cut along a tubular neighborhood of
$\varepsilon$ relative to $F_\cup$. We will consider the boundary of $F$ in
a natural way as
\[
\partial F=(\overline{W_0}\times\{0\})\cup(W_1\times\{1\})
\cup\left(\bigsqcup_{p\in\varepsilon}p\times[0,1]\right).
\]
The orientation of $W_1$ and $\varepsilon\times[0,1]$ agrees with the
orientation of the facets of $F$, and the orientation of $W_0$ is opposite
to the orientation of the facets of $F$. A coloring of $F$ is a coloring of
$F_\cup$ that respects the cutting we performed.

When $\varepsilon\neq\emptylist$, we call $F$ a \emph{foam with corners}.
\end{definition}

We will occasionally refer to a foam from $W_0$ to $W_1$ as a
$(W_0,W_1)$-foam or simply as foam when the boundary webs are clear from the context, just as in the case of foams with boundary.

We present a way of concatenating foams with corners in the next definition, also allowing linear combinations of foams.

\begin{definition}\label{Foam-Compo-Conca}
Given webs $W_0,W_1,$ and $W_2$ with the same boundary $\varepsilon$,
consider a foam $M_0$ from $W_0$ to $W_1$ and a foam $M_1$ from $W_1$ to
$W_2$. Define the composition of $M_1$ with $M_0$, denoted as multiplication by
$M_1\cdot M_0$ or $M_1 M_0$, as the foam from $W_0$ to $W_2$ obtained by stacking $M_1$
on top of $M_0$. Consider a finite $\Symf$-linear combination
$f_0=\sum_i u_i\cdot\leftindex^i M_0$ of foams $\leftindex^i M_0$ from
$W_0$ to $W_1$ and a finite $\Symf$-linear combination
$f_1=\sum_j v_j\cdot\leftindex^j M_1$ of foams $\leftindex^j M_1$ from
$W_1$ to $W_2$; define the composition of $f_1$ with $f_0$ as
\[
f_1\cdot f_0=\sum_{i,j}v_ju_i\cdot
\leftindex^j M_1\cdot\leftindex^i M_0.
\]
\end{definition}

We will now adapt Definition~\ref{dfn:degreefoam}
to foams with corners. The original definition is due to
\cite[Definition 3.3]{QR2016}.

\begin{lemma}\label{lem:euler-char-boundary-conditions}
Let
$\varepsilon=a_1^{\epsilon_1}a_2^{\epsilon_2}\cdots a_k^{\epsilon_k}
\in\Seq_N$. Let $W$ be a web with boundary $\varepsilon$ and choose any
coloring $c$ for the identity foam $F:=\idfoam_W$. Then
\[
\sum_{1\leq i<j\leq N}\chi(F_{ij})
=\frac12\sum_{l=1}^k a_l(N-a_l).
\]
In particular, this quantity depends only on the boundary sequence
$\varepsilon$.
\end{lemma}

Lemma~\ref{lem:euler-char-boundary-conditions} justifies the following
definition and notation.

\begin{definition}\label{def:deg of web with boundary}
 For an admissible boundary $\varepsilon =a_1^{\epsilon_1}a_2^{\epsilon_2}\cdots a_k^{\epsilon_k}$ sequence let
\[
\boundarydeg(\varepsilon):=\frac12\sum_{l=1}^k a_l(N-a_l).
\]
\end{definition}

This quantity appears in \cite[Section 2]{Robert2015} and corresponds to the
dimension of complex Grassmannians. It will serve as a correction term in
the degree of a foam with corners.

\begin{definition}\label{cor:degree foam with corner alt}
Let $W_0,W_1$ be webs with boundary $\varepsilon$. Let $F$ be a foam from
$W_0$ to $W_1$. Then
\[
\foamdeg(F)
=-\sum_{1\leq i<j\leq N}\chi(F_{ij}(c))
 +\sum_{f\mathrel{\mathrm{facet\ of\ }}F}\deg(P_f)
 +\boundarydeg(\varepsilon).
\]
\end{definition}

This definition ensures that the identity foam on any web $W$ with
boundary $\varepsilon$ has degree $0$.

\begin{remark}
\label{rem:bending-up-foamdeg}
Let $W, W'$ be webs with boundary $\varepsilon$, and let $F$ be a $(W,W')$-foam.
Let $F_\cup$ be the foam with boundary in $\Hom_{\Foam}(\emptyweb, \overline{W} W')$ 
associated to $F$ (see Definition \ref{def:foamfromto}).
Let $c$ be a coloring of $F$.
Lemma \ref{lem:euler-char-boundary-conditions} shows that
\[
    \foamdeg(F_\cup)= -\sum_{1 \leq i < j \leq N} \chi(F_{ij}(c)) + \sum_{f \textrm{ facet of $F$}} \deg(P_f).
\]
\end{remark}

\begin{remark}\label{rmk:degree with decomposition}
    The degree of a foam can alternatively be computed by decomposing the foam into a collection of local models. See \cite[Definitions 2.12 and 2.13]{ETW2018}.
\end{remark}

\begin{definition}
Fix an admissible boundary sequence $\varepsilon\in\Seq_N$. The
\emph{$\mathfrak{gl}_N$ foam category with boundary $\varepsilon$}, denoted
$\Foam^\varepsilon$, is defined as follows. Its objects are the webs with
boundary $\varepsilon$. Given two such webs $W_0$ and $W_1$, a morphism from
$W_0$ to $W_1$ is a finite formal $\Symf$-linear combination of foams from
$W_0$ to $W_1$,
\[
f=\sum_i k_iM_i,
\]
where each $k_i\in\Symf$ and each $M_i$ is a foam from $W_0$ to $W_1$.
These morphism spaces are taken modulo the following evaluation relation:
for foams $M_i$ from $W_0$ to $W_1$ and coefficients $a_i\in\Symf$, we set
$\sum_i a_iM_i=0$ if and only if
$\sum_i a_i\langle VM_i\rangle=0$ for every foam $V$ with boundary
$(\overline{W_0}W_1,\emptyweb)$, where
$\langle VM_i\rangle\in\Symf$ denotes the evaluation of the closed foam
$VM_i$.
\end{definition}

\begin{remark}
When the boundary sequence $\varepsilon$ is clear from the context or when it does not play an important role, we allow ourselves to simply write $\Foam$ instead of $\Foam^\varepsilon$. 

In contrast to Remark \ref{rmk:web with langle in web category}, we will write simply $W$ (without $\langle\_\rangle$) for a web when considered as an object in $\Foam^\varepsilon$.
\end{remark}

The category $\Foam^\varepsilon$ is graded. The degree of an individual term
$k_iM_i$ is
\[
\deg(k_iM_i):=\deg(k_i)+\foamdeg(M_i),
\]
and for each $d\in\mathbb{Z}$ the homogeneous degree-$d$ morphism space
$\Hom_{\Foam^\varepsilon}^d(W_0,W_1)$ is spanned by those finite sums
$\sum_i k_iM_i$ in which every nonzero term has total degree $d$. Thus
\[
\Hom_{\Foam^\varepsilon}(W_0,W_1)
=\bigoplus_{d\in\mathbb{Z}}\Hom_{\Foam^\varepsilon}^d(W_0,W_1)
\]
is a graded $\Symf$-module.

Note that in comparison to Definition~\ref{def:foam between closed webs}, we
omit the analogue of $\mathcal{F}$ in order to make the text in the next
sections cleaner. We will use $\mathcal{E}$ for objects in the Karoubi
envelope, which we define now.

\begin{definition}\label{def:Karoubi of cornered foams}
Let $\Kar(\Foam^\varepsilon)$ denote the \emph{Karoubi
completion} (or \emph{Karoubi envelope}) of $\Foam^\varepsilon$. An object
of $\Kar(\Foam^\varepsilon)$ is a pair $(W,P)$, where $W$ is
a web with boundary $\varepsilon$ and
\[
P\in
\Hom^0_{\Foam^\varepsilon}(W,W)
\]
is a homogeneous degree-zero idempotent, that is,
$P\cdot P=P$.
We denote this object by $\mathcal{E}(P)$. In particular, if
$P=\idfoam_W$, we write
\[
\mathcal{E}(W):=\mathcal{E}(\idfoam_W).
\]

A morphism
\[
f\colon\mathcal{E}(P_0)=(W_0,P_0)
\longrightarrow\mathcal{E}(P_1)=(W_1,P_1)
\]
is a morphism $f\in\Hom_{\Foam^\varepsilon}(W_0,W_1)$ satisfying
\[
P_1\cdot f=f=f\cdot P_0.
\]

There exists a formal grading of objects in
$\Kar(\Foam^\varepsilon)$ which respects the degree of
morphisms between the objects. An object $\mathcal{E}(P)$ with formal
grading $i\in\mathbb{Z}$ is denoted by
\[
q^i\mathcal{E}(P).
\]
Similarly, write
\[
q^i\mathcal{E}(W):=q^i\mathcal{E}(\idfoam_W),
\]
and for any polynomial $p=\sum_ip_iq^i \in \C(q)$, 
\[
p\mathcal{E}(W):=\osum_i(\osum_{p_i}q^i\mathcal{E}(W)).
\]

The morphism space between objects $q^i\mathcal{E}(P_0)$ and
$q^j\mathcal{E}(P_1)$ is defined by
\[
\Hom_{\Kar(\Foam^\varepsilon)}
\bigl(q^i\mathcal{E}(P_0),q^j\mathcal{E}(P_1)\bigr)
=P_1\,\Hom_{\Foam^\varepsilon}^{\,i-j}(W_0,W_1)\,P_0.
\]

In other words, a morphism
\[
f\colon q^i\mathcal{E}(P_0)
\longrightarrow q^j\mathcal{E}(P_1)
\]
is a homogeneous morphism $f\in\Hom_{\Foam^\varepsilon}(W_0,W_1)$ such
that
\[
\deg(f)=i-j
\qquad\text{and}\qquad
P_1\cdot f=f=f\cdot P_0.
\]
\end{definition}

% \begin{remark}
% The formal grading of objects in $\Kar(\Foam^\varepsilon)$
% matches the grading shift of the corresponding state spaces. To be more
% precise, given $\mathcal{E}(P_0)=(W_0,P_0)$ and
% $\mathcal{E}(P_1)=(W_1,P_1)\in\Kar(\Foam^\varepsilon)$,
% consider a morphism
% \[
% f\colon\mathcal{E}(P_0)\{i\}
% \longrightarrow\mathcal{E}(P_1)\{j\},
% \]
% where $f=M$ and $\deg(f)=\foamdeg(M)=i-j$. Given any web $W_2$ with
% boundary $\varepsilon$, since
% $\foamdeg(\overline{\idfoam_{W_2}}\ast M)=\foamdeg(M)=i-j$, we have the
% grading-preserving homomorphism
% \[
% \mathcal{F}(M):
% \mathcal{F}(\overline{\idfoam_{W_2}}\ast P_0)
% \bigl(\mathcal{F}(\overline{W_2}\ast W_0)\bigr)\{i\}
% \longrightarrow
% \mathcal{F}(\overline{\idfoam_{W_2}}\ast P_1)
% \bigl(\mathcal{F}(\overline{W_2}\ast W_1)\bigr)\{j\}.
% \]
% \end{remark}

\begin{remark}
In the web category viewed after decategorification, a formal linear
combination of webs corresponds categorically to a direct sum of the
corresponding objects in $\Foam^\varepsilon$ or in
$\Kar(\Foam^\varepsilon)$. More precisely, a sum
\[
W_1+\cdots+W_r
\]
is categorified by the direct sum
\[
\mathcal{E}(W_1)\oplus\cdots\oplus\mathcal{E}(W_r).
\]

If grading shifts are present, then a graded term $q^iW$ corresponds to the
shifted object $q^i\mathcal{E}(W)$. Thus formal sums of webs in the
decategorified web category are realized by direct sums of the associated
objects in the foam category or its Karoubi completion.
\end{remark}

\section{Kekulé webs and facet insertion morphisms}

Fix a choice of $N$ in the category $\Web$ and $\Foam$ throughout the rest of the paper. 

Figure \ref{relations} shows how to simplify a rectangle appearing in a web. In \cite[Section 1.1]{Mor2007} Morrison showed how any polygon appearing in a web can be simplified. These relations are known as Kekulé relations.

The proceding definition includes many auxiliary quantities and naming conventions which we will use extensively.

\begin{definition}\label{def:kekule}
Fix $m \in \Z$, $m > 0$ and let $a,b \in \Z^m$ be two sequences of nonnegative integers, such that $a_j\leq b_j,b_{j-1}$ and set $b_0=b_m$. Introduce notation $l=\mathrm{min}\,a+N-\mathrm{max}\, b$. For each $j$, such that $\sum b\leq j\leq \sum a +N-1$, define the \emph{Kekulé relation} as a morphism in the category $\Web$:
\begin{equation}
  \label{KekuleRelation}
  \sum_{k=\mathrm{max}\,b-\sum b}^{\mathrm{min}\,a-\sum a+1}(-1)^{j+k}
  \begin{bmatrix}
    j+k-\mathrm{max}\,b\\
    j-\sum b
  \end{bmatrix}_q
  \begin{bmatrix}
    \mathrm{min}\,a+N-j-k\\
    \sum a+N-1-j
  \end{bmatrix}_q
  \left\langle
  \begin{tikzpicture}[baseline=-3,use Hobby shortcut, scale=0.8]
    \draw[thick] (1.7,0)--(1.7,0.7);
    \draw[thick,-Stealth] (1.7,0.7)--(1.1,1.3);
    \draw[thick] (1.2,1.2)--(0.7,1.7);
    \draw[thick,-Stealth] (0.7,1.7)--(-0.1,1.7);
    \draw[thick,-Stealth] (0,1.7)--(-0.7,1.7)--(-1.3,1.1);
    \draw[thick,-Stealth] (-1.2,1.2)--(-1.7,0.7)--(-1.7,-0.1);
    \draw[thick,-Stealth] (-1.7,0)--(-1.7,-0.7)--(-1.1,-1.3);
    \draw[thick,-Stealth] (-1.2,-1.2)--(-0.7,-1.7)--(0.1,-1.7);
    \draw[thick,-Stealth,dashed] (0,-1.7)--(0.7,-1.7)--(1.3,-1.1); 
    \draw[thick,-Stealth,dashed] (1.2,-1.2)--(1.7,-0.7)--(1.7,0.1); 

    \draw[thick,-Stealth] (1.7,0.7)--(1.7*1.4,0.7*1.4);
    \draw[thick] (1.7*1.3,0.7*1.3)--(1.7*1.7,0.7*1.7);
    \draw[thick,-Stealth] (0.7*1.7,1.7*1.7)--(0.7*1.3,1.7*1.3);
    \draw[thick] (0.7,1.7)--(0.7*1.3,1.7*1.3);
    \draw[thick] (-0.7*1.7,1.7*1.7)--(-0.7*1.3,1.7*1.3);
    \draw[thick,-Stealth] (-0.7,1.7)--(-0.7*1.4,1.7*1.4);
    \draw[thick] (-1.7,0.7)--(-1.7*1.4,0.7*1.4);
    \draw[thick,Stealth-] (-1.7*1.3,0.7*1.3)--(-1.7*1.7,0.7*1.7);
    \draw[thick,-Stealth] (-1.7,-0.7)--(-1.7*1.4,-0.7*1.4);
    \draw[thick] (-1.7*1.3,-0.7*1.3)--(-1.7*1.7,-0.7*1.7);
    \draw[thick,-Stealth] (-0.7*1.7,-1.7*1.7)--(-0.7*1.3,-1.7*1.3);
    \draw[thick] (-0.7,-1.7)--(-0.7*1.3,-1.7*1.3);
    \draw[thick,dashed] (0.7*1.7,-1.7*1.7)--(0.7*1.3,-1.7*1.3);
    \draw[thick,-Stealth,dashed] (0.7,-1.7)--(0.7*1.4,-1.7*1.4);
    \draw[thick,dashed] (1.7,-0.7)--(1.7*1.4,-0.7*1.4);
    \draw[thick,Stealth-,dashed] (1.7*1.3,-0.7*1.3)--(1.7*1.7,-0.7*1.7);
    
    \draw[thick,blue] (0,0) circle (1.2);
    \draw[thick,blue,-Stealth] (0,1.2)--(-0.1,1.2);
    \node[blue] at (0,0.7) {$j+k$};
    \draw[thick,blue,-Stealth] (-0.1,2.2)--(0,2.2);
    \node[blue] at (0,2.7) {$a_1$};
    \draw[thick,blue,-Stealth] (-2.2,0)--(-2.2,0.1);
    \node[blue] at (-2.7,0) {$a_2$};
    \draw[thick,blue,-Stealth] (1.6,1.6)--(1.7,1.5);
    \node[blue] at (2,2) {$b_m$};
    \draw[thick,blue,-Stealth] (-1.6,1.6)--(-1.5,1.7);
    \node[blue] at (-2,2) {$b_1$};
    \draw[thick,blue,-Stealth] (-1.6,-1.6)--(-1.7,-1.5);
    \node[blue] at (-2,-2) {$b_2$};
    \begin{knot}[clip width=3.7,flip crossing/.list={1},]
      \strand[thick,blue] (0.7,1.7*1.7)..++(-0.01,-0.02)..(0,2.2)..++(-0.01,0)..(-0.7,1.7*1.7)..++(-0.01,0.02);
      \strand[thick,blue] (-1.7*1.7,-0.7)..++(0.02,0.01)..(-2.2,0)..++(0,0.01)..(-1.7*1.7,0.7)..++(-0.02,0.01);
      \strand[thick,blue] (2.65,1.55)..++(-0.02,-0.01)..(1.6,1.6)..++(-0.01,0.01)..(1.55,2.65)..++(0.01,0.02);
      \strand[thick,blue] (-2.65,1.55)..++(0.02,-0.01)..(-1.6,1.6)..++(0.01,0.01)..(-1.55,2.65)..++(-0.01,0.02);
      \strand[thick,blue] (-2.65,-1.55)..++(0.02,0.01)..(-1.6,-1.6)..++(0.01,-0.01)..(-1.55,-2.65)..++(-0.01,-0.02);
    \end{knot}
  \end{tikzpicture}
  \right\rangle
\end{equation}
Each label is understood to be the signed sum of the blue arrows on either side.
%and when the label is a negative integer $-c$, we understand it as if the arrow had the opposite orientation with label $c$.

The webs in Equation~\eqref{KekuleRelation} are called \emph{Kekulé webs}, and we denote them as $\Gamma_{j+k-\mathrm{max}\,b}$ respective to the index in the equation. We will call the web edges forming the cycle graph around the blue circle labeled $j+k$ \emph{internal (edges)} and we will call the other edges \emph{boundary (edges)} of a Kekulé web. 
%The same notation also expresses the corresponding $\Symf$ module as in Definition \ref{def:foam between closed webs}. 
We will use the notation $\Gamma_i$ for an arbitrary Kekulé web when the sequences $a,b$ are clear from the context or not relevant to the discussion. 
\end{definition}

Note that $\Gamma_{i}=0$ when $i> l$ or $i<0$. Moreover, note that $i$ represents the minimum label on the non boundary edges of the web $\Gamma_i$.

\begin{remark}
    Morrison showed that all Kekulé relations correspond to the zero morphism of $\Web$ (see \cite[ Theorem 5.3.2]{Mor2007}). That is, the formula in Equation \eqref{KekuleRelation} equals $0$.
\end{remark}

\begin{definition}\label{def:smoothing}
Pick a choice of $m, a,$ and $b$ from Definition \ref{def:kekule}. Let $\Gamma_i$ be some Kekulé web for that choice of sequences. Let $c$ be a coloring of the foam $\idfoam_{\Gamma_i}$. Any such coloring induces a coloring of $\Gamma_i$ in the obvious way. The tuple $(\Gamma_i,c)$ is called a \emph{colored web}. 

Locally, a colored web looks like the left side of Figure~\ref{fig:Smoothing}. Any colored web can be smoothed into a union of (not neccesary disjoint) curves in the plane, as shown in Figure \ref{fig:Smoothing}. We will call the collection of these curves the \emph{smoothing} of $(\Gamma_i,c)$ and denote it by $\smoothing(\Gamma_i,c)$. 
\begin{figure}[htbp]
  \centering
  \begin{tikzpicture}[
    x=1cm,
    y=1cm,
    stroke/.style={
        line width=1pt,
        line cap=round,
        line join=round
    },
    >={Stealth[length=2.4mm,width=1.8mm]},
    every node/.style={font=\large}
]

% 左侧分叉图
\coordinate (O) at (0,0);

\draw[stroke,->] (O) -- (-1.05,1.05);
\draw[stroke,->] (O) -- ( 1.05,1.05);
\draw[stroke,-<] (O) -- (0,-1.25);

\node at (-1.23,1.28) {$i$};
\node at ( 1.23,1.28) {$j$};

\node at (-0.88,0.43) {$A$};
\node at ( 0.88,0.43) {$B$};
\node at ( 0.78,-0.82) {$A\sqcup B$};
\node at (0,-1.62) {$i+j$};

% 中间箭头
\draw[stroke,->] (1.75,0) -- (3.55,0);

% 右侧展开图
\begin{scope}[xshift=5.1cm]

    % 属于 A 的支路
    \draw[stroke,<-<]
        (-0.85,1.05) -- (0,0) -- (0,-1.05);

    \draw[stroke,<-<]
        (0.38,1.05) -- (1.23,0) -- (1.23,-1.05);

    % 属于 B 的支路
    \draw[stroke,<-<]
        (2.58,1.05) -- (1.68,0) -- (1.68,-1.05);

    \draw[stroke,<-<]
        (3.82,1.05) -- (2.92,0) -- (2.92,-1.05);

    % 省略号
    \node at (0.52,0) {$\cdots$};
    \node at (2.40,0) {$\cdots$};

    % 顶部标记
    \node at (0.20,1.48)
        {$a_1,\cdots,a_i\in A$};

    \node at (3.8,1.48)
        {$b_1,\cdots,b_j\in B$};

\end{scope}

\end{tikzpicture}
\caption{The smoothing of a colored web}
\label{fig:Smoothing}
\end{figure}
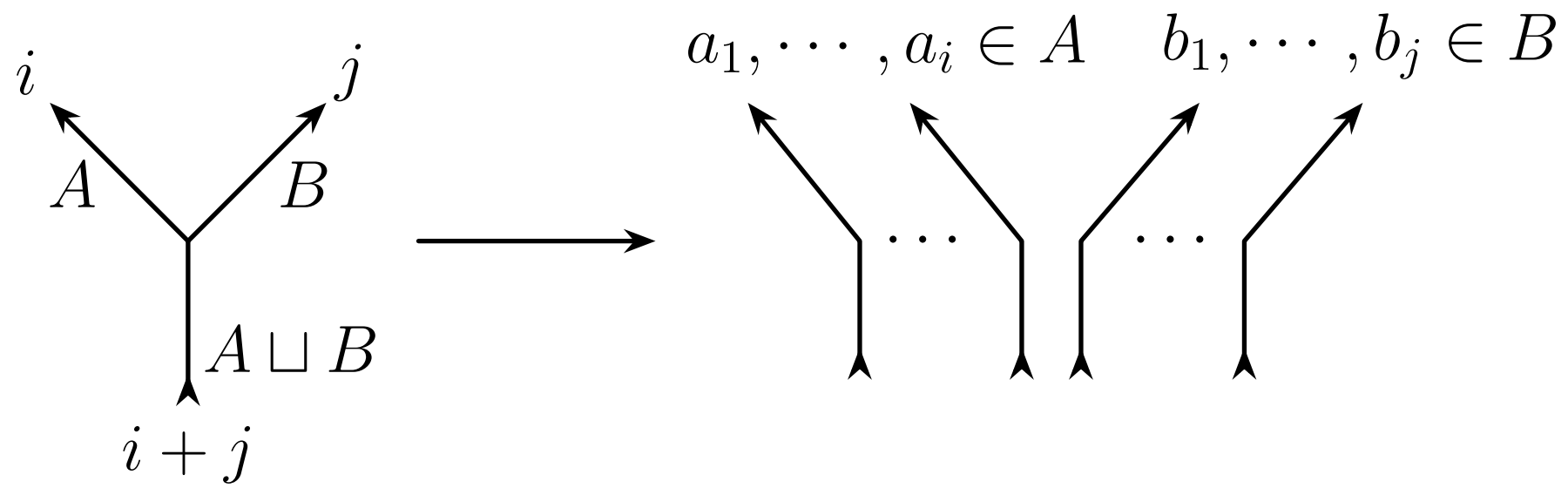
Each color from $c$ corresponds to some curve in the smoothing of $(\Gamma_i,c)$. We introduce notation for \emph{color sets of $(\Gamma_i,c)$} (any of which may be empty):
\begin{itemize}
  %\item $O$ denotes the subset of colors not used in $c$, 
  \item $A$ denotes the subset of colors of $c$, which correspond to the circles in $\smoothing(\Gamma_i,c)$,
  \item $B_k$ denotes the subset of colors of $c$, where each of the colors corresponds to $k$ many disjoint open curves, 
  \item $O$ denotes the subset of remaining colors, which are not used in the coloring $c$.
\end{itemize}
Hence the set of all colors (that is $\{1, \ldots, N\}$) equals $O\sqcup A\sqcup(\sqcup_k B_k)$. We will sometimes understand the color sets as sequences instead of sets.

We interpret $B:=\sqcup_k B_k$ as a multiset with multiplicity induced by the disjoint union. In other words, the colors in $B_k$ have multiplicity $k$ and $|B|=\sum b-\sum a$. Lastly, we introduce the notation $r:=N-1-|B|$. 
\end{definition}

The following example will be built upon throughout the rest of the paper.

\begin{example}
  \label{EgSmoothing}
  Consider $\Gamma_1$ with $m=3$, $a=(0,0,0)$, $b=(1,1,1)$, adn $N \geq 4$. Consider the coloring $c$ of $\Gamma_1$ depicted on the left of Figure \ref{fig:EgSmoothing}. The set $\smoothing(\Gamma_1, c)$ consists of one circle and three curves as depicted on the right of Figure \ref{fig:EgSmoothing}. Hence, $O = \{1, \ldots, N\}-\{a,b_1,b_2,b_3\}$, $A=\{a\}$, $B=B_1=\{b_1,b_2,b_3\}$.
  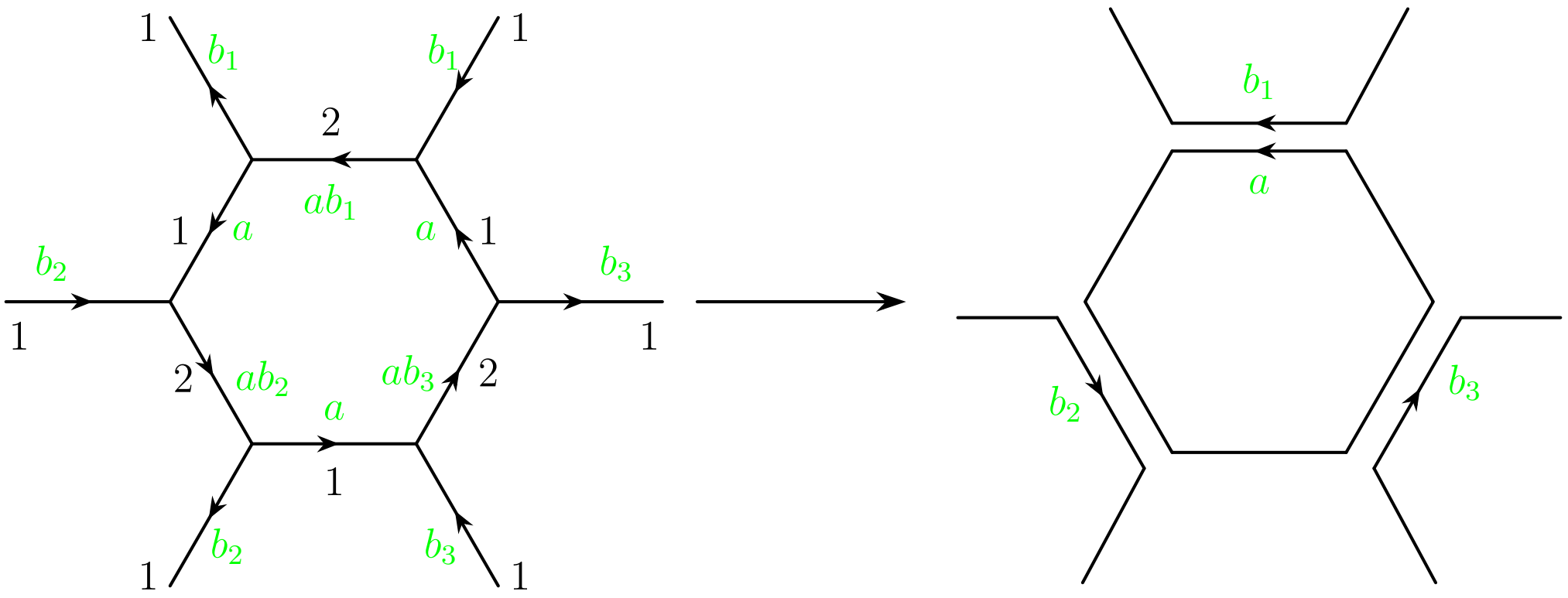
\begin{figure}[htbp]
    \centering
    \newcommand{\ExSmoothing}{0.85}
    \scalebox{\ExSmoothing}{%
    \begin{tikzpicture}[
  scale=1,
  line/.style={
    line width=0.9pt,
    line cap=round,
    line join=round
  },
  flow/.style={
    line,
    postaction={decorate},
    decoration={
      markings,
      mark=at position 0.53 with
      {\arrow{Stealth[length=2.2mm,width=1.7mm]}}
    }
  },
  every node/.style={
    font=\normalsize,
    inner sep=1.5pt
  }
]

%==================================================
% 左图
%==================================================
\begin{scope}

% 正六边形的外接圆半径
\def\R{1.65}
\pgfmathsetmacro{\H}{sqrt(3)*\R/2}

% 正六边形顶点
\coordinate (LT) at (-0.5*\R, \H);
\coordinate (RT) at ( 0.5*\R, \H);
\coordinate (RM) at ( \R,0);
\coordinate (RB) at ( 0.5*\R,-\H);
\coordinate (LB) at (-0.5*\R,-\H);
\coordinate (LM) at (-\R,0);

% 外侧端点
\coordinate (LTout) at (-1*\R, 2*\H);
\coordinate (RTout) at ( 1*\R, 2*\H);
\coordinate (RMout) at ( 2*\R,0);
\coordinate (RBout) at ( 1*\R,-2*\H);
\coordinate (LBout) at (-1*\R,-2*\H);
\coordinate (LMout) at (-2*\R,0);

% 六边形上的有向边
\draw[flow] (RT) -- (LT)
  node[pos=.52,above=5pt] {$2$}
  node[pos=.52,below=5pt,green] {$ab_1$};

\draw[flow] (LT) -- (LM)
  node[pos=.50,right=4pt,green] {$a$}
  node[pos=.50,left=4pt]  {$1$};

\draw[flow] (LM) -- (LB)
  node[pos=.54,right=4pt,green] {$ab_2$}
  node[pos=.54,left=4pt]  {$2$};

\draw[flow] (LB) -- (RB)
  node[pos=.50,above=5pt,green] {$a$}
  node[pos=.50,below=5pt] {$1$};

\draw[flow] (RB) -- (RM)
  node[pos=.50,left=4pt,green]  {$ab_3$}
  node[pos=.50,right=4pt] {$2$};

\draw[flow] (RM) -- (RT)
  node[pos=.50,left=4pt,green]  {$a$}
  node[pos=.50,right=4pt] {$1$};

% 六条向外的支路
\draw[flow] (LT) -- (LTout)
  node[pos=.76,right=3pt,green] {$b_1$}
  node[pos=.93,left=3pt]  {$1$};

\draw[flow] (RTout) -- (RT)
  node[pos=.24,left=3pt,green]  {$b_1$}
  node[pos=.07,right=3pt] {$1$};

\draw[flow] (LMout) -- (LM)
  node[pos=.28,above=4pt,green] {$b_2$}
  node[pos=.08,below=4pt] {$1$};

\draw[flow] (RM) -- (RMout)
  node[pos=.72,above=4pt,green] {$b_3$}
  node[pos=.92,below=4pt] {$1$};

\draw[flow] (LB) -- (LBout)
  node[pos=.72,right=3pt,green] {$b_2$}
  node[pos=.93,left=3pt]  {$1$};

\draw[flow] (RBout) -- (RB)
  node[pos=.28,left=3pt,green]  {$b_3$}
  node[pos=.07,right=3pt] {$1$};

\end{scope}

%==================================================
% 中间箭头
%==================================================
\draw[
  line,
  -{Stealth[length=3mm,width=2mm]}
] (3.65,0) -- (5.75,0);

%==================================================
% 右图
%==================================================
\begin{scope}[xshift=9.3cm]

% 正六边形半径
\def\RR{1.75}
\pgfmathsetmacro{\HH}{sqrt(3)*\RR/2}

% 内部正六边形顶点
\coordinate (ILT) at (-0.5*\RR, \HH);
\coordinate (IRT) at ( 0.5*\RR, \HH);
\coordinate (IRM) at ( \RR,0);
\coordinate (IRB) at ( 0.5*\RR,-\HH);
\coordinate (ILB) at (-0.5*\RR,-\HH);
\coordinate (ILM) at (-\RR,0);

% 内部正六边形
\draw[line] (ILT) -- (ILM) -- (ILB);
\draw[line] (ILB) -- (IRB);
\draw[line] (IRB) -- (IRM) -- (IRT);

% 顶边 a
\draw[flow] (IRT) -- (ILT)
  node[pos=.50,below=5pt,green] {$a$};

% 顶部平行边 b_1
\coordinate (OLT) at (-0.5*\RR,\HH+0.28);
\coordinate (ORT) at ( 0.5*\RR,\HH+0.28);

\draw[flow] (ORT) -- (OLT)
  node[pos=.50,above=5pt,green] {$b_1$};

% 顶部外支路
\draw[line]
  (OLT) -- ++(-0.62,1.15);

\draw[line]
  (ORT) -- ++(0.62,1.15);

% 左下平行边 b_2
\coordinate (LouterM) at (-\RR-0.28,-0.16);
\coordinate (LouterB) at (-0.5*\RR-0.28,-\HH-0.16);

\draw[line]
  (LouterM) -- ++(-1.0,0);

\draw[flow] (LouterM) -- (LouterB)
  node[pos=.57,left=5pt,green] {$b_2$};

\draw[line]
  (LouterB) -- ++(-0.62,-1.15);

% 右下平行边 b_3
\coordinate (RouterB) at (0.5*\RR+0.28,-\HH-0.16);
\coordinate (RouterM) at (\RR+0.28,-0.16);

\draw[flow] (RouterB) -- (RouterM)
  node[pos=.57,right=5pt,green] {$b_3$};

\draw[line]
  (RouterM) -- ++(1.0,0);

\draw[line]
  (RouterB) -- ++(0.62,-1.15);

\end{scope}

\end{tikzpicture}}
  \caption{An example of the smoothing of a Kekulé web.}
  \label{fig:EgSmoothing}
  \end{figure}
\end{example}

Next we consider a special set of morphisms between two Kekulé webs with the same defining sequences $a,b$.

\begin{definition}\label{def:slice morphism}
  Fix a choice of $m, a, b$ from Definition \ref{def:kekule}. Let $\varepsilon$ be the corresponding boundary sequence, and $0\leq i\neq j\leq l $. A \emph{facet insertion morphism (foam)} is a foam $F_i^j \in \Hom_{\Foam^\varepsilon}(\Gamma_i, \Gamma_j)$, constructed as a cylinder with an $|i-j|$-facet in the middle (see Example \ref{eg:slice morphism}). We call the middle facet \emph{the inserted facet}. For $j=i\pm 1$, we call $F_i^j$ \emph{a simple facet insertion morphism (foam)}. If $i$ or $j$ are outside of the specified bounds, $F_i^j$ is understood to be $0$. For $f \in \mathrm{Sym}_{|i-j|}$, we denote by $F_i^j(f)$ a facet insertion morphism $F_i^j$ which carries $f$ as a decoration on the inserted facet. It is to be understood that $F_i^j = F_i^j(1)$.
\end{definition} 

In Section \ref{sec:munu}, the so-called ``duality'' (insertion) morphisms will play an important role. The name is chosen to parallel the isomorphism between the fundamental representations of $U_q(\GL_N)$ given by (see \cite[Section 3.1]{CKM2014})
\begin{align*}
    \wedge_q^k\mathbb{C}_q^N \cong \wedge_q^N\mathbb{C}_q^N \otimes \left(\wedge_q^{N-k}\mathbb{C}_q^N\right)^*.
\end{align*}

\begin{definition}\label{def:duality morphism}
    Fix a choice of $m, a, b$ from Definition \ref{def:kekule}. Let $\varepsilon$ be the corresponding boundary sequence, and $0\leq i,j\leq l $. We denote by $\Gamma_i'$ \emph{the web dual to} $\Gamma_i$. That is a web obtained from $\Gamma_i$ by replacing each internal edge labeled with $i$ by an oppositely oriented edge labeled by $N -i$. A \emph{duality morphisms (foam)} is a foam $D$ from $\Gamma_i$ to $\Gamma_i'$, given by inserting a facet of thickness $N$, oriented oppositely compared to orientation of the facets inserted in Definition \ref{def:slice morphism}.
\end{definition}

\begin{remark}\label{rmk:D4 = id}
    It is shown in \cite[Lemma 3.21]{RW2020} that $D^4 = \idfoam$.
\end{remark}

\begin{example}\label{eg:slice morphism}
  Consider $\Gamma_1$ with $m=3$, $a=(0,0,0)$, $b=(1,1,1)$ as in Example \ref{EgSmoothing}. A facet insertion morphism from $\Gamma_1$ to $\Gamma_0$ is depicted in Figure~\ref{fig:ExSliceFoam} as a foam, and in Figure~\ref{fig:ExMovieStyle} as a movie. Maps are read from bottom to top or left to right. 

  \begin{figure}[htbp]
    \centering 
    \newcommand{\ExSliceFoamScale}{0.7}
    \scalebox{\ExSliceFoamScale}{%
    \begin{tikzpicture}[
  x=.9cm,y=.9cm,
  base/.style={
    line width=.75pt,
    line cap=round,
    line join=round
  },
  dline/.style args={#1}{
    base,
    postaction={decorate},
    decoration={
      markings,
      mark=at position #1 with
      {\arrow{Stealth[length=2.2mm,width=1.5mm]}}
    }
  },
  dline/.default=.5,
  twoarrows/.style={
    base,
    postaction={decorate},
    decoration={
      markings,
      mark=at position .28 with
      {\arrow{Stealth[length=2.2mm,width=1.5mm]}},
      mark=at position .72 with
      {\arrow{Stealth[length=2.2mm,width=1.5mm]}}
    }
  },
  every node/.style={
    font=\small,
    inner sep=1.5pt
  }
]

%==================================================
% 坐标
%==================================================

% 左右两个外侧矩形
\coordinate (LoutT) at (-6, 3.75);
\coordinate (LinT)  at (-4.5,3.75);
\coordinate (LinB)  at (-4.5,-3.75);
\coordinate (LoutB) at (-6,-3.75);

\coordinate (RinT)  at (4.5,3.75);
\coordinate (RoutT) at (6,3.75);
\coordinate (RoutB) at (6,-3.75);
\coordinate (RinB)  at (4.5,-3.75);

% 第一组中央矩形
\coordinate (TL) at (-1.9, 5.25);
\coordinate (TR) at ( 1.9, 5.25);
\coordinate (BL) at (-1.9,-2.25);
\coordinate (BR) at ( 1.9,-2.25);

% 第二组中央矩形
\coordinate (UL) at (-2.65, 2.25);
\coordinate (UR) at ( 2.65, 2.25);
\coordinate (DL) at (-2.65,-5.25);
\coordinate (DR) at ( 2.65,-5.25);

% 顶部两侧结构
\coordinate (LcapTop) at (-3,6);
\coordinate (LcapLow) at (-3,-1.5);
\coordinate (RcapTop) at ( 3,6);
\coordinate (RcapLow) at ( 3,-1.5);

% 下部两侧结构
\coordinate (LsideTop) at (-3.75, 1.5);
\coordinate (LsideBot) at (-3.75,-6);
\coordinate (RsideTop) at ( 3.75, 1.5);
\coordinate (RsideBot) at ( 3.75,-6);

% 中央绿色六边形
\coordinate (GL)  at (-4.5,0);
\coordinate (GTL) at (-1.9,1.5);
\coordinate (GTR) at ( 1.9,1.5);
\coordinate (GR)  at ( 4.5,0);
\coordinate (GBR) at ( 2.65,-1.5);
\coordinate (GBL) at (-2.65,-1.5);

%==================================================
% 中央绿色区域
%==================================================

\fill[green!35,opacity=.35]
  (GL)--(GTL)--(GTR)--(GR)--(GBR)--(GBL)--cycle;

%==================================================
% 左右外侧矩形
%==================================================

\draw[dline=.5] (LoutT)--(LinT)
  node[pos=.5,below=4pt] {$1$};

\draw[twoarrows] (LinT)--(LinB);

\draw[dline=.5] (LoutB)--(LinB)
  node[pos=.5,below=4pt] {$1$};

\draw[base] (LoutB)--(LoutT);

\draw[dline=.5] (RinT)--(RoutT)
  node[pos=.5,above=4pt] {$1$};

\draw[twoarrows] (RinB)--(RinT);

\draw[dline=.5] (RinB)--(RoutB)
  node[pos=.5,below=4pt] {$1$};

\draw[base] (RoutT)--(RoutB);

%==================================================
% 两个相交的中央矩形
%==================================================

\draw[dline=.5] (TR)--(TL)
  node[pos=.5,above=4pt] {$1$};

\draw[twoarrows] (BL)--(TL);
\draw[twoarrows] (TR)--(BR);

\draw[dline=.5] (BR)--(BL)
  node[pos=.5,below=5pt] {$2$};

\draw[base] (UL)--(UR)
  node[pos=.5,above=4pt] {$0$};

\draw[twoarrows] (DL)--(UL);
\draw[twoarrows] (UR)--(DR);

\draw[dline=.5] (DL)--(DR)
  node[pos=.5,below=4pt] {$1$};

%==================================================
% 顶部两侧结构
%==================================================

\draw[base] (LcapTop)--(LcapLow);

\draw[dline=.55] (LcapTop)--(TL)
  node[pos=.5,above=3pt] {$1$};

\draw[dline=.55] (LcapLow)--(BL)
  node[pos=.5,below=3pt] {$1$};

\draw[base] (RcapLow)--(RcapTop);

\draw[dline=.55] (RcapTop)--(TR)
  node[pos=.5,above=3pt] {$1$};

\draw[dline=.55] (RcapLow)--(BR)
  node[pos=.5,below=3pt] {$1$};

%==================================================
% 上部对角线
%==================================================

\draw[base] (LinT)--(TL)
  node[pos=.48,above=4pt] {$0$};

\draw[dline=.54] (LinT)--(UL)
  node[pos=.52,below=4pt] {$1$};

\draw[base] (TR)--(RinT)
  node[pos=.52,above=4pt] {$0$};

\draw[dline=.54] (UR)--(RinT)
  node[pos=.52,below=4pt] {$1$};

%==================================================
% 两侧细长结构
%==================================================

\draw[base] (LsideTop)--(LsideBot);

\draw[dline=.52] (UL)--(LsideTop)
  node[pos=.55,below=4pt] {$1$};

\draw[dline=.52] (DL)--(LsideBot)
  node[pos=.56,below=4pt] {$1$};

\draw[base] (RsideBot)--(RsideTop);

\draw[dline=.52] (RsideTop)--(UR)
  node[pos=.55,above=4pt] {$1$};

\draw[dline=.52] (RsideBot)--(DR)
  node[pos=.56,below=4pt] {$1$};

%==================================================
% 下部对角线
%==================================================

\draw[dline=.52] (BL)--(LinB)
  node[pos=.52,above=4pt] {$1$};

\draw[dline=.52] (LinB)--(DL)
  node[pos=.52,below=4pt] {$2$};

\draw[dline=.52] (RinB)--(BR)
  node[pos=.52,above=4pt] {$1$};

\draw[dline=.52] (DR)--(RinB)
  node[pos=.52,below=4pt] {$2$};

%==================================================
% 绿色六边形边界
%==================================================

\draw[dline=.53] (GL)--(GTL);
\draw[dline=.50] (GTL)--(GTR);
\draw[dline=.53] (GTR)--(GR);
\draw[dline=.53] (GR)--(GBR);
\draw[dline=.50] (GBR)--(GBL);
\draw[dline=.53] (GBL)--(GL);

\node[green] at (0,0) {$1$};

\end{tikzpicture}}
    \caption{A foam depicting a facet insertion morphism. }
    \label{fig:ExSliceFoam}
  \end{figure}
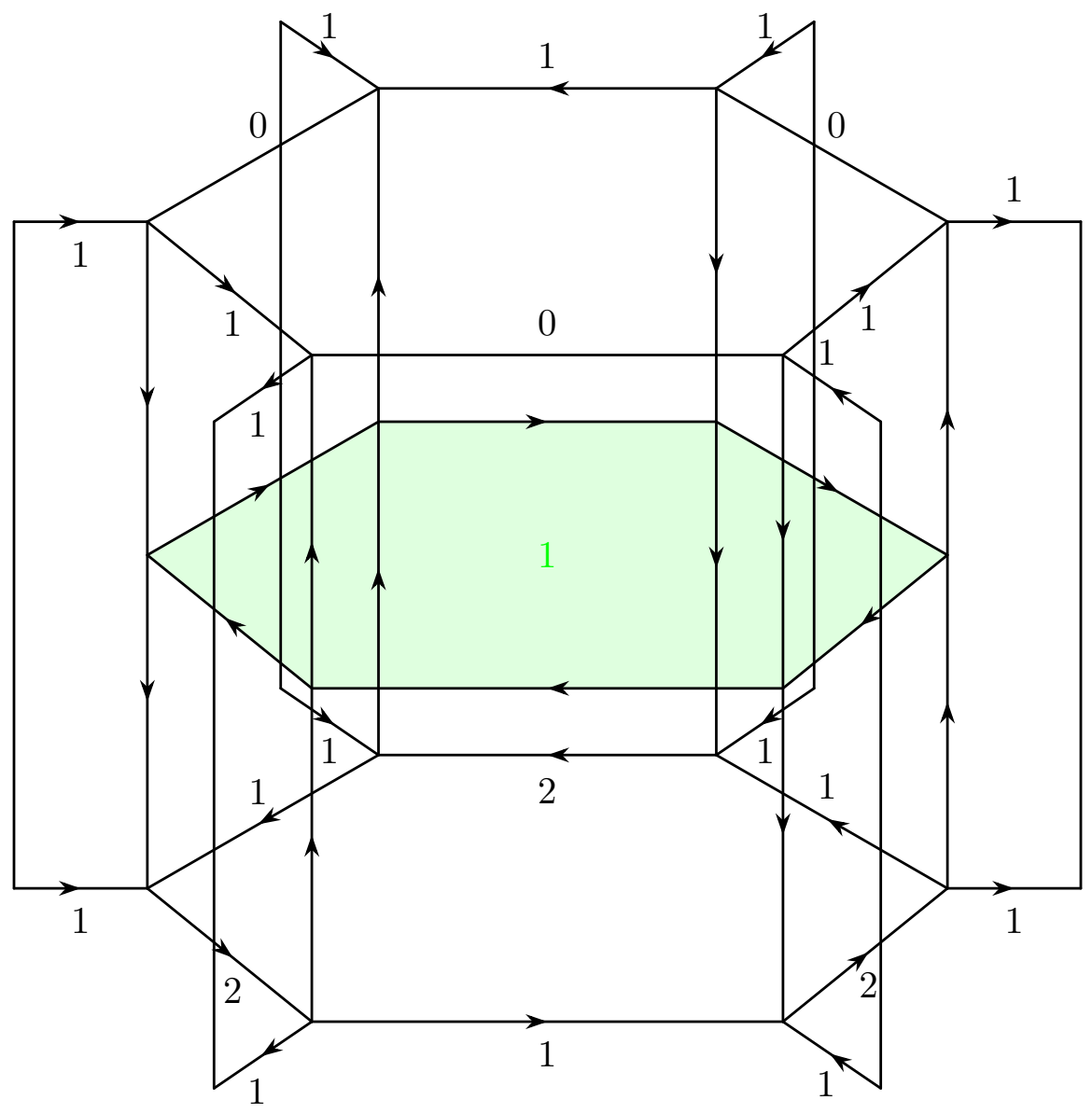

  \begin{figure}[htbp]
    \centering
    \newcommand{\ExMovieStyleScale}{0.85}
    \scalebox{\ExMovieStyleScale}{%
    \begin{tikzpicture}[
  scale=0.75,
  line/.style={
    line width=0.9pt,
    line cap=round,
    line join=round
  },
  flow/.style={
    line,
    postaction={decorate},
    decoration={
      markings,
      mark=at position 0.53 with
      {\arrow{Stealth[length=2.2mm,width=1.7mm]}}
    }
  },
  every node/.style={
    font=\normalsize,
    inner sep=1.5pt
  }
]

%==================================================
% 左图
%==================================================
\begin{scope}

% 正六边形的外接圆半径
\def\R{1.65}
\pgfmathsetmacro{\H}{sqrt(3)*\R/2}

% 正六边形顶点
\coordinate (LT) at (-0.5*\R, \H);
\coordinate (RT) at ( 0.5*\R, \H);
\coordinate (RM) at ( \R,0);
\coordinate (RB) at ( 0.5*\R,-\H);
\coordinate (LB) at (-0.5*\R,-\H);
\coordinate (LM) at (-\R,0);

% 外侧端点
\coordinate (LTout) at (-1*\R, 2*\H);
\coordinate (RTout) at ( 1*\R, 2*\H);
\coordinate (RMout) at ( 2*\R,0);
\coordinate (RBout) at ( 1*\R,-2*\H);
\coordinate (LBout) at (-1*\R,-2*\H);
\coordinate (LMout) at (-2*\R,0);

% 六边形上的有向边
\draw[flow] (RT) -- (LT)
  node[pos=.52,above=5pt] {$2$}
  node[pos=.52,below=5pt,green] {$ab_1$};

\draw[flow] (LT) -- (LM)
  node[pos=.50,right=4pt,green] {$a$}
  node[pos=.50,left=4pt]  {$1$};

\draw[flow] (LM) -- (LB)
  node[pos=.54,right=4pt,green] {$ab_2$}
  node[pos=.54,left=4pt]  {$2$};

\draw[flow] (LB) -- (RB)
  node[pos=.50,above=5pt,green] {$a$}
  node[pos=.50,below=5pt] {$1$};

\draw[flow] (RB) -- (RM)
  node[pos=.50,left=4pt,green]  {$ab_3$}
  node[pos=.50,right=4pt] {$2$};

\draw[flow] (RM) -- (RT)
  node[pos=.50,left=4pt,green]  {$a$}
  node[pos=.50,right=4pt] {$1$};

% 六条向外的支路
\draw[flow] (LT) -- (LTout)
  node[pos=.76,right=3pt,green] {$b_1$}
  node[pos=.93,left=3pt]  {$1$};

\draw[flow] (RTout) -- (RT)
  node[pos=.24,left=3pt,green]  {$b_1$}
  node[pos=.07,right=3pt] {$1$};

\draw[flow] (LMout) -- (LM)
  node[pos=.28,above=4pt,green] {$b_2$}
  node[pos=.08,below=4pt] {$1$};

\draw[flow] (RM) -- (RMout)
  node[pos=.72,above=4pt,green] {$b_3$}
  node[pos=.92,below=4pt] {$1$};

\draw[flow] (LB) -- (LBout)
  node[pos=.72,right=3pt,green] {$b_2$}
  node[pos=.93,left=3pt]  {$1$};

\draw[flow] (RBout) -- (RB)
  node[pos=.28,left=3pt,green]  {$b_3$}
  node[pos=.07,right=3pt] {$1$};

\end{scope}

\begin{scope}[xshift=7cm]

% 正六边形的外接圆半径
\def\R{1.65}
\pgfmathsetmacro{\H}{sqrt(3)*\R/2}

\node at (-0.3,0) {$1$};
\node[green] at (0.3,0) {$a$};

% 正六边形顶点
\coordinate (LT) at (-0.5*\R, \H);
\coordinate (RT) at ( 0.5*\R, \H);
\coordinate (RM) at ( \R,0);
\coordinate (RB) at ( 0.5*\R,-\H);
\coordinate (LB) at (-0.5*\R,-\H);
\coordinate (LM) at (-\R,0);

% 外侧端点
\coordinate (LTout) at (-1*\R, 2*\H);
\coordinate (RTout) at ( 1*\R, 2*\H);
\coordinate (RMout) at ( 2*\R,0);
\coordinate (RBout) at ( 1*\R,-2*\H);
\coordinate (LBout) at (-1*\R,-2*\H);
\coordinate (LMout) at (-2*\R,0);

% 六边形上的有向边
\draw[flow] (LT) -- (RT);
\draw[flow] (LM) -- (LT);
\draw[flow] (LB) -- (LM);
\draw[flow] (RB) -- (LB);
\draw[flow] (RM) -- (RB);
\draw[flow] (RT) -- (RM);

% 六条向外的支路
\draw[line] (LT) -- (LTout);
\draw[line] (RTout) -- (RT);
\draw[line] (LMout) -- (LM);
\draw[line] (RM) -- (RMout);
\draw[line] (LB) -- (LBout);
\draw[line] (RBout) -- (RB);
\end{scope}

\begin{scope}[xshift=14cm]

% 正六边形的外接圆半径
\def\R{1.65}
\pgfmathsetmacro{\H}{sqrt(3)*\R/2}

% 正六边形顶点
\coordinate (LT) at (-0.5*\R, \H);
\coordinate (RT) at ( 0.5*\R, \H);
\coordinate (RM) at ( \R,0);
\coordinate (RB) at ( 0.5*\R,-\H);
\coordinate (LB) at (-0.5*\R,-\H);
\coordinate (LM) at (-\R,0);

% 外侧端点
\coordinate (LTout) at (-1*\R, 2*\H);
\coordinate (RTout) at ( 1*\R, 2*\H);
\coordinate (RMout) at ( 2*\R,0);
\coordinate (RBout) at ( 1*\R,-2*\H);
\coordinate (LBout) at (-1*\R,-2*\H);
\coordinate (LMout) at (-2*\R,0);

% 六边形上的有向边
\draw[flow] (RT) -- (LT)
  node[pos=.52,above=5pt] {$1$}
  node[pos=.52,below=5pt,green] {$b_1$};

\draw[flow] (LT) -- (LM)
  node[pos=.50,left=4pt]  {$0$};

\draw[flow] (LM) -- (LB)
  node[pos=.54,right=4pt,green] {$b_2$}
  node[pos=.54,left=4pt]  {$1$};

\draw[flow] (LB) -- (RB)
  node[pos=.50,below=5pt] {$0$};

\draw[flow] (RB) -- (RM)
  node[pos=.50,left=4pt,green]  {$b_3$}
  node[pos=.50,right=4pt] {$1$};

\draw[flow] (RM) -- (RT)
  node[pos=.50,right=4pt] {$0$};

% 六条向外的支路
\draw[flow] (LT) -- (LTout)
  node[pos=.76,right=3pt,green] {$b_1$}
  node[pos=.93,left=3pt]  {$1$};

\draw[flow] (RTout) -- (RT)
  node[pos=.24,left=3pt,green]  {$b_1$}
  node[pos=.07,right=3pt] {$1$};

\draw[flow] (LMout) -- (LM)
  node[pos=.28,above=4pt,green] {$b_2$}
  node[pos=.08,below=4pt] {$1$};

\draw[flow] (RM) -- (RMout)
  node[pos=.72,above=4pt,green] {$b_3$}
  node[pos=.92,below=4pt] {$1$};

\draw[flow] (LB) -- (LBout)
  node[pos=.72,right=3pt,green] {$b_2$}
  node[pos=.93,left=3pt]  {$1$};

\draw[flow] (RBout) -- (RB)
  node[pos=.28,left=3pt,green]  {$b_3$}
  node[pos=.07,right=3pt] {$1$};

\end{scope}
\end{tikzpicture}}
    \caption{A movie between colored webs.}
    \label{fig:ExMovieStyle}
  \end{figure}
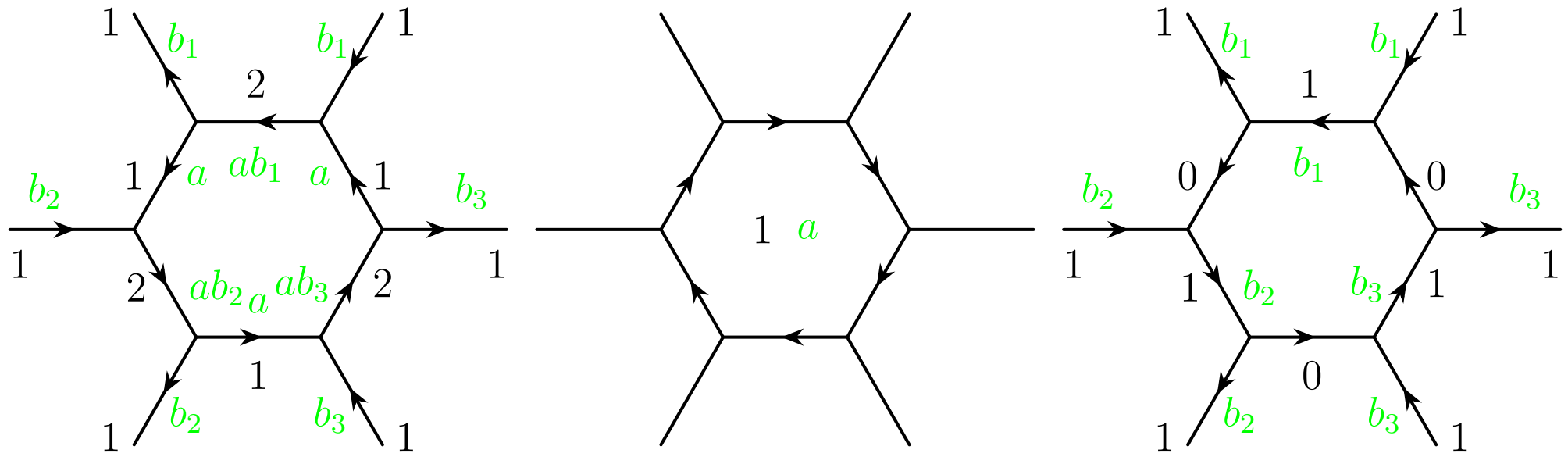

  For clarity, we omit the facet labels besides the label of the inserted facet, and their colors in Figure~\ref{fig:ExSliceFoam}. While the movie depiction allows us to keep the colors as in Figure \ref{fig:ExMovieStyle}. Movies are read as in the work of Bar-Natan \cite{Bar2005}. Recall that the orientations of the hexagons agree with whether they increase (or decrease) the labels on the neighboring frames in a movie. For instance in Figure \ref{fig:ExMovieStyle}, the middle frame decreases the labels of the neighboring ones when read from left to right. 

  Using these choices we depict the duality morphism in Figure \ref{fig:duality}.
 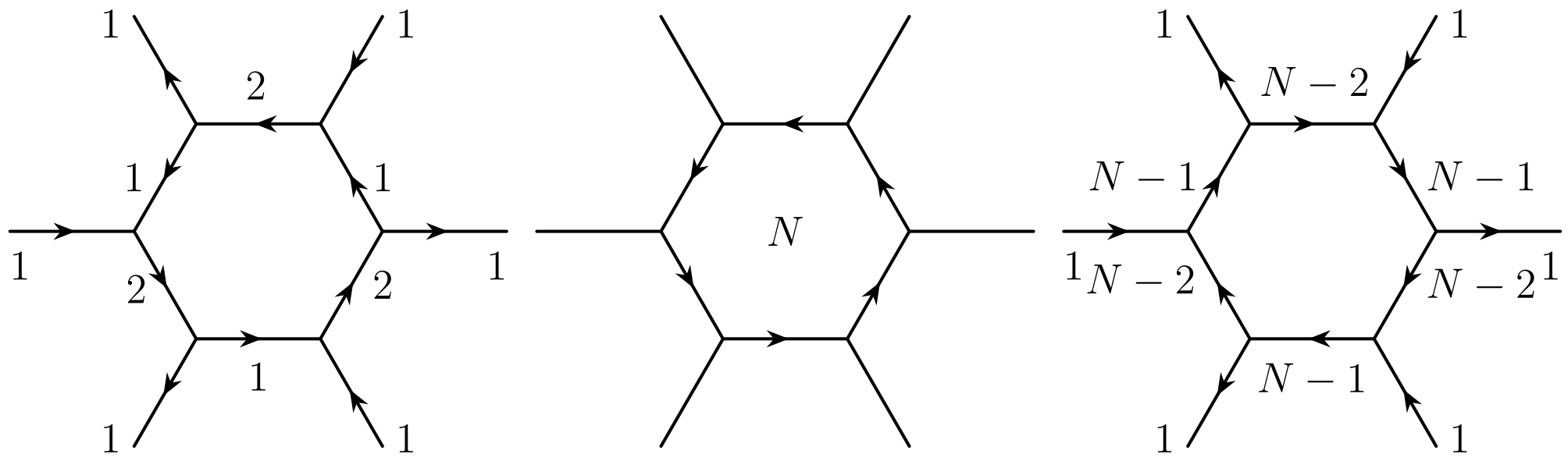
\begin{figure}[htbp]
    \centering
    \newcommand{\ExMovieStyleScale}{0.85}
    \scalebox{\ExMovieStyleScale}{%
    \begin{tikzpicture}[
  scale=0.75,
  line/.style={
    line width=0.9pt,
    line cap=round,
    line join=round
  },
  flow/.style={
    line,
    postaction={decorate},
    decoration={
      markings,
      mark=at position 0.53 with
      {\arrow{Stealth[length=2.2mm,width=1.7mm]}}
    }
  },
  every node/.style={
    font=\normalsize,
    inner sep=1.5pt
  }
]

%==================================================
% 左图
%==================================================
\begin{scope}

% 正六边形的外接圆半径
\def\R{1.65}
\pgfmathsetmacro{\H}{sqrt(3)*\R/2}

% 正六边形顶点
\coordinate (LT) at (-0.5*\R, \H);
\coordinate (RT) at ( 0.5*\R, \H);
\coordinate (RM) at ( \R,0);
\coordinate (RB) at ( 0.5*\R,-\H);
\coordinate (LB) at (-0.5*\R,-\H);
\coordinate (LM) at (-\R,0);

% 外侧端点
\coordinate (LTout) at (-1*\R, 2*\H);
\coordinate (RTout) at ( 1*\R, 2*\H);
\coordinate (RMout) at ( 2*\R,0);
\coordinate (RBout) at ( 1*\R,-2*\H);
\coordinate (LBout) at (-1*\R,-2*\H);
\coordinate (LMout) at (-2*\R,0);

% 六边形上的有向边
\draw[flow] (RT) -- (LT)
  node[pos=.52,above=5pt] {$2$};

\draw[flow] (LT) -- (LM)
  node[pos=.50,left=4pt]  {$1$};

\draw[flow] (LM) -- (LB)
  node[pos=.54,left=4pt]  {$2$};

\draw[flow] (LB) -- (RB)
  node[pos=.50,below=5pt] {$1$};

\draw[flow] (RB) -- (RM)
  node[pos=.50,right=4pt] {$2$};

\draw[flow] (RM) -- (RT)
  node[pos=.50,right=4pt] {$1$};

% 六条向外的支路
\draw[flow] (LT) -- (LTout)
  node[pos=.93,left=3pt]  {$1$};

\draw[flow] (RTout) -- (RT)
  node[pos=.07,right=3pt] {$1$};

\draw[flow] (LMout) -- (LM)
  node[pos=.08,below=4pt] {$1$};

\draw[flow] (RM) -- (RMout)
  node[pos=.92,below=4pt] {$1$};

\draw[flow] (LB) -- (LBout)
  node[pos=.93,left=3pt]  {$1$};

\draw[flow] (RBout) -- (RB)
  node[pos=.07,right=3pt] {$1$};

\end{scope}

\begin{scope}[xshift=7cm]

% 正六边形的外接圆半径
\def\R{1.65}
\pgfmathsetmacro{\H}{sqrt(3)*\R/2}

\node at (0,0) {$N$};

% 正六边形顶点
\coordinate (LT) at (-0.5*\R, \H);
\coordinate (RT) at ( 0.5*\R, \H);
\coordinate (RM) at ( \R,0);
\coordinate (RB) at ( 0.5*\R,-\H);
\coordinate (LB) at (-0.5*\R,-\H);
\coordinate (LM) at (-\R,0);

% 外侧端点
\coordinate (LTout) at (-1*\R, 2*\H);
\coordinate (RTout) at ( 1*\R, 2*\H);
\coordinate (RMout) at ( 2*\R,0);
\coordinate (RBout) at ( 1*\R,-2*\H);
\coordinate (LBout) at (-1*\R,-2*\H);
\coordinate (LMout) at (-2*\R,0);

% 六边形上的有向边
\draw[flow] (RT) -- (LT);
\draw[flow] (LT) -- (LM);
\draw[flow] (LM) -- (LB);
\draw[flow] (LB) -- (RB);
\draw[flow] (RB) -- (RM);
\draw[flow] (RM) -- (RT);

% 六条向外的支路
\draw[line] (LT) -- (LTout);
\draw[line] (RTout) -- (RT);
\draw[line] (LMout) -- (LM);
\draw[line] (RM) -- (RMout);
\draw[line] (LB) -- (LBout);
\draw[line] (RBout) -- (RB);
\end{scope}

\begin{scope}[xshift=14cm]

% 正六边形的外接圆半径
\def\R{1.65}
\pgfmathsetmacro{\H}{sqrt(3)*\R/2}

% 正六边形顶点
\coordinate (LT) at (-0.5*\R, \H);
\coordinate (RT) at ( 0.5*\R, \H);
\coordinate (RM) at ( \R,0);
\coordinate (RB) at ( 0.5*\R,-\H);
\coordinate (LB) at (-0.5*\R,-\H);
\coordinate (LM) at (-\R,0);

% 外侧端点
\coordinate (LTout) at (-1*\R, 2*\H);
\coordinate (RTout) at ( 1*\R, 2*\H);
\coordinate (RMout) at ( 2*\R,0);
\coordinate (RBout) at ( 1*\R,-2*\H);
\coordinate (LBout) at (-1*\R,-2*\H);
\coordinate (LMout) at (-2*\R,0);

% 六边形上的有向边
\draw[flow] (LT) -- (RT)
  node[pos=.52,above=5pt] {$N - 2$};

\draw[flow] (LM) -- (LT)
  node[pos=.50,left=4pt]  {$N - 1$};

\draw[flow] (LB) -- (LM)
  node[pos=.54,left=4pt]  {$N - 2$};

\draw[flow] (RB) -- (LB)
  node[pos=.50,below=5pt] {$N - 1$};

\draw[flow] (RM) -- (RB)
  node[pos=.50,right=4pt] {$N - 2$};

\draw[flow] (RT) -- (RM)
  node[pos=.50,right=4pt] {$N - 1$};

% 六条向外的支路
\draw[flow] (LT) -- (LTout)
  node[pos=.93,left=3pt]  {$1$};

\draw[flow] (RTout) -- (RT)
  node[pos=.07,right=3pt] {$1$};

\draw[flow] (LMout) -- (LM)
  node[pos=.08,below=4pt] {$1$};

\draw[flow] (RM) -- (RMout)
  node[pos=.92,below=4pt] {$1$};

\draw[flow] (LB) -- (LBout)
  node[pos=.93,left=3pt]  {$1$};

\draw[flow] (RBout) -- (RB)
  node[pos=.07,right=3pt] {$1$};

\end{scope}
\end{tikzpicture}}
    \caption{A movie representing the duality morphism on $\Gamma_1$.}
    \label{fig:duality}
\end{figure}
\end{example}

\begin{lemma}\label{lem:deg compuation}
 Fix a choice of $m, a, b$ from Definition \ref{def:kekule}.  Let $F_i^{i\pm 1}:\Gamma_i\to\Gamma_{i\pm 1}$, be a simple facet insertion morphism. Then
  \[\foamdeg(F_i^{i\pm 1}(h_s))=2s-r.\]
  \label{DegreeLem}
\end{lemma}
\begin{proof}
    This follows by a degree computation. The computation is straightforward by following the approach in Remark \ref{rmk:degree with decomposition}.
\end{proof}

\begin{example}
    For the choice of $m,a$ and $b$ from Example \ref{EgSmoothing}, it holds that $|B| = 3, r = N - 4$ and $l = N - 1$. Hence the simple facet insertion morphisms (without any decorations) for this choice of parameters all have degree $4 - N$.
\end{example}

In the following sections, we will categorify a simpler set of web relations that arise naturally from our constructions. These relations appear in the next lemma, where it is shown that they are equivalent to the Kekulé relations. 

\begin{lemma}\label{lem:true kekule}
Fix $m \in \Z$, $m > 0$ and let $a,b \in \Z_{\geq 0}^m$ be two sequences of nonnegative integers. The collection of web relations for $l \geq i \geq l - r$
\begin{align}\label{eq:true kekule}
    0 = \sum_{j=0}^{i}(-1)^j\begin{bmatrix}
    r+i-j\\
    r
  \end{bmatrix}_q\langle\Gamma_j\rangle,
\end{align}
is equivalent to
\[\sum_{i=s}^{l-r+s}(-1)^{i}\begin{bmatrix}
    i\\s
  \end{bmatrix}_q\begin{bmatrix}
    l-i\\r-s
  \end{bmatrix}_q\langle\Gamma_i\rangle=0,\]
  and agrees with the Kekulé web relations from Definition \ref{def:kekule}.
\end{lemma}
\begin{proof}
The equations \eqref{eq:true kekule} can be presented as a finite system of equations where the $\Gamma_i$'s are regarded as formal symbols:
  \[ \left(\begin{array}{ccccc}
    %%first row
    \begin{bmatrix}
    l\\
    r
    \end{bmatrix}_q
    &
    \begin{bmatrix}
    l-1\\
    r
    \end{bmatrix}_q
    &
    \begin{bmatrix}
    l-2\\
    r
    \end{bmatrix}_q
    &
    \cdots 
    &
    \begin{bmatrix}
    0\\
    r
    \end{bmatrix}_q
    \\
    %second row
    \begin{bmatrix}
    l+1\\
    r
    \end{bmatrix}_q
    &
    \begin{bmatrix}
    l\\
    r
    \end{bmatrix}_q
    &
    \begin{bmatrix}
    l-1\\
    r
    \end{bmatrix}_q
    &
    \cdots 
    &
    \begin{bmatrix}
    1\\
    r
    \end{bmatrix}_q
    \\
    %third row
    \begin{bmatrix}
    l+2\\
    r
    \end{bmatrix}_q
    &
    \begin{bmatrix}
    l+1\\
    r
    \end{bmatrix}_q
    &
    \begin{bmatrix}
    l\\
    r
    \end{bmatrix}_q
    &
    \cdots 
    &
    \begin{bmatrix}
    2\\
    r
    \end{bmatrix}_q
    \\
    %...
    \vdots&\vdots&\vdots& &\vdots\\
    %last row
    \begin{bmatrix}
    l+r \\
    r 
    \end{bmatrix}_q
    &
    \begin{bmatrix}
    l+r-1\\
    r
    \end{bmatrix}_q
    &
    \begin{bmatrix}
    l+r-2\\
    r
    \end{bmatrix}_q
    &
    \cdots 
    &
    \begin{bmatrix}
    r \\
    r
    \end{bmatrix}_q
  \end{array}\right)
  \left(\begin{array}{c}
    \langle\Gamma_0\rangle\\-\langle\Gamma_1\rangle\\\langle\Gamma_2\rangle\\\vdots\\(-1)^l\langle\Gamma_l\rangle
  \end{array}\right)=0.\]
  Denote the matrix on the left by $M$. We wish to perform Gaussian elimination on $M$ to make it upper trapezoidal. This requires using the equation $q^{-n}[m-n]_q+q^{m-n}[n]_q=[m]_q$, and the quantum version of the Pascal's formula:
  \begin{align*}
      \begin{bmatrix}
    m \\
    n
    \end{bmatrix}_q=q^{-n}\begin{bmatrix}
    m-1 \\
    n
    \end{bmatrix}_q+q^{m-n}\begin{bmatrix}
    m-1 \\
    n-1
    \end{bmatrix}_q.   \end{align*}
   To reduce $M$, multiply each row of $M$, except the first, by $q^{r}$, and subtract the $i$-th row from the $i+1$-st row of $M$. Using the Pascal's formula this gives:
    \[M \sim \left(\begin{array}{ccccc}
    %%first row
    \begin{bmatrix}
    l\\
    r
    \end{bmatrix}_q
    &
    \begin{bmatrix}
    l-1\\
    r
    \end{bmatrix}_q
    &
    \begin{bmatrix}
    l-2\\
    r
    \end{bmatrix}_q
    &
    \cdots 
    &
    \begin{bmatrix}
    0\\
    r
    \end{bmatrix}_q
    \\
    %second row
    q^{l+1}\begin{bmatrix}
    l\\
    r-1
    \end{bmatrix}_q
    &
    q^{l}\begin{bmatrix}
    l-1\\
    r-1
    \end{bmatrix}_q
    &
    q^{l-1}\begin{bmatrix}
    l-2\\
    r-1
    \end{bmatrix}_q
    &
    \cdots 
    &
    q\begin{bmatrix}
    0\\
    r-1
    \end{bmatrix}_q
    \\
    %third row
    q^{l+2}\begin{bmatrix}
    l+1\\
    r-1
    \end{bmatrix}_q
    &
    q^{l+1}\begin{bmatrix}
    l\\
    r-1
    \end{bmatrix}_q
    &
    q^l\begin{bmatrix}
    l-1\\
    r-1
    \end{bmatrix}_q
    &
    \cdots 
    &
    q^2\begin{bmatrix}
    1\\
    r-1
    \end{bmatrix}_q
    \\
    %...
    \vdots&\vdots&\vdots& &\vdots\\
    %last row
    q^{l+r}\begin{bmatrix}
    l+r-1 \\
    r -1
    \end{bmatrix}_q
    &
    q^{l+r-1}\begin{bmatrix}
    l+r-2\\
    r-1
    \end{bmatrix}_q
    &
q^{l+r-2}\begin{bmatrix}l+r-3\\r-1\end{bmatrix}_q
    &
    \cdots 
    &
    q^{r}\begin{bmatrix}
    r-1 \\
    r-1
    \end{bmatrix}_q
  \end{array}\right).\]
  Next, divide the $i$-th row by $q^i$. This ensures that the powers of $q$ are the same for the last $r$ entries of each column. This allows us to use Pascal's formula again. Repeating these steps we obtain
    \[M \sim \left(\begin{array}{ccccc}
    %%first row
    \begin{bmatrix}
    l\\
    r
    \end{bmatrix}_q
    &
    \begin{bmatrix}
    l-1\\
    r
    \end{bmatrix}_q
    &
    \begin{bmatrix}
    l-2\\
    r
    \end{bmatrix}_q
    &
    \cdots 
    &
    \begin{bmatrix}
    0\\
    r
    \end{bmatrix}_q
    \\
    %second row
    q^{l}\begin{bmatrix}
    l\\
    r-1
    \end{bmatrix}_q
    &
    q^{l-1}\begin{bmatrix}
    l-1\\
    r-1
    \end{bmatrix}_q
    &
    q^{l-2}\begin{bmatrix}
    l-2\\
    r-1
    \end{bmatrix}_q
    &
    \cdots 
    &
    \begin{bmatrix}
    0\\
    r-1
    \end{bmatrix}_q
    \\
    %third row
    q^{2l}\begin{bmatrix}
    l\\
    r-2
    \end{bmatrix}_q
    &
    q^{2(l-1)}\begin{bmatrix}
    l-1\\
    r-2
    \end{bmatrix}_q
    &
    q^{2(l-2)}\begin{bmatrix}
    l-2\\
    r-2
    \end{bmatrix}_q
    &
    \cdots 
    &
    \begin{bmatrix}
    0\\
    r-2
    \end{bmatrix}_q
    \\
    %...
    \vdots&\vdots&\vdots& &\vdots\\
    %last row
    q^{rl}\begin{bmatrix}
    l \\
    0
    \end{bmatrix}_q
    &
    q^{r(l-1)}\begin{bmatrix}
    l-1\\
    0
    \end{bmatrix}_q
    &
    q^{r(l-2)}\begin{bmatrix}
    l-2\\
    0
    \end{bmatrix}_q
    &
    \cdots 
    &
    \begin{bmatrix}
    0 \\
    0
    \end{bmatrix}_q
  \end{array}\right).\]

  Next, we eliminate the first entry of the last $r$ rows. Note that for $k=1,2,\cdots,r$ and $j=0,1,\cdots,l$ we have
  \begin{align*}
    &q^{k(l-j)}\begin{bmatrix}
      l-j\\r-k
    \end{bmatrix}_q\begin{bmatrix}
      l\\r-k+1
    \end{bmatrix}_q-q^{kl-(k-1)j}\begin{bmatrix}
      l-j\\r-k+1
    \end{bmatrix}_q\begin{bmatrix}
      l\\r-k
    \end{bmatrix}_q\\
    =&q^{(k-1)(l-j)+r-k}\left[j\right]_q\begin{bmatrix}
      l-j\\r-k
    \end{bmatrix}_q\frac{\left[l\right]_q!}{\left[r-k+1\right]_q!\left[l-r+k\right]_q!}.
  \end{align*}
  Besides from the first row, we can eliminate the first entry of the $i$-th row by using the above equality and the $i-1$-st row. Doing so on all rows we get:
  \[M \sim \left(\begin{array}{ccccc}
    %%first row
    \begin{bmatrix}
    l\\
    r
    \end{bmatrix}_q
    &
    \begin{bmatrix}
    l-1\\
    r
    \end{bmatrix}_q
    &
    \begin{bmatrix}
    l-2\\
    r
    \end{bmatrix}_q
    &
    \cdots 
    &
    \begin{bmatrix}
    0\\
    r
    \end{bmatrix}_q
    \\
    %second row
    q^{0}\left[0\right]_q\begin{bmatrix}
    l\\
    r-1
    \end{bmatrix}_q
    &
    q^{0}\left[1\right]_q\begin{bmatrix}
    l-1\\
    r-1
    \end{bmatrix}_q
    &
    q^{0}\left[2\right]_q\begin{bmatrix}
    l-2\\
    r-1
    \end{bmatrix}_q
    &
    \cdots 
    &
    q^{0}\left[l\right]_q\begin{bmatrix}
    0\\
    r-1
    \end{bmatrix}_q
    \\
    %third row
    q^{l}\left[0\right]_q\begin{bmatrix}
    l\\
    r-2
    \end{bmatrix}_q
    &
    q^{l-1}\left[1\right]_q\begin{bmatrix}
    l-1\\
    r-2
    \end{bmatrix}_q
    &
    q^{l-2}\left[2\right]_q\begin{bmatrix}
    l-2\\
    r-2
    \end{bmatrix}_q
    &
    \cdots 
    &
    q^{0}\left[l\right]_q\begin{bmatrix}
    0\\
    r-2
    \end{bmatrix}_q
    \\
    %...
    \vdots&\vdots&\vdots& &\vdots\\
    %last row
    q^{(r-1)l}\left[0\right]_q\begin{bmatrix}
    l \\
    0
    \end{bmatrix}_q
    &
    q^{(r-1)(l-1)}\left[1\right]_q\begin{bmatrix}
    l-1\\
    0
    \end{bmatrix}_q
    &
    q^{(r-1)(l-2)}\left[2\right]_q\begin{bmatrix}
    l-2\\
    0
    \end{bmatrix}_q
    &
    \cdots 
    &
    q^{0}\left[l\right]_q\begin{bmatrix}
    0 \\
    0
    \end{bmatrix}_q
  \end{array}\right).\]
  Repeating the above procedure for other columns gives:
    \[M \sim \left(\begin{array}{ccccc}
    %%first row
    \begin{bmatrix}
    l\\
    r
    \end{bmatrix}_q
    &
    \begin{bmatrix}
    l-1\\
    r
    \end{bmatrix}_q
    &
    \begin{bmatrix}
    l-2\\
    r
    \end{bmatrix}_q
    &
    \cdots 
    &
    \begin{bmatrix}
    0\\
    r
    \end{bmatrix}_q
    \\
    %second row
    0
    &
    \left[1\right]_q\begin{bmatrix}
    l-1\\
    r-1
    \end{bmatrix}_q
    &
    \left[2\right]_q\begin{bmatrix}
    l-2\\
    r-1
    \end{bmatrix}_q
    &
    \cdots 
    &
    \left[l\right]_q\begin{bmatrix}
    0\\
    r-1
    \end{bmatrix}_q
    \\
    %third row
    0
    &
    0
    &
    \left[2\right]_q\left[1\right]_q\begin{bmatrix}
    l-2\\
    r-2
    \end{bmatrix}_q
    &
    \cdots 
    &
    \left[l\right]_q\left[l-1\right]_q\begin{bmatrix}
    0\\
    r-2
    \end{bmatrix}_q
    \\
    %...
    \vdots& & & &\vdots\\
    %last row
    0
    &
    \cdots 
    &
    \left[r\right]_q\cdots\left[1\right]_q\begin{bmatrix}
    l-r\\
    0
    \end{bmatrix}_q
    &
    \cdots 
    &
    \left[l\right]_q\cdots\left[l-r+1\right]_q\begin{bmatrix}
    0 \\
    0
    \end{bmatrix}_q
  \end{array}\right).\]
  It follows that for $0\leq s\leq r$, we have
  \[\sum_{i=s}^{l-r+s}(-1)^{i}\begin{bmatrix}
    i\\s
  \end{bmatrix}_q\begin{bmatrix}
    l-i\\r-s
  \end{bmatrix}_q\langle\Gamma_i\rangle=0,\]
  which agrees with the Kekulé relations from Equation \eqref{KekuleRelation} by setting $i=j+k-\max b$, $l=N+\min a-\max b$, $r=N-1-\sum b+\sum a$, and $s=j-\sum b$ and adjusting the sign by $(-1)^{\max \, b}$. 
\end{proof}

\begin{remark}
    The previous lemma can be understood as a change of basis for the space $\mathcal{APR}_{a,b}^m$ from \cite[Theorem 5.3.2.]{Mor2007}.
\end{remark}

Justified by the above, we have the following definition of transformed Kekulé relations, which are more suitable for categorifying.
\begin{definition}\label{def:true kekule}
    Fix $m \in \Z$, $m > 0$ and let $a,b \in \Z_{\geq 0}^m$ be two sequences of nonnegative integers. For $l \geq i \geq l - r$, we introduce a collection of web relations, called \emph{transformed Kekulé relations}
\begin{align}\label{eq:true kekuleTwo}
    0 = \sum_{j=0}^{i}(-1)^j\begin{bmatrix}
    r+i-j\\
    r
  \end{bmatrix}_q\langle\Gamma_j\rangle.
\end{align}
\end{definition}

\begin{example}
    Let $m,a$ and $b$ be as in the Example \ref{EgSmoothing}, additionally, set $N = 4$. Then Lemma \ref{lem:true kekule} says that the only Kekulé relation is
    \begin{align*}
          0 = \sum_{k=-2}^{1}(-1)^{3+k}
  \begin{bmatrix}
    3+k-1\\
    3-3
  \end{bmatrix}_q
  \begin{bmatrix}
    1-k\\
    0
  \end{bmatrix}_q \langle\Gamma_{3+k-1}\rangle = \sum_{k=-2}^{1}(-1)^{3+k}\langle\Gamma_{2+k}\rangle 
    \end{align*}
     and it is equivalent to
    \begin{align*}
    0 = \sum_{j=0}^{3}(-1)^j\begin{bmatrix}
    3-j\\
    0
  \end{bmatrix}_q\langle\Gamma_j\rangle = \sum_{j=0}^{3}(-1)^j\langle\Gamma_j\rangle,
\end{align*}
which in this simple case also holds by reindexing the summations and multiplying by an overall sign.
\end{example}

\section{Categorification of the Kekulé relations in the homotopy category}\label{sec:hocat categorificaiton}

Fix a choice of  $m$, $a$, and $b$ throughout the rest of the paper. The potentially non zero webs are $\Gamma_0,\Gamma_1,\cdots,\Gamma_l$. Hence, for $0\leq i\leq l-1$, there is a simple facet insertion morphism from $\Gamma_i$ to $\Gamma_{i+1}$ and $\Gamma_{i+1}$ to $\Gamma_i$. 

In this section, we construct a collection of long exact sequences, whose graded Euler characteristics recover the Kekulé relations. Before we introduce the main result of this section, we need to define the analogue of elementary and complete homogeneous symmetric polynomials in the category $\Foam^\varepsilon$. For the purposes of the next sections we need to do so for multiple sets of variables.

\begin{definition}\label{def:e,h on foam}
    $E_u(x_B)$ is an endomorphism of $\Gamma_{i}$ (where $i = j+k-\mathrm{max}\,b$ is understood to run from $0$ to $l$). The notation $E_u(x_B)$ hides the sequences $a$ and $b,$ and the numbers $m$ and $i$, those will be understood from the context (we hide them to simplify the notation). For the particular choice of $i$, $E_u(x_B)$ is given by a movie
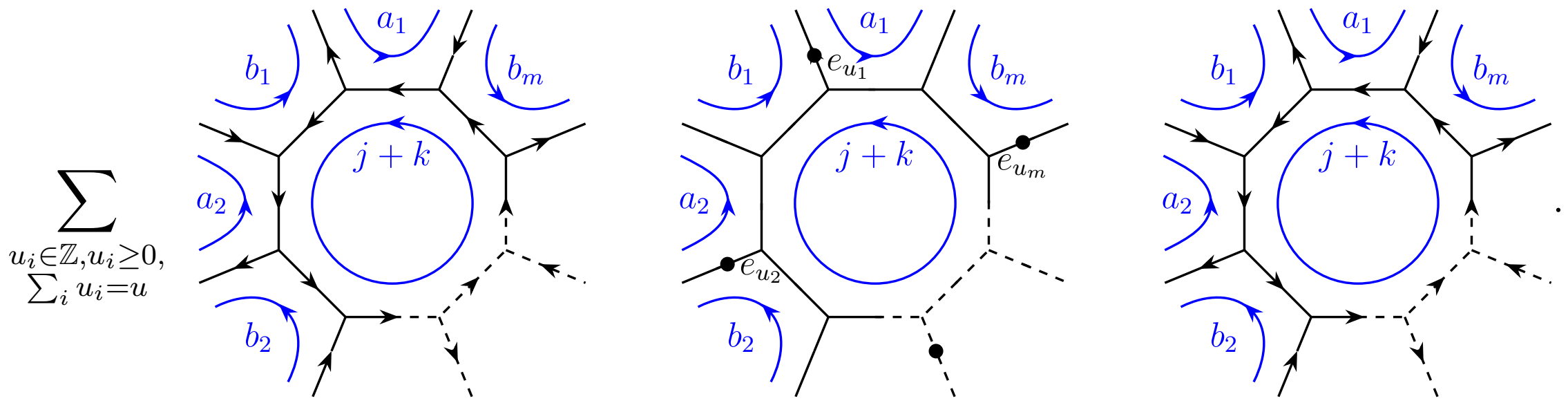
\begin{equation}\label{ek(xB)}
  \sum_{\substack{u_i \in \Z, u_i \geq 0, \\ \sum_iu_i=u}}
  \scalebox{0.8}{
  \begin{tikzpicture}[baseline=-3,use Hobby shortcut, scale=0.8]
    \draw[thick] (1.7,0)--(1.7,0.7);
    \draw[thick,-Stealth] (1.7,0.7)--(1.1,1.3);
    \draw[thick] (1.2,1.2)--(0.7,1.7);
    \draw[thick,-Stealth] (0.7,1.7)--(-0.1,1.7);
    \draw[thick,-Stealth] (0,1.7)--(-0.7,1.7)--(-1.3,1.1);
    \draw[thick,-Stealth] (-1.2,1.2)--(-1.7,0.7)--(-1.7,-0.1);
    \draw[thick,-Stealth] (-1.7,0)--(-1.7,-0.7)--(-1.1,-1.3);
    \draw[thick,-Stealth] (-1.2,-1.2)--(-0.7,-1.7)--(0.1,-1.7);
    \draw[thick,-Stealth,dashed] (0,-1.7)--(0.7,-1.7)--(1.3,-1.1); 
    \draw[thick,-Stealth,dashed] (1.2,-1.2)--(1.7,-0.7)--(1.7,0.1); 

    \draw[thick,-Stealth] (1.7,0.7)--(1.7*1.4,0.7*1.4);
    \draw[thick] (1.7*1.3,0.7*1.3)--(1.7*1.7,0.7*1.7);
    \draw[thick,-Stealth] (0.7*1.7,1.7*1.7)--(0.7*1.3,1.7*1.3);
    \draw[thick] (0.7,1.7)--(0.7*1.3,1.7*1.3);
    \draw[thick] (-0.7*1.7,1.7*1.7)--(-0.7*1.3,1.7*1.3);
    \draw[thick,-Stealth] (-0.7,1.7)--(-0.7*1.4,1.7*1.4);
    \draw[thick] (-1.7,0.7)--(-1.7*1.4,0.7*1.4);
    \draw[thick,Stealth-] (-1.7*1.3,0.7*1.3)--(-1.7*1.7,0.7*1.7);
    \draw[thick,-Stealth] (-1.7,-0.7)--(-1.7*1.4,-0.7*1.4);
    \draw[thick] (-1.7*1.3,-0.7*1.3)--(-1.7*1.7,-0.7*1.7);
    \draw[thick,-Stealth] (-0.7*1.7,-1.7*1.7)--(-0.7*1.3,-1.7*1.3);
    \draw[thick] (-0.7,-1.7)--(-0.7*1.3,-1.7*1.3);
    \draw[thick,dashed] (0.7*1.7,-1.7*1.7)--(0.7*1.3,-1.7*1.3);
    \draw[thick,-Stealth,dashed] (0.7,-1.7)--(0.7*1.4,-1.7*1.4);
    \draw[thick,dashed] (1.7,-0.7)--(1.7*1.4,-0.7*1.4);
    \draw[thick,Stealth-,dashed] (1.7*1.3,-0.7*1.3)--(1.7*1.7,-0.7*1.7);
    
    \draw[thick,blue] (0,0) circle (1.2);
    \draw[thick,blue,-Stealth] (0,1.2)--(-0.1,1.2);
    \node[blue] at (0,0.7) {$j+k$};
    \draw[thick,blue,-Stealth] (-0.1,2.2)--(0,2.2);
    \node[blue] at (0,2.7) {$a_1$};
    \draw[thick,blue,-Stealth] (-2.2,0)--(-2.2,0.1);
    \node[blue] at (-2.7,0) {$a_2$};
    \draw[thick,blue,-Stealth] (1.6,1.6)--(1.7,1.5);
    \node[blue] at (2,2) {$b_m$};
    \draw[thick,blue,-Stealth] (-1.6,1.6)--(-1.5,1.7);
    \node[blue] at (-2,2) {$b_1$};
    \draw[thick,blue,-Stealth] (-1.6,-1.6)--(-1.7,-1.5);
    \node[blue] at (-2,-2) {$b_2$};
    \begin{knot}[clip width=3.7,flip crossing/.list={1},]
      \strand[thick,blue] (0.7,1.7*1.7)..++(-0.01,-0.02)..(0,2.2)..++(-0.01,0)..(-0.7,1.7*1.7)..++(-0.01,0.02);
      \strand[thick,blue] (-1.7*1.7,-0.7)..++(0.02,0.01)..(-2.2,0)..++(0,0.01)..(-1.7*1.7,0.7)..++(-0.02,0.01);
      \strand[thick,blue] (2.65,1.55)..++(-0.02,-0.01)..(1.6,1.6)..++(-0.01,0.01)..(1.55,2.65)..++(0.01,0.02);
      \strand[thick,blue] (-2.65,1.55)..++(0.02,-0.01)..(-1.6,1.6)..++(0.01,0.01)..(-1.55,2.65)..++(-0.01,0.02);
      \strand[thick,blue] (-2.65,-1.55)..++(0.02,0.01)..(-1.6,-1.6)..++(0.01,-0.01)..(-1.55,-2.65)..++(-0.01,-0.02);
    \end{knot}
  \end{tikzpicture}
  \qquad
\begin{tikzpicture}[baseline=-3,use Hobby shortcut, scale=0.8]
    \draw[thick] (1.7,0)--(1.7,0.7);
    \draw[thick] (1.7,0.7)--(1.1,1.3);
    \draw[thick] (1.2,1.2)--(0.7,1.7);
    \draw[thick] (0.7,1.7)--(-0.1,1.7);
    \draw[thick] (0,1.7)--(-0.7,1.7)--(-1.3,1.1);
    \draw[thick] (-1.2,1.2)--(-1.7,0.7)--(-1.7,-0.1);
    \draw[thick] (-1.7,0)--(-1.7,-0.7)--(-1.1,-1.3);
    \draw[thick] (-1.2,-1.2)--(-0.7,-1.7)--(0.1,-1.7);
    \draw[thick,dashed] (0,-1.7)--(0.7,-1.7)--(1.3,-1.1); 
    \draw[thick,dashed] (1.2,-1.2)--(1.7,-0.7)--(1.7,0.1); 

    \draw[thick] (1.7,0.7)--(1.7*1.4,0.7*1.4);
    \draw[thick] (1.7*1.3,0.7*1.3)--(1.7*1.7,0.7*1.7);
    \draw[fill=black] (1.7*1.3,0.7*1.3) circle[radius=0.1];
    \node[below] at (1.7*1.3,0.7*1.3) {$e_{u_m}$};
    \draw[thick] (0.7*1.7,1.7*1.7)--(0.7*1.3,1.7*1.3);
    \draw[thick] (0.7,1.7)--(0.7*1.3,1.7*1.3);

    \draw[thick] (-0.7*1.7,1.7*1.7)--(-0.7*1.3,1.7*1.3);
    \draw[thick,-Stealth] (-0.7,1.7)--(-0.7*1.4,1.7*1.4);
    \draw[fill=black] (-0.7*1.3,1.7*1.3) circle[radius=0.1];
    \node[right] at (-0.7*1.3,1.7*1.2) {$e_{u_1}$};
    \draw[thick] (-1.7,0.7)--(-1.7*1.4,0.7*1.4);
    \draw[thick] (-1.7*1.3,0.7*1.3)--(-1.7*1.7,0.7*1.7);
    \draw[thick] (-1.7,-0.7)--(-1.7*1.4,-0.7*1.4);
    \draw[thick] (-1.7*1.3,-0.7*1.3)--(-1.7*1.7,-0.7*1.7);
    \draw[fill=black] (-1.7*1.3,-0.7*1.3) circle[radius=0.1];
    \node[right] at (-1.7*1.3,-0.7*1.45) {$e_{u_2}$};
    
    \draw[thick] (-0.7*1.7,-1.7*1.7)--(-0.7*1.3,-1.7*1.3);
    \draw[thick] (-0.7,-1.7)--(-0.7*1.3,-1.7*1.3);
    \draw[thick,dashed] (0.7*1.7,-1.7*1.7)--(0.7*1.3,-1.7*1.3);
    \draw[thick,dashed] (0.7,-1.7)--(0.7*1.4,-1.7*1.4);
    \draw[fill=black] (0.7*1.3,-1.7*1.3) circle[radius=0.1];
    \draw[thick,dashed] (1.7,-0.7)--(1.7*1.4,-0.7*1.4);
    \draw[thick,dashed] (1.7*1.3,-0.7*1.3)--(1.7*1.7,-0.7*1.7);
    
    \draw[thick,blue] (0,0) circle (1.2);
    \draw[thick,blue,-Stealth] (0,1.2)--(-0.1,1.2);
    \node[blue] at (0,0.7) {$j+k$};
    \draw[thick,blue,-Stealth] (-0.1,2.2)--(0,2.2);
    \node[blue] at (0,2.7) {$a_1$};
    \draw[thick,blue,-Stealth] (-2.2,0)--(-2.2,0.1);
    \node[blue] at (-2.7,0) {$a_2$};
    \draw[thick,blue,-Stealth] (1.6,1.6)--(1.7,1.5);
    \node[blue] at (2,2) {$b_m$};
    \draw[thick,blue,-Stealth] (-1.6,1.6)--(-1.5,1.7);
    \node[blue] at (-2,2) {$b_1$};
    \draw[thick,blue,-Stealth] (-1.6,-1.6)--(-1.7,-1.5);
    \node[blue] at (-2,-2) {$b_2$};
    \begin{knot}[clip width=3.7,flip crossing/.list={1},]
      \strand[thick,blue] (0.7,1.7*1.7)..++(-0.01,-0.02)..(0,2.2)..++(-0.01,0)..(-0.7,1.7*1.7)..++(-0.01,0.02);
      \strand[thick,blue] (-1.7*1.7,-0.7)..++(0.02,0.01)..(-2.2,0)..++(0,0.01)..(-1.7*1.7,0.7)..++(-0.02,0.01);
      \strand[thick,blue] (2.65,1.55)..++(-0.02,-0.01)..(1.6,1.6)..++(-0.01,0.01)..(1.55,2.65)..++(0.01,0.02);
      \strand[thick,blue] (-2.65,1.55)..++(0.02,-0.01)..(-1.6,1.6)..++(0.01,0.01)..(-1.55,2.65)..++(-0.01,0.02);
      \strand[thick,blue] (-2.65,-1.55)..++(0.02,0.01)..(-1.6,-1.6)..++(0.01,-0.01)..(-1.55,-2.65)..++(-0.01,-0.02);
    \end{knot}
  \end{tikzpicture}
  \qquad
    \begin{tikzpicture}[baseline=-3,use Hobby shortcut, scale=0.8]
    \draw[thick] (1.7,0)--(1.7,0.7);
    \draw[thick,-Stealth] (1.7,0.7)--(1.1,1.3);
    \draw[thick] (1.2,1.2)--(0.7,1.7);
    \draw[thick,-Stealth] (0.7,1.7)--(-0.1,1.7);
    \draw[thick,-Stealth] (0,1.7)--(-0.7,1.7)--(-1.3,1.1);
    \draw[thick,-Stealth] (-1.2,1.2)--(-1.7,0.7)--(-1.7,-0.1);
    \draw[thick,-Stealth] (-1.7,0)--(-1.7,-0.7)--(-1.1,-1.3);
    \draw[thick,-Stealth] (-1.2,-1.2)--(-0.7,-1.7)--(0.1,-1.7);
    \draw[thick,-Stealth,dashed] (0,-1.7)--(0.7,-1.7)--(1.3,-1.1); 
    \draw[thick,-Stealth,dashed] (1.2,-1.2)--(1.7,-0.7)--(1.7,0.1); 

    \draw[thick,-Stealth] (1.7,0.7)--(1.7*1.4,0.7*1.4);
    \draw[thick] (1.7*1.3,0.7*1.3)--(1.7*1.7,0.7*1.7);
    \draw[thick,-Stealth] (0.7*1.7,1.7*1.7)--(0.7*1.3,1.7*1.3);
    \draw[thick] (0.7,1.7)--(0.7*1.3,1.7*1.3);
    \draw[thick] (-0.7*1.7,1.7*1.7)--(-0.7*1.3,1.7*1.3);
    \draw[thick,-Stealth] (-0.7,1.7)--(-0.7*1.4,1.7*1.4);
    \draw[thick] (-1.7,0.7)--(-1.7*1.4,0.7*1.4);
    \draw[thick,Stealth-] (-1.7*1.3,0.7*1.3)--(-1.7*1.7,0.7*1.7);
    \draw[thick,-Stealth] (-1.7,-0.7)--(-1.7*1.4,-0.7*1.4);
    \draw[thick] (-1.7*1.3,-0.7*1.3)--(-1.7*1.7,-0.7*1.7);
    \draw[thick,-Stealth] (-0.7*1.7,-1.7*1.7)--(-0.7*1.3,-1.7*1.3);
    \draw[thick] (-0.7,-1.7)--(-0.7*1.3,-1.7*1.3);
    \draw[thick,dashed] (0.7*1.7,-1.7*1.7)--(0.7*1.3,-1.7*1.3);
    \draw[thick,-Stealth,dashed] (0.7,-1.7)--(0.7*1.4,-1.7*1.4);
    \draw[thick,dashed] (1.7,-0.7)--(1.7*1.4,-0.7*1.4);
    \draw[thick,Stealth-,dashed] (1.7*1.3,-0.7*1.3)--(1.7*1.7,-0.7*1.7);
    
    \draw[thick,blue] (0,0) circle (1.2);
    \draw[thick,blue,-Stealth] (0,1.2)--(-0.1,1.2);
    \node[blue] at (0,0.7) {$j+k$};
    \draw[thick,blue,-Stealth] (-0.1,2.2)--(0,2.2);
    \node[blue] at (0,2.7) {$a_1$};
    \draw[thick,blue,-Stealth] (-2.2,0)--(-2.2,0.1);
    \node[blue] at (-2.7,0) {$a_2$};
    \draw[thick,blue,-Stealth] (1.6,1.6)--(1.7,1.5);
    \node[blue] at (2,2) {$b_m$};
    \draw[thick,blue,-Stealth] (-1.6,1.6)--(-1.5,1.7);
    \node[blue] at (-2,2) {$b_1$};
    \draw[thick,blue,-Stealth] (-1.6,-1.6)--(-1.7,-1.5);
    \node[blue] at (-2,-2) {$b_2$};
    \begin{knot}[clip width=3.7,flip crossing/.list={1},]
      \strand[thick,blue] (0.7,1.7*1.7)..++(-0.01,-0.02)..(0,2.2)..++(-0.01,0)..(-0.7,1.7*1.7)..++(-0.01,0.02);
      \strand[thick,blue] (-1.7*1.7,-0.7)..++(0.02,0.01)..(-2.2,0)..++(0,0.01)..(-1.7*1.7,0.7)..++(-0.02,0.01);
      \strand[thick,blue] (2.65,1.55)..++(-0.02,-0.01)..(1.6,1.6)..++(-0.01,0.01)..(1.55,2.65)..++(0.01,0.02);
      \strand[thick,blue] (-2.65,1.55)..++(0.02,-0.01)..(-1.6,1.6)..++(0.01,0.01)..(-1.55,2.65)..++(-0.01,0.02);
      \strand[thick,blue] (-2.65,-1.55)..++(0.02,0.01)..(-1.6,-1.6)..++(0.01,-0.01)..(-1.55,-2.65)..++(-0.01,-0.02);
    \end{knot}
  \end{tikzpicture}}.
\end{equation}

It is understood that $E_0(x_B) = \idfoam_{\Gamma_{i}}$. We also define $H_u(x_B)$ inductively by 
\[H_u(x_B):=\sum_{v=1}^u(-1)^{v-1}E_v(x_B)H_{u-v}(x_B),\quad H_0(x_B):=\idfoam,\]
and denote
\begin{align*}
E_u(x_{\overline{B}})&:=\sum_{v=0}^u(-1)^{u-v}H_{u-v}(x_B)E_{v}(x_1,\ldots x_N),\\
H_u(x_{\overline{B}})&:=\sum_{v=0}^u(-1)^{u-v}E_{u-v}(x_B)H_{v}(x_1,\ldots x_N).
\label{def:e/h on xB,xBbar}
\end{align*}
All four of the defined foams are understood to be $0$ for negative indices.
\end{definition}

\begin{remark}
    We can deduce from Equation~\ref{ek(xB)} that $E_u(x_B)$, $H_u(x_B)$, $E_u(x_{\overline{B}})$ and $H_u(x_{\overline{B}})$ are in the center of the endomorphism algebra of $\Gamma_{i}$. Moreover, they commute with all the facet insertion morphisms.
    
    When all of the colors in $B$ are different, the definition can be understood by reading $x_B$ as $x_1,\cdots,x_{|B|}$ and $x_{\overline{B}}$ as $x_{|B|+1},\cdots,x_N$. 
\end{remark}

Using simple facet insertion morphisms we can construct two sets of endomorphisms of any $\Gamma_i$ given by $F_{i+1}^i(h_t)F_i^{i+1}(h_s)$ and $F_{i-1}^i(h_s)F_i^{i-1}(h_t)$. The following theorem tells us precisely how these two endomporhisms are related. In later sections we will modify the facet insertion morphisms, so that the right hand side of Equation \eqref{eq:Fii+1 Fi+1i}, always equals $\idfoam_{\Gamma_i}$ (see Lemma \ref{lem:Li,Ri prop}).

\begin{theorem}\label{Null-homotopyTh}
  For any  $i,t,s \in \Z$, such that $0\leq i\leq l$, $t\geq0, s\geq 0$, we have: 
  \begin{equation}\label{eq:Fii+1 Fi+1i}
    F_{i+1}^i(h_t)F_i^{i+1}(h_s)+F_{i-1}^i(h_s)F_i^{i-1}(h_t)=(-1)^{|B| + 1}H_{s+t-r}(x_{\overline{B}}).
  \end{equation}
\end{theorem}
\begin{proof}
  First, consider the map $F_{i+1}^i(h_t)F_i^{i+1}(h_s)$. Let $c_d$ be a coloring of this foam, for which the $1$-labeled facets are colored by the same color, see Figure~\ref{fig:ito(i+1)toi} for an illustration. 
  \begin{figure}[htbp]
    \centering
    \begin{tikzpicture}[
  scale=0.65,
  line/.style={
    line width=0.9pt,
    line cap=round,
    line join=round
  },
  flow/.style={
    line,
    postaction={decorate},
    decoration={
      markings,
      mark=at position 0.53 with
      {\arrow{Stealth[length=2.2mm,width=1.7mm]}}
    }
  },
  every node/.style={
    font=\normalsize,
    inner sep=1.5pt
  }
    ]
    \begin{scope}
      \node[green] at (2,2.5) {$b\cdots$};
      \draw[line] (0,2)--(4,2);
      \draw[line] (0,0)--(1,2);
      \draw[line] (3,2)--(4,0);
      \draw[flow] (0.5,1)--(0.4,0.8);
      \draw[flow] (2,2)--(1.8,2);
      \draw[flow] (3.5,1)--(3.4,1.2);
    \end{scope}

    \begin{scope}[xshift=5cm]
      \node at (2,1) {$1$};
      \node[green] at (2.5,0.5) {$d$};
      \draw[line,red] (1,2)--(3,2);
      \node[red] at (2,2.5) {$bd$}; 
      \draw[line] (0,2)--(1,2);
      \draw[line] (3,2)--(4,2);
      \draw[line] (0,0)--(1,2);
      \draw[line] (3,2)--(4,0);
      \draw[flow] (0.5,1)--(0.4,0.8);
      \draw[flow] (2,2)--(1.8,2);
      \draw[flow] (3.5,1)--(3.4,1.2);
      \draw[fill=black] (1.5,0.8) circle[radius=0.1];
      \node at (1.2,0.5) {$h_s$};
    \end{scope}

    \begin{scope}[xshift=10cm]
      \node[green] at (2,2.5) {$bd\cdots$};
      % \node[red] at (1,2) {$\odot$};
      % \node[red] at (0.7,2.5) {$bd$};
      % \node[red] at (3,2) {$\otimes$};
      % \node[red] at (3.3,2.5) {$bd$};
      \draw[line] (0,2)--(4,2);
      \draw[line] (0,0)--(1,2);
      \draw[line] (3,2)--(4,0);
      \draw[flow] (0.5,1)--(0.4,0.8);
      \draw[flow] (2,2)--(1.8,2);
      \draw[flow] (3.5,1)--(3.4,1.2);
    \end{scope}

    \begin{scope}[xshift=15cm]
      \node at (2,1) {$1$};
      \node[green] at (2.5,0.5) {$d$};
      \draw[line,red] (1,2)--(3,2);
      \node[red] at (2,2.5) {$bd$}; 
      \draw[line] (0,2)--(1,2);
      \draw[line] (3,2)--(4,2);
      \draw[line] (0,0)--(1,2);
      \draw[line] (3,2)--(4,0);
      \draw[flow] (0.5,1)--(0.6,1.2);
      \draw[flow] (2,2)--(2.2,2);
      \draw[flow] (3.5,1)--(3.6,0.8);
      \draw[fill=black] (1.5,0.8) circle[radius=0.1];
      \node at (1.2,0.5) {$h_t$};
    \end{scope}

    \begin{scope}[xshift=20cm]
      \node[green] at (2,2.5) {$b\cdots$};
      \draw[line] (0,2)--(4,2);
      \draw[line] (0,0)--(1,2);
      \draw[line] (3,2)--(4,0);
      \draw[flow] (0.5,1)--(0.4,0.8);
      \draw[flow] (2,2)--(1.8,2);
      \draw[flow] (3.5,1)--(3.4,1.2);
    \end{scope}
    \end{tikzpicture}
    \caption{The movie representing $\Gamma_i\xrightarrow{F_i^{i+1}(h_s)}\Gamma_{i+1}\xrightarrow{F_{i+1}^i(h_t)}\Gamma_i$ colored by $c_d$.}
    \label{fig:ito(i+1)toi}
  \end{figure}
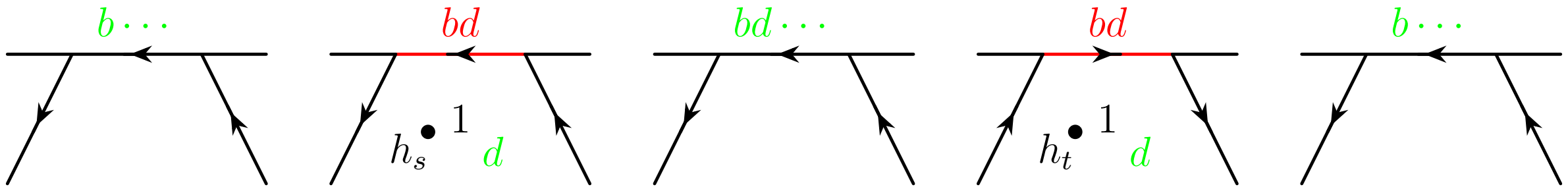
  Let $A,B,O$ be the color sets of $(\Gamma_i, c_d)$ as in Definition \ref{def:smoothing}. Note that $d\in O$. Now we compare $F_{i+1}^i(h_t)F_i^{i+1}(h_s)$ with $\idfoam_{\Gamma_i}$. For readability, we write $G =F_{i+1}^i(h_t)F_i^{i+1}(h_s) -  \idfoam_{\Gamma_i}$, by abuse of notation we use $c_d$ to denote the coloring of $G$, where $F_{i+1}^i(h_t)F_i^{i+1}(h_s)$ is colored by $c_d$ and $\idfoam_{\Gamma_i}$ is colored by $c'$, which is the coloring appearing on the leftmost web of Figure \ref{fig:ito(i+1)toi}. Observe that: 
  \begin{itemize}
    \item For any $o\in O- d$, $(-1)^{\theta^+_{do}(G, c_d)}=1$, and $\chi(G_{do}, c_d)=2$. 
    \item For any $a\in A$, $(-1)^{\theta^+_{ad}(G, c_d)}=1$, and $\chi(G_{ad}, c_d)=2$. 
    \item For any $b\in B_k$, $(-1)^{\theta^+_{bd}(G, c_d)}=(-1)^{\mathbbm{1}\{d<b\}\cdot k}$, and $\chi(G_{bd}, c_d)=-2(k-1)$. 
    \item For $d$ itself, $\chi(G_d, c_d)=2$.  
    \item Clearly, $\Delta P : = P(F_{i+1}^i(h_t)F_i^{i+1}(h_s), c_d)/P(\idfoam_{\Gamma_i}, c') =x_d^{s+t}$.
  \end{itemize}
  Denote $\Delta Q : = Q(F_{i+1}^i(h_t)F_i^{i+1}(h_s), c_d)/Q(\idfoam_{\Gamma_i},c')$ and compute: 
  \begin{align*}
    (-1)^{\Delta s}\Delta Q
    =&(-1)^{d+\sum_{k\in\mathbb{N}_+}\sum_{b\in B_k}\mathbbm{1}\{d<b\}\cdot k}\prod_{o\in O- \{d\}}(-1)^{\mathbbm{1}\{d>o\}}(x_d-x_o)\\
    &\prod_{a\in A- \{d\}}(-1)^{\mathbbm{1}\{d>a\}}(x_d-x_a)\prod_{k\in\mathbb{N}_+}\prod_{b\in B_k}[(-1)^{\mathbbm{1}\{d>b\}}(x_d-x_b)]^{-(k-1)}\\
    =&(-1)^d\prod_{c\in C- \{d\}}(-1)^{\mathbbm{1}\{d>c\}}(x_d-x_c)\prod_{k\in\mathbb{N}_+}\prod_{b\in B_k}(x_b-x_d)^{-k}.
  \end{align*}
  
Let $c_a$ be the coloring of $F_{i-1}^i(h_s)F_i^{i-1}(h_t)$ that agrees with $c_d$ from Figure \ref{fig:ito(i+1)toi} on the boundary, but whose (oppositely oriented when compared to Figure \ref{fig:ito(i+1)toi}) $1$-labeled facets are colored by some $a \in A$. All of the above computed values remain the same for $(F_{i-1}^i(h_s)F_i^{i-1}(h_t), c_a)$ when compared to $(F_{i+1}^i(h_t)F_i^{i+1}(h_s), c_d)$. 

For the purposes of the next computation, we understand that $(F_{i-1}^i(h_s)F_i^{i-1}(h_t), c_j) = 0$ for $j\notin A$ and $(F_{i+1}^i(h_t)F_i^{i+1}(h_s), c_j) = 0$ for $j\notin O$. Therefore, it suffices to consider the set of pigments $A \sqcup O$, since for any other color the following term would vanish:
  \begin{align*}
    \sum_{d\in O\sqcup A} & \frac{\langle F_{i-1}^i(h_s)F_i^{i-1}(h_t), c_i\rangle  +  \langle F_{i+1}^i(h_t)F_i^{i+1}(h_s), c_i \rangle}   {\langle \idfoam_{\Gamma_i}, c' \rangle} \\
    =& \sum_{d\in O\sqcup A}(-1)^d\frac{x_d^{s+t}\prod_{k\geq1}\prod_{b\in B_k}(x_b-x_d)^{k}}{\prod_{c\in C- \{d\}}(-1)^{\mathbbm{1}\{d>c\}}(x_d-x_c)}\\
    =&\sum_{d\in O\sqcup A}(-1)^{d+N-1}\frac{x_d^{s+t}\prod_{k\geq1}\prod_{b\in B_k}(x_b-x_d)^{k}}{\prod\limits_{i<j\atop i=d\text{ or }j=d}(x_j-x_i)}\\
    =&\prod_{1\leq i<j\leq N}(x_j-x_i)\sum_{d\in O\sqcup A}(-1)^{d+N-1}\sum_{u=0}^{|B|}(-1)^u E_{|B|-u}(x_B) x_d^{s+t+u}\prod\limits_{i<j\atop i,j\neq d}(x_j-x_i)\\
    =&\prod_{1\leq i<j\leq N}\sum_{u=0}^{|B|}(-1)^{u+1} E_{|B|-u}(x_B) 
    \begin{vmatrix}
    1 & \cdots & x_1^{N-2} & x_1^{s+t+u}\\
    1 & \cdots & x_2^{N-2} & x_2^{s+t+u}\\
    \vdots & & \vdots & \vdots \\
    1 & \cdots & x_n^{N-2} & x_n^{s+t+u}
    \end{vmatrix}\\
    =&\sum_{u=0}^{|B|}(-1)^{u+1}E_{|B|-u}(x_B)H_{s+t+u+1-N}(x_1, \ldots, x_N)\\
    =&(-1)^{|B| + 1}H_{s+t-r}(x_{\overline{B}}).
  \end{align*}
 
  This is precisely the right-hand side of Equation \ref{eq:Fii+1 Fi+1i}. It remains to show the vanishing of the sum of colorings of $F_{i+1}^i(h_t)F_i^{i+1}(h_s)$ with the two horizontal $1$-labeled facets colored by different colors. 
  
  \begin{figure}[htbp]
    \centering
    \begin{tikzpicture}[
  scale=0.65,
  line/.style={
    line width=0.9pt,
    line cap=round,
    line join=round
  },
  flow/.style={
    line,
    postaction={decorate},
    decoration={
      markings,
      mark=at position 0.53 with
      {\arrow{Stealth[length=2.2mm,width=1.7mm]}}
    }
  },
  every node/.style={
    font=\normalsize,
    inner sep=1.5pt
  }
    ]
    \begin{scope}
      \node[green] at (2,2.5) {$ab\cdots$};
      \draw[line] (0,2)--(4,2);
      \draw[line] (0,0)--(1,2);
      \draw[line] (3,2)--(4,0);
      \draw[flow] (0.5,1)--(0.4,0.8);
      \draw[flow] (2,2)--(1.8,2);
      \draw[flow] (3.5,1)--(3.4,1.2);
    \end{scope}

    \begin{scope}[xshift=5cm]
      \node at (2,1) {$1$};
      \node[green] at (2.5,0.5) {$d$};
      \draw[line,red] (1,2)--(3,2);
      \node[red] at (2,2.5) {$bd$}; 
      \draw[line] (0,2)--(1,2);
      \draw[line] (3,2)--(4,2);
      \draw[line] (0,0)--(1,2);
      \draw[line] (3,2)--(4,0);
      \draw[flow] (0.5,1)--(0.4,0.8);
      \draw[flow] (2,2)--(1.8,2);
      \draw[flow] (3.5,1)--(3.4,1.2);
      \draw[fill=black] (1.5,0.8) circle[radius=0.1];
      \node at (1.2,0.5) {$h_s$};
    \end{scope}

    \begin{scope}[xshift=10cm]
      \node[green] at (2,2.5) {$abd\cdots$};
      \draw[line] (0,2)--(4,2);
      \draw[line] (0,0)--(1,2);
      \draw[line] (3,2)--(4,0);
      \draw[flow] (0.5,1)--(0.4,0.8);
      \draw[flow] (2,2)--(1.8,2);
      \draw[flow] (3.5,1)--(3.4,1.2);
    \end{scope}

    \begin{scope}[xshift=15cm]
      \node at (2,1) {$1$};
      \node[green] at (2.5,0.5) {$a$};
      \draw[line,red] (1,2)--(3,2);
      \node[red] at (2,2.5) {$ab$}; 
      \draw[line] (0,2)--(1,2);
      \draw[line] (3,2)--(4,2);
      \draw[line] (0,0)--(1,2);
      \draw[line] (3,2)--(4,0);
      \draw[flow] (0.5,1)--(0.6,1.2);
      \draw[flow] (2,2)--(2.2,2);
      \draw[flow] (3.5,1)--(3.6,0.8);
      \draw[fill=black] (1.5,0.8) circle[radius=0.1];
      \node at (1.2,0.5) {$h_t$};
    \end{scope}

    \begin{scope}[xshift=20cm]
      \node[green] at (2,2.5) {$bd\cdots$};
      \draw[line] (0,2)--(4,2);
      \draw[line] (0,0)--(1,2);
      \draw[line] (3,2)--(4,0);
      \draw[flow] (0.5,1)--(0.4,0.8);
      \draw[flow] (2,2)--(1.8,2);
      \draw[flow] (3.5,1)--(3.4,1.2);
    \end{scope}
    \end{tikzpicture}
    $$\Gamma_i\xrightarrow{F_i^{i+1}(h_s)}\Gamma_{i+1}\xrightarrow{F_{i+1}^i(h_t)}\Gamma_i$$
    \begin{tikzpicture}[
  scale=0.65,
  line/.style={
    line width=0.9pt,
    line cap=round,
    line join=round
  },
  flow/.style={
    line,
    postaction={decorate},
    decoration={
      markings,
      mark=at position 0.53 with
      {\arrow{Stealth[length=2.2mm,width=1.7mm]}}
    }
  },
  every node/.style={
    font=\normalsize,
    inner sep=1.5pt
  }
    ]
    \begin{scope}
      \node[green] at (2,2.5) {$ab\cdots$};
      \draw[line] (0,2)--(4,2);
      \draw[line] (0,0)--(1,2);
      \draw[line] (3,2)--(4,0);
      \draw[flow] (0.5,1)--(0.4,0.8);
      \draw[flow] (2,2)--(1.8,2);
      \draw[flow] (3.5,1)--(3.4,1.2);
    \end{scope}

    \begin{scope}[xshift=5cm]
      \node at (2,1) {$1$};
      \node[green] at (2.5,0.5) {$a$};
      \draw[line,red] (1,2)--(3,2);
      \node[red] at (2,2.5) {$ab$}; 
      \draw[line] (0,2)--(1,2);
      \draw[line] (3,2)--(4,2);
      \draw[line] (0,0)--(1,2);
      \draw[line] (3,2)--(4,0);
      \draw[flow] (0.5,1)--(0.6,1.2);
      \draw[flow] (2,2)--(2.2,2);
      \draw[flow] (3.5,1)--(3.6,0.8);
      \draw[fill=black] (1.5,0.8) circle[radius=0.1];
      \node at (1.2,0.5) {$h_t$};
    \end{scope}

    \begin{scope}[xshift=10cm]
      \node[green] at (2,2.5) {$b\cdots$};
      \draw[line] (0,2)--(4,2);
      \draw[line] (0,0)--(1,2);
      \draw[line] (3,2)--(4,0);
      \draw[flow] (0.5,1)--(0.4,0.8);
      \draw[flow] (2,2)--(1.8,2);
      \draw[flow] (3.5,1)--(3.4,1.2);
    \end{scope}

    \begin{scope}[xshift=15cm]
      \node at (2,1) {$1$};
      \node[green] at (2.5,0.5) {$d$};
      \draw[line,red] (1,2)--(3,2);
      \node[red] at (2,2.5) {$bd$}; 
      \draw[line] (0,2)--(1,2);
      \draw[line] (3,2)--(4,2);
      \draw[line] (0,0)--(1,2);
      \draw[line] (3,2)--(4,0);
      \draw[flow] (0.5,1)--(0.4,0.8);
      \draw[flow] (2,2)--(1.8,2);
      \draw[flow] (3.5,1)--(3.4,1.2);
      \draw[fill=black] (1.5,0.8) circle[radius=0.1];
      \node at (1.2,0.5) {$h_s$};
    \end{scope}

    \begin{scope}[xshift=20cm]
      \node[green] at (2,2.5) {$bd\cdots$};
      \draw[line] (0,2)--(4,2);
      \draw[line] (0,0)--(1,2);
      \draw[line] (3,2)--(4,0);
      \draw[flow] (0.5,1)--(0.4,0.8);
      \draw[flow] (2,2)--(1.8,2);
      \draw[flow] (3.5,1)--(3.4,1.2);
    \end{scope}
    \end{tikzpicture}
    $$\Gamma_i\xrightarrow{F_i^{i-1}(h_t)}\Gamma_{i-1}\xrightarrow{F_{i-1}^i(h_s)}\Gamma_i$$
    \caption{When horizontal facets are colored differently}
    \label{fig:2colors}
  \end{figure}
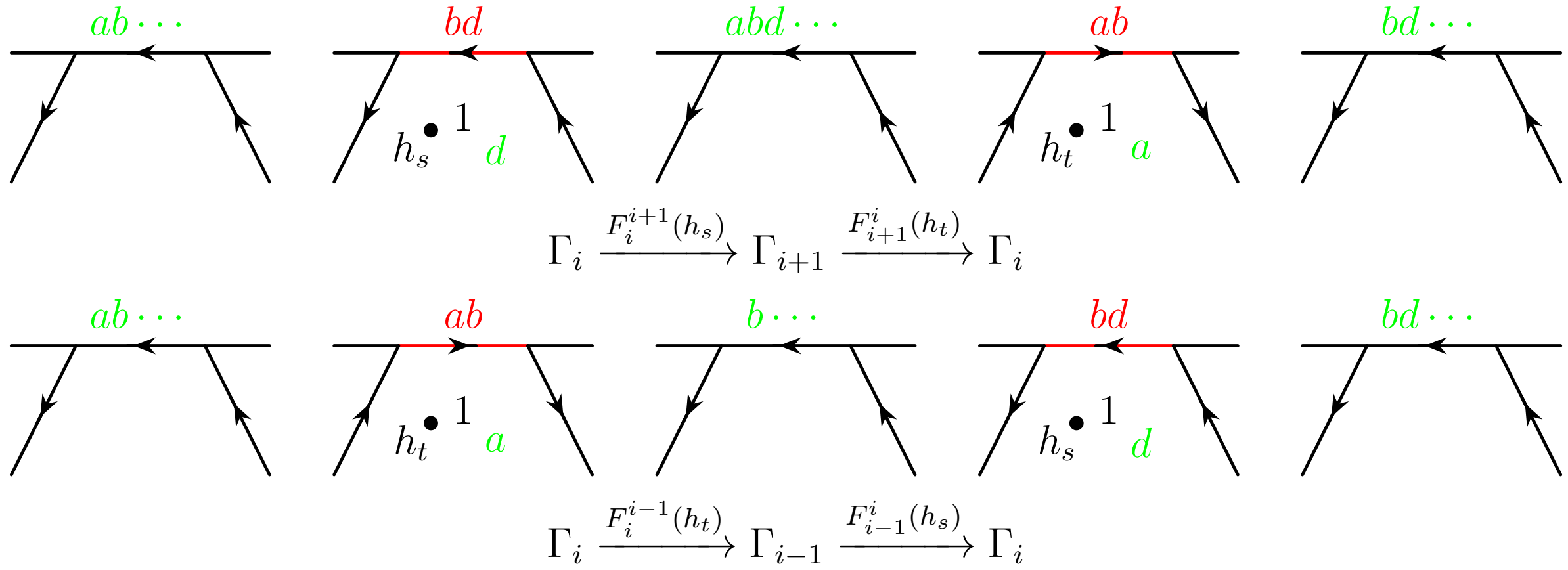

  % The two colorings in Figure~\ref{fig:2colors} are related by a Kempe move along a sphere, hence their contribution cancels out by \cite[Lemma 2.19]{RW2020} and its proof.
  
  Compare the two colored foams in Figure~\ref{fig:2colors}, as before, we can check that they have the same monochrome and bichrome surfaces. The only difference in the sign contribution comes from the pair $\{a,d\}$. The upper foam in Figure \ref{fig:2colors} has two additional circles (bounding the inserted facets), exactly one of which is positive. The two contributions therefore differ by a sign and cancel in the foam evaluation.
\end{proof}

The following special case plays an important role in what follows, therefore we state it separately.

\begin{corollary}\label{cor:null homotopy without overline}
 For any $0\leq i\leq l$, and $s \in \Z$, $r \geq s\geq 0$, we have: 
  \begin{equation}
    F_{i+1}^i(h_{r-s})F_i^{i+1}(h_s)+F_{i-1}^i(h_s)F_i^{i-1}(h_{r-s})=(-1)^{|B|+1}\idfoam_{\Gamma_i}.
  \label{null-homotopyFormula}
  \end{equation}
\end{corollary}

The next theorem unveils how the decorations on simple facet insertion morphisms add together to form facet insertion morphisms that need not be simple. By induction this leads to a more general addition formula. Alternatively, the next theorem can be understood as a decomposition of facet insertion morphisms into simple facet insertion morphisms.
\begin{theorem}\label{thm:composition of slice morphisms}
  For any $i,s,t \in \Z$, such that  $0\leq i\leq l$, and $s,t \geq 0$, we have:
  \begin{equation}
    F^{i\pm 2}_{i\pm 1}(h_t)F^{i\pm 1}_{i}(h_s)=\pm F^{i\pm 2}_{i}(\frac{x_1^sx_2^t-x_1^tx_2^s}{x_1-x_2}).
    \label{DifferentialFormulaOne}
  \end{equation}
\end{theorem}
\begin{proof}
For any coloring $c$ of $F_i^{i+2}(\frac{x_1^sx_2^t-x_1^tx_2^s}{x_1-x_2})$ as displayed in the top row of Figure \ref{fig:double1facet}, there are two colorings of $F^{i\pm 2}_{i\pm 1}(h_t)F^{i\pm 1}_{i}(h_s)$ corresponding to $c$. One of those, call it $c'$, is displayed in the bottom row Figure \ref{fig:double1facet}, the other one, call it $c''$, is obtained from $c'$ by swapping $a$ and $a'$ in said figure. 
  \begin{figure}[htbp]
    \centering 
    \begin{tikzpicture}[
  scale=0.60,
  line/.style={
    line width=0.9pt,
    line cap=round,
    line join=round
  },
  flow/.style={
    line,
    postaction={decorate},
    decoration={
      markings,
      mark=at position 0.53 with
      {\arrow{Stealth[length=2.2mm,width=1.7mm]}}
    }
  },
  every node/.style={
    font=\normalsize,
    inner sep=1.5pt
  }
    ]
    \begin{scope}
      %\node[green] at (2,2.5) {$ab\cdots$};
      \draw[line] (0,2)--(4,2);
      \draw[line] (0,0)--(1,2);
      \draw[line] (3,2)--(4,0);
      \draw[flow] (0.5,1)--(0.4,0.8);
      \draw[flow] (2,2)--(1.8,2);
      \draw[flow] (3.5,1)--(3.4,1.2);
    \end{scope}

    \begin{scope}[xshift=5cm]
      \node at (2,1) {$1$};
      \node[green] at (2.5,0.5) {$a$};
      \draw[line] (1,2)--(3,2);
      %\node[red] at (2,2.5) {$bd$}; 
      \draw[line] (0,2)--(1,2);
      \draw[line] (3,2)--(4,2);
      \draw[line] (0,0)--(1,2);
      \draw[line] (3,2)--(4,0);
      \draw[flow] (0.5,1)--(0.4,0.8);
      \draw[flow] (2,2)--(1.8,2);
      \draw[flow] (3.5,1)--(3.4,1.2);
      \draw[fill=black] (1.5,0.8) circle[radius=0.1];
      \node at (1.2,0.5) {$h_s$};
    \end{scope}

    \begin{scope}[xshift=10cm]
      %\node[green] at (2,2.5) {$abd\cdots$};
      \draw[line] (0,2)--(4,2);
      \draw[line] (0,0)--(1,2);
      \draw[line] (3,2)--(4,0);
      \draw[flow] (0.5,1)--(0.4,0.8);
      \draw[flow] (2,2)--(1.8,2);
      \draw[flow] (3.5,1)--(3.4,1.2);
    \end{scope}

    \begin{scope}[xshift=15cm]
      \node at (2,1) {$1$};
      \node[green] at (2.5,0.5) {$a'$};
      \draw[line] (1,2)--(3,2);
      %\node[red] at (2,2.5) {$ab$}; 
      \draw[line] (0,2)--(1,2);
      \draw[line] (3,2)--(4,2);
      \draw[line] (0,0)--(1,2);
      \draw[line] (3,2)--(4,0);
      \draw[flow] (0.5,1)--(0.4,0.8);
      \draw[flow] (2,2)--(1.8,2);
      \draw[flow] (3.5,1)--(3.4,1.2);
      \draw[fill=black] (1.5,0.8) circle[radius=0.1];
      \node at (1.2,0.5) {$h_t$};
    \end{scope}

    \begin{scope}[xshift=20cm]
      %\node[green] at (2,2.5) {$bd\cdots$};
      \draw[line] (0,2)--(4,2);
      \draw[line] (0,0)--(1,2);
      \draw[line] (3,2)--(4,0);
      \draw[flow] (0.5,1)--(0.4,0.8);
      \draw[flow] (2,2)--(1.8,2);
      \draw[flow] (3.5,1)--(3.4,1.2);
    \end{scope}
    \end{tikzpicture}
    $$\Gamma_i\xrightarrow{F_i^{i+1}(h_s)}\Gamma_{i+1}\xrightarrow{F_{i+1}^{i+2}(h_t)}\Gamma_{i+2}$$
    \begin{tikzpicture}[
  scale=0.60,
  line/.style={
    line width=0.9pt,
    line cap=round,
    line join=round
  },
  flow/.style={
    line,
    postaction={decorate},
    decoration={
      markings,
      mark=at position 0.53 with
      {\arrow{Stealth[length=2.2mm,width=1.7mm]}}
    }
  },
  every node/.style={
    font=\normalsize,
    inner sep=1.5pt
  }
    ]
    \begin{scope}
      %\node[green] at (2,2.5) {$ab\cdots$};
      \draw[line] (0,2)--(4,2);
      \draw[line] (0,0)--(1,2);
      \draw[line] (3,2)--(4,0);
      \draw[flow] (0.5,1)--(0.4,0.8);
      \draw[flow] (2,2)--(1.8,2);
      \draw[flow] (3.5,1)--(3.4,1.2);
    \end{scope}

    \begin{scope}[xshift=10cm]
      \node at (2.5,1.3) {$2$};
      \node[green] at (3.2,0) {$aa'$};
      \draw[line] (1,2)--(3,2);
      %\node[red] at (2,2.5) {$bd$}; 
      \draw[line] (0,2)--(1,2);
      \draw[line] (3,2)--(4,2);
      \draw[line] (0,0)--(1,2);
      \draw[line] (3,2)--(4,0);
      \draw[flow] (0.5,1)--(0.4,0.8);
      \draw[flow] (2,2)--(1.8,2);
      \draw[flow] (3.5,1)--(3.4,1.2);
      \draw[fill=black] (1.5,0.8) circle[radius=0.1];
      \node at (1.5,0.2) {$\frac{x_1^sx_2^t-x_1^tx_2^s}{x_1-x_2}$};
    \end{scope}

    \begin{scope}[xshift=20cm]
      %\node[green] at (2,2.5) {$bd\cdots$};
      \draw[line] (0,2)--(4,2);
      \draw[line] (0,0)--(1,2);
      \draw[line] (3,2)--(4,0);
      \draw[flow] (0.5,1)--(0.4,0.8);
      \draw[flow] (2,2)--(1.8,2);
      \draw[flow] (3.5,1)--(3.4,1.2);
    \end{scope}
    \end{tikzpicture}
    $$\Gamma_i\xrightarrow{F_i^{i+2}(\frac{x_1^sx_2^t-x_1^tx_2^s}{x_1-x_2})}\Gamma_{i+2}$$
    \caption{Composition of two simple facet insertion morphisms and a facet insertion morphism with a $2$-facet.}
    \label{fig:double1facet}
  \end{figure}
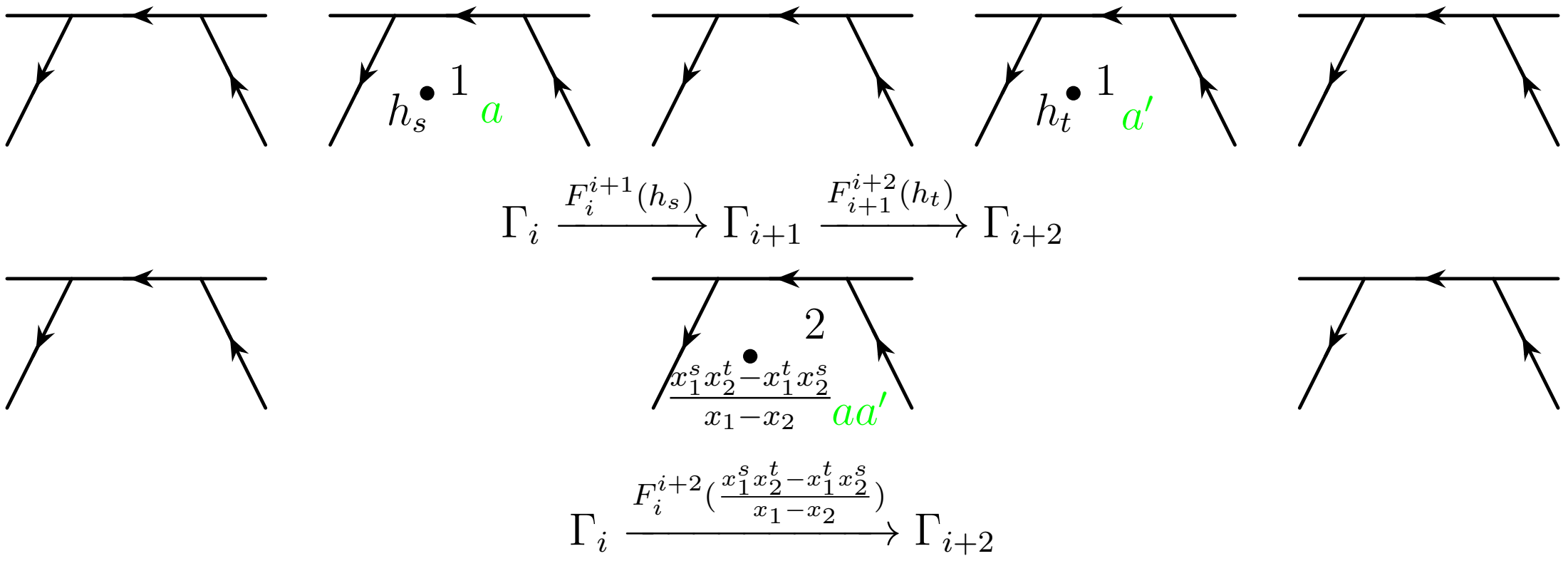
    In the foam evaluation, the only nontrivial contributions come from the bichrome surface associated with $a$ and $a'$ and an additional circle, which is positive precisely when $a>a'$. Taking into account the decorations we can compute: 
  \[ \frac{\langle F^{i\pm 2}_{i\pm 1}(h_t)F^{i\pm 1}_{i}(h_s),c' \rangle}{\langle F_i^{i+2}(\frac{x_1^sx_2^t-x_1^tx_2^s}{x_1-x_2}), c\rangle}  =(-1)^{\mathbbm{1}\{a>a'\}}\cdot \frac{1}{(-1)^{\mathbbm{1}\{a>a'\}}(x_a-x_{a'})}\cdot \frac{x_a^sx_{a'}^t(x_a-x_{a'})}{x_a^sx_{a'}^t-x_{a'}^sx_{a}^t}=\frac{x_a^sx_{a'}^t}{x_a^sx_{a'}^t-x_{a'}^sx_{a}^t}.\] 
  
  An analogous computation gives:
\[\frac{\langle F^{i\pm 2}_{i\pm 1}(h_t)F^{i\pm 1}_{i}(h_s),c'' \rangle}{\langle F_i^{i+2}(\frac{x_1^sx_2^t-x_1^tx_2^s}{x_1-x_2}), c\rangle}=\frac{x_{a'}^sx_{a}^t}{x_{a'}^sx_{a}^t-x_{a}^sx_{a'}^t}.\]
  
  Summing the contributions, we obtain: 
  \[\frac{x_a^sx_{a'}^t}{x_a^sx_{a'}^t-x_{a'}^sx_{a}^t}+\frac{x_{a'}^sx_{a}^t}{x_{a'}^sx_{a}^t-x_{a}^sx_{a'}^t}=1.\]
  This proves Equation~\ref{DifferentialFormulaOne} when the sign is positive. The other case can be proved in a similar way.
\end{proof}

Direct consequences of the above are two useful corollaries. The first one tells us how the composition of simple facet insertion morphisms works in the special cases of the simplest possible decoration pairs and the second one states the formula for the composition of a simple facet insertion foam with a generic facet insertion foam.

\begin{corollary}\label{cor:composition of slice morphisms}
  For any $0\leq i\leq l$,  $t,s \in \Z$, $t\geq0, s\geq 0$, we have: 
  \begin{equation}
    F^{i\pm 2}_{i\pm 1}(h_t)F^{i\pm 1}_{i}(h_s)=-F^{i\pm 2}_{i\pm 1}(h_s)F^{i\pm 1}_{i}(h_t), 
    \label{ExchangeDecorationFormula}
  \end{equation}
  \begin{equation}
    F^{i\pm 2}_{i\pm 1}(h_s)F^{i\pm 1}_{i}(h_s)=0,
    \label{differentialFormulaTwo}
  \end{equation}
\end{corollary}

\begin{corollary}\label{cor:iterated composition of slice morphisms}
      For any $0\leq i, k\leq l$, such that $i+k \leq l$ and $k \geq 1$,
  \begin{equation}
    F^{i\pm k}_{i\pm 1}( s_{\overline{\lambda}} )F^{i\pm 1}_{i}(h_{\lambda_1 + k - 1})= \pm F^{i\pm k}_{i}(s_{\lambda}).
    \label{DifferentialFormulaThree}
  \end{equation}

  Where $s_\lambda$ is a Schur polynomial associated to a Young tableau (in English notation) $\lambda$, $\overline{\lambda}$ denotes the Young tableau obtained from $\lambda$ by erasing the top row, and $\lambda_1$ denotes the size of the top row $\lambda$. See \cite[Section 1.3]{McD} for precise definitions.
\end{corollary}
\begin{proof}
    This follows from a simple induction using the result from Theorem \ref{thm:composition of slice morphisms}.
\end{proof}

The following theorem provides a categorification of the Kekulé relations in the homotopy category of $\Foam$.

\begin{theorem}\label{thm:hocat categorificaiton}
For any $s\in \{0,1,\cdots, r\}$, let $C_s$ be the bidirectional chain complex in the homotopy category of $\Foam$ 
\begin{equation}
  0\rightleftharpoons \Gamma_0\xrightleftharpoons[F_1^0(h_{r-s})]{F_0^1(h_s)}q^{r-2s}\Gamma_1\xrightleftharpoons[F_2^1(h_{r-s})]{F_1^2(h_s)} \cdots\xrightleftharpoons[F_{l}^{l-1}(h_{r-s})]{F_{l-1}^l(h_s)} q^{l(r-2s)}\Gamma_l\rightleftharpoons 0, 
  \label{LESsF}
\end{equation}
Then $C_s$'s are long exact sequences (when read in either direction), their Euler characteristics $\sum_{i=0}^l(-1)^iq^{i(r-2s)}\langle\Gamma_i\rangle$ are equal to $0$ and they recover the transformed Kekulé relations from Definition \ref{def:true kekule} and by Lemma \ref{lem:true kekule} also the Kekulé relations of Definition \ref{def:kekule}. 
\label{th:1stCat}
\end{theorem}
\begin{proof}
    Note that $C_s$'s are well defined by Equation~\ref{differentialFormulaTwo}. Moreover, Equation~\ref{null-homotopyFormula} says that the identity chain maps on $C_s$'s are homotopic to the zero maps, that is $C_s$'s are null-homotopic and therefore exact.

  Let $\gamma(t):=\sum_{i=0}^{l}(-1)^i\langle\Gamma_i\rangle t^i$ be a generating function in $t$ for $\langle\Gamma_i\rangle$'s (understood as formal symbols). The Euler characteristic of $C_s$ equals $\gamma(q^{r-2s})$ for any $s$. Since $C_s$'s are exact, their Euler characteristics equal $0$, therefore $q^{r-2s}$'s are roots of $\gamma$ (closed under inversion). It follows that:
  \begin{align*}
      \gamma(t)=G(t)\prod_{s=0}^{r}(t-q^{r-2s})=(\sum_{i=0}^{w}g_i t^i)\prod_{s=0}^{r}(t-q^{r-2s}),
  \end{align*}
  where $G(t)$ is a polynomial of degree $\deg G(t) = w\leq l-r-1$ with coefficients $g_i$. Hence,
  \[G(t)=\gamma(t)\prod_{s=0}^{r}\frac{1}{t-q^{r-2s}}=(\sum_{i=0}^{l}(-1)^i\langle\Gamma_i\rangle t^i)(\sum_{i=0}^{\infty}h_i(q^{-r},q^{2-r},\cdots,q^r)t^i).\]
  Recall that $h_i(q^{-r},q^{2-r},\cdots,q^r)=\begin{bmatrix}
    r+i\\
    r
  \end{bmatrix}_q$, which allows us to express the coefficients of $G(t)$ as 
  \[g_i=\sum_{j=0}^{i}(-1)^j\langle\Gamma_j\rangle h_{i-j}(q^{-r},q^{2-r},\cdots,q^r)=\sum_{j=0}^{i}(-1)^j\begin{bmatrix}
    r+i-j\\
    r
  \end{bmatrix}_q\langle\Gamma_j\rangle.\]
  As $g_i=0$ for $i\geq l-r$, the result follows from Lemma \ref{lem:true kekule}
\end{proof}

The above theorem justifies the definition of a third version of Kekulé relations.

\begin{definition}\label{def:simplified kekule}
    Fix $m \in \Z$, $m > 0$ and let $a,b \in \Z_{\geq 0}^m$ be two sequences of nonnegative integers. For $l \geq i \geq l - r$, we introduce a collection of web relations, called \emph{simplified Kekulé relations}
\begin{align}\label{eq:true kekuleThree}
0 = \sum_{i=0}^l(-1)^iq^{i(r-2s)}\langle\Gamma_i\rangle.
\end{align}
\end{definition}

\section{Categorification of the Kekulé relations in $\Kar(\Foam)$}\label{sec:kar categorification}

Recall that we fixed an arbitrary choice of $m$, $a$, and $b$ in the beginning of Section \ref{sec:hocat categorificaiton}.

\subsection{The orthogonal idempotents $\Li{s}{i}$ and $\Ri{s}{i}$}
We wish to exhibit explicit chain contractions for the bidirectional long exact sequence from Equation \eqref{LESsF}. We will modify the morphisms $F_{i+1}^i$ to $\overline{F_{i+1}^i}$, so that when the replacement is made in Equation \eqref{eq:Fii+1 Fi+1i}, the right-hand side vanishes for $s+t\neq r$.

\begin{definition}\label{def:Li,Ri}
    For any $s\in \{0,1,\cdots, r\}$ and $i \in \{0,1, \ldots l\}$, let
    \begin{align*}
        \overline{F_{i+1}^i} (h_s)= (-1)^{|B| + 1}\sum_{v = 0}^s(-1)^v E_v(x_{\overline{B}}) F_{i+1}^i(h_{s-v}).
    \end{align*}
    % We call this \emph{the weighted facet insertion morphism (foam)}. 
    Additionally we introduce the following (which are idempotents according to Lemma \ref{lem:Li,Ri prop}):
    \begin{align*}
        \Li{s}{i} &:= F_{i-1}^i(h_{r-s} )\overline{F_{i}^{i-1}}(h_s), \\
        \Ri{s}{i} &:= \overline{F_{i+1}^i} (h_s)F_{i}^{i+1} (h_{r-s}).
    \end{align*}
\end{definition}

\begin{remark}
    Instead of modifying $F_{i+1}^i(h_s)$ to $\overline{F_{i+1}^i} (h_s)$, we could modify $F_{i-1}^i(h_{s} )$, or both $F_{i-1}^i(h_{s} )$ and $F_{i+1}^i(h_s)$. Any of these modifications would lead to an analogous theory.
\end{remark}

Replacing $F_{i+1}^i$ with $\overline{F_{i+1}^i}$ simplifies some of the results from the previous section. Notably, we have the following two lemmas.

\begin{lemma}\label{lem:Li,Ri prop}
    For any $s,t,y \in \{0,1,\cdots, r\}$, $s \neq t, s+y \neq r$ and any $i-1, i, i+1, i+2 \in \{0,1, \ldots l\}$, it holds that
\begin{align}
    \overline{F_{i+1}^i} (h_s) \overline{F_{i+2}^{i+1}} (h_s) & = 0, \\
     \overline{F_{i+1}^i} (h_s) \overline{F_{i+2}^{i+1}} (h_t) & = -\overline{F_{i+1}^i} (h_t) \overline{F_{i+2}^{i+1}} (h_s), \\
    F_{i-1}^i (h_s) \overline{F_{i}^{i-1}} (h_y) & = -\overline{F_{i+1}^{i}} (h_y)  F_{i}^{i+1} (h_s) , \label{eq:main L,R the rest vanish} \\
     \Li{s}{i} + \Ri{s}{i} &= \idfoam_{\Gamma_i},\label{eq:main L,R equation}  \\
     \Li{s}{i} \cdot \Li{t}{i} &= \Li{t}{i} \cdot \Li{s}{i}, \\
     \Ri{s}{i} \cdot \Ri{t}{i} &= \Ri{t}{i} \cdot \Ri{s}{i}.
\end{align}
Moreover, for each $i$, $\Li{s}{i},$ and $\Ri{s}{i}$ are mutually orthogonal idempotents of degree $0$. 
\end{lemma}
\begin{proof}
    The degree result follows readily from Lemma \ref{DegreeLem}. The rest of the statements follow by a direct computation using the results from Theorems \ref{Null-homotopyTh}, \ref{thm:composition of slice morphisms}, and Corollary \ref{cor:composition of slice morphisms}. We will work out \eqref{eq:main L,R equation} and \eqref{eq:main L,R the rest vanish} in detail and omit the rest, as they follow using similar techniques.

    First, we show that Equation \eqref{eq:main L,R equation} holds, by expanding and applying  Theorem \ref{Null-homotopyTh}:
    \begin{align*}
         \Li{s}{i} + \Ri{s}{i} &= F_{i-1}^i(h_{r-s} )\overline{F_{i}^{i-1}}(h_s) + \overline{F_{i+1}^i} (h_s)F_{i}^{i+1} (h_{r-s}) \\
         &= (-1)^{|B| + 1}\sum_{v = 0}^s(-1)^v E_v(x_{\overline{B}}) \left(F_{i-1}^i(h_{r-s} )F_{i}^{i-1}(h_{s-v}) + F_{i+1}^i(h_{s-v})F_{i}^{i+1} (h_{r-s}) \right) \\
         &= (-1)^{2(|B| + 1)} \idfoam_{\Gamma_i} = \idfoam_{\Gamma_i}.
    \end{align*}

    Next, we focus on Equation \eqref{eq:main L,R the rest vanish}:
    \begin{align*}
         F_{i-1}^i (h_s) \overline{F_{i}^{i-1}} (h_y) &+\overline{F_{i+1}^{i}} (h_y)  F_{i}^{i+1} (h_s) \\
         &=(-1)^{|B| + 1} \sum_{v = 0}^y(-1)^v E_v(x_{\overline{B}}) \left(F_{i-1}^i(h_{s} )F_{i}^{i-1}(h_{y-v}) + F_{i+1}^i(h_{y-v})F_{i}^{i+1} (h_{s}) \right) \\
         &= \sum_{v = 0}^y(-1)^v E_v(x_{\overline{B}}) H_{s + y - v - r}(x_{\overline{B}}) \\
         &= 0.
    \end{align*}
The last equality holds because $y+s \neq r$ by Definition \ref{def:e,h on foam}.
\end{proof}

It is important to note that Equation \eqref{eq:main L,R equation} will play the central role in the rest of the paper. All of the foam relations are going to be essentially deduced from this equation. See for instance Proposition \ref{prop:exact-word-projectors}, Theorem \ref{thm:upper-factor-through}, and Lemma \ref{lem:six vertex categorification}. Naive attempts to obtain analogues of said equation using $F_i^{i+k}(h_s)$ instead of $F_i^{i+1}(h_s)$ in Definition \ref{def:Li,Ri} and Corollary \ref{cor:null homotopy without overline} fail due to the decomposition given in Corollary \ref{cor:iterated composition of slice morphisms}. A more detailed study is warranted to determine the entire structure of the ($0$-degree part of the) endomorphism algebra of $\Gamma_i$. 

\begin{lemma}\label{lem:prod of Ri Li vanishing}
    Let $s_0, s_1, \ldots s_i \in \{0,1,\cdots, r\}$ be pairwise distinct. Then 
    \begin{align*}
        \prod_{k = 0}^i \Li{s_k}{i} = 0.
    \end{align*}
    Similarly, for pairwise distinct $t_0, t_1, \ldots t_{l-i} \in \{0,1,\cdots, r\}$,
    \begin{align*}
        \prod_{k = 0}^{l-i} \Ri{t_k}{i} = 0.
    \end{align*}
\end{lemma}
\begin{proof}
    We will prove the first equality. By repeatedly applying Lemma \ref{lem:Li,Ri prop}, we can expand (we understand $\prod_{i = a}^b c_i = c_b c_{b-1}\cdots c_a$):
    \begin{align*}
        \prod_{k = 0}^i \Li{s_k}{i} & = \prod_{k = 0}^i F_{-1 + k}^{k} (h_{k}) \prod_{l = 0}^i \overline{F_{i- l }^{i-l-1}}(h_l),
    \end{align*}
    for some decorations $(h_k)_k, (h_l)_l$ on the inserted facets. Hence this map factors through $\Gamma_{-1} = 0$ and the entire product vanishes.
\end{proof}

\subsection{Corners and isomorphisms}\label{sec:corner}

We will first define corner algebras in our context. The original definition comes from matrix algebras and noncommutative algebras, see \cite[Section 2.7]{Bresar2014}. Using morphisms between different corners will lead us to a categorification of the Kekulé relations in $\Kar(\Foam)$.

\begin{definition}\label{def:corner algebra}
    For a degree zero idempotent $e \in\Hom_{\Foam^\varepsilon}(\Gamma_i,\Gamma_i)$, its associated \emph{corner algebra} is
\[
 e\Hom_{\Foam^\varepsilon}(\Gamma_i,\Gamma_i)e.
\]
Its degree $0$ part is the endomorphism algebra of
$\mathcal{E}(e)$ in the Karoubi envelope. The grading shifts from Equation~\eqref{LESsF} will be suppressed in the corner algebra, as formal shifts do not change the algebra structure.
\end{definition}

The following definition introduces a compact notation for products of the $\Li{s}{i}$ and $\Ri{s}{i}$ idempotents. Their properties are summarized in Proposition \ref{prop:exact-word-projectors}. 

\begin{definition}\label{def:def Eis}
    For $0\leq i\leq l$ and $S\subseteq\{0,\ldots,r\}$, define
\begin{align}\label{eq:exact-word-projector}
 \mathfrak{E}_{i,S}:=
 \left(\prod_{s\in S}\Li{s}{i}\right)
 \left(\prod_{t\in\{0,\ldots,r\}- S}\Ri{t}{i}\right).
\end{align}
\end{definition}

Note that many of the $\mathfrak{E}_{i,S}$ are $0$ due to Lemma \ref{lem:prod of Ri Li vanishing}. The following proposition gives a decomposition of the identity foam. It is not a priori clear how this decomposition recovers the Kekulé relations. Section \ref{sec:iterated-transfer} will reinterpret this equality and make the connection clear. We will accomplish this by ``transferring" each $\mathcal{E}(\mathfrak{E}_{i,S})$ to some $\mathcal{E}(\mathfrak{E}_{i',\emptyset})$ (up to a grading shift) using products of facet insertion morphisms (see Lemma \ref{lem:choice-transfer}).

\begin{proposition}\label{prop:exact-word-projectors}
For $0\leq i\leq l$, the set $(\mathfrak{E}_{i,S})_{S\subseteq \{0,\ldots,r\}}$ is a complete set of idempotents of $\Gamma_i$, all of which are mutually orthogonal and
\begin{equation}\label{eq:exact-word-sum}
 \idfoam_{\Gamma_i}=\sum_{S\subseteq\{0,\ldots,r\}}\mathfrak{E}_{i,S}.
\end{equation}
The idempotent $\mathfrak{E}_{i,S}$ is zero when
\begin{equation}\label{eq:exact-word-cardinality}
 |S| < i+r+1-l \text{ or } |S| > i .
\end{equation}
Consequently, in the Karoubi envelope the web $\Gamma_i$ decomposes as
\begin{equation}\label{eq:web-karoubi-decomposition}
 \mathcal{E}(\Gamma_i)\cong
 \osum_{\substack{S\subseteq\{0,\ldots,r\}\\\mathfrak{E}_{i,S}\neq0}}
 \mathcal{E}(\mathfrak{E}_{i,S}),
\end{equation}
\end{proposition}

\begin{proof}
It is straightforward to check that each of the $\mathfrak{E}_{i,S}$ is idempotent from the definition. If $S\neq T$, then there exists some $s$ in the symmetric difference of $S$ and $T$, for which $\Li{s}{i} \cdot \Ri{s}{i}=0$ or  $0 =  \Ri{s}{i}\cdot \Li{s}{i} \cdot $ appear in $\mathfrak{E}_{i,S} \mathfrak{E}_{i, T}$. So $\mathfrak{E}_{i,S}$'s are mutually orthogonal.

The decomposition of $\idfoam_{\Gamma_i}$ follows from Lemma~\ref{lem:Li,Ri prop},
\[
 \prod_{s=0}^{r}\bigl(\Li{s}{i}+\Ri{s}{i}\bigr)
 =\sum_{S\subseteq\{0,\ldots,r\}}\mathfrak{E}_{i,S}
 =\idfoam_{\Gamma_i}.
\]

If $|S|\geq i+1$, the $\overline L$-part of $\mathfrak{E}_{i,S}$ contains a product of $i+1$ distinct factors and vanishes by Lemma \ref{lem:prod of Ri Li vanishing}. Similarly, if $r+1-|S|\geq l-i+1$, the $\overline R$-part of $\mathfrak{E}_{i,S}$ vanishes. This together with $0\leq|S|\leq r+1$ proves Equation~\eqref{eq:exact-word-cardinality}. Similar arguments show that $(\mathfrak{E}_{i,S})_{S\subseteq \{0,\ldots,r\}}$ is a complete set of idempotents of $\Gamma_i$ and therefore Equation~\eqref{eq:web-karoubi-decomposition} holds.
\end{proof}

The next proposition gives isomorphisms between the corner algebras of $\Gamma_{i-1}$ and $\Gamma_i$. We focus on the simplest case of using just $\Li{s}{i}$ or $\Ri{s}{i}$ for some $0\leq s\leq r$, as that suffices for the proof of Lemma \ref{lem:choice-transfer} and its applications. A more general statement could be made by using more complicated $\mathfrak{E}_{i,S}$'s by extending the arguments with induction.

\begin{proposition}\label{prop:corner-transfer}

For $1\leq i\leq l$ and $0\leq s\leq r$, there is an 
algebra isomorphism
\[
 \omega: \Ri{s}{i-1}\Hom_{\Foam^\varepsilon}(\Gamma_{i-1},\Gamma_{i-1}) \Ri{s}{i-1}
 \ \xrightarrow{\ \sim\ }\
 \Li{s}{i}\Hom_{\Foam^\varepsilon}(\Gamma_i,\Gamma_i)\Li{s}{i},
\]
given by
\begin{equation}\label{eq:corner-transfer}
 \omega: e\longmapsto   F_{i-1}^{i}(h_{r-s})\,e\,\overline{F_i^{i-1}}(h_s),
\end{equation}
with its inverse given as
\begin{equation}\label{eq:corner-transfer-inverse}
 \omega^{-1}:e\longmapsto  \overline{F_i^{i-1}}(h_s)\,e\,F_{i-1}^{i}(h_{r-s}).
\end{equation}
We again suppress grading shifts from Equation~\eqref{LESsF}.
\end{proposition}

\begin{proof}
Since $\Ri{s}{i-1} =  \overline{F_i^{i-1}}(h_s)F_{i-1}^{i}(h_{r-s})$ and $\Li{s}{i} = F_{i-1}^{i}(h_{r-s}) \overline{F_i^{i-1}}(h_s)$ it holds that $\omega$ and $\omega^{-1}$ are indeed inverses by \cite[Exercise 3.4]{Bresar2014}. 
\end{proof}

The following lemma is an auxiliary computation that will be used in the proof of Proposition \ref{prop:iterated-transfer}. The main conclusions are Equations \eqref{eq:exact-word-lower-transfer} and \eqref{eq:exact-word-transfer}, which relate 

\begin{lemma}\label{lem:choice-transfer}
For $1\leq i\leq l$ and $0\leq s,t\leq r$ with $t\neq s$, the
isomorphism from Proposition~\ref{prop:corner-transfer} acts as:
\begin{align}
 &F_{i-1}^{i}(h_{r-s})  \bigl(\Ri{s}{i-1}\cdot \Li{t}{i-1} \cdot\Ri{s}{i-1}\bigr)  \overline{F_i^{i-1}}(h_s)  =\Li{s}{i} \cdot\Li{t}{i},
 \label{eq:L-choice-transfer} \\
 &F_{i-1}^{i}(h_{r-s})  \bigl(\Ri{s}{i-1}\cdot \Ri{t}{i-1}\cdot \Ri{s}{i-1}\bigr)  \overline{F_i^{i-1}}(h_s)  =\Li{s}{i} \cdot \Ri{t}{i},
 \label{eq:R-choice-transfer} \\
  &\overline{F_i^{i-1}}(h_s)  \bigl(\Li{s}{i}\cdot \Li{t}{i}\cdot\Li{s}{i}\bigr)  F_{i-1}^{i}(h_{r-s})  =\Ri{s}{i-1}\cdot\Li{t}{i-1},
 \label{eq:L-choice-lower-transfer}\\
 &\overline{F_i^{i-1}}(h_s)  \bigl(\Li{s}{i}\cdot \Ri{t}{i}\cdot\Li{s}{i}\bigr)  F_{i-1}^{i}(h_{r-s})  =\Ri{s}{i-1}\cdot\Ri{t}{i-1}.
 \label{eq:R-choice-lower-transfer}
\end{align}
Consequently, if $S\subseteq\{0,\ldots,r\}$ and $s\notin S$, then
\begin{align}
 F_{i-1}^{i}(h_{r-s}) \cdot \mathfrak{E}_{i-1, S}\cdot \overline{F_i^{i-1}}(h_s)&=\mathfrak{E}_{i, S\cup\{s\}},\label{eq:exact-word-transfer} \\
  \overline{F_i^{i-1}}(h_s)\cdot \mathfrak{E}_{i, S\cup\{s\}}  \cdot  F_{i-1}^{i}(h_{r-s})&= \mathfrak{E}_{i-1, S}.\label{eq:exact-word-lower-transfer}
\end{align}
\end{lemma}

\begin{proof}
We prove \eqref{eq:L-choice-transfer}. The other identities follow analogously. Using Lemma~\ref{lem:Li,Ri prop} compute that:
\begin{align*}
 &F_{i-1}^{i}(h_{r-s}) \left( \Li{t}{i-1} \right) \overline{F_i^{i-1}}(h_s)\\
 &\quad=
 F_{i-1}^{i}(h_{r-s})F_{i-2}^{i-1}(h_{r-t})
 \overline{F_{i-1}^{i-2}}(h_t)\overline{F_i^{i-1}}(h_s)\\
 &\quad=-
 F_{i-1}^{i}(h_{r-s})F_{i-2}^{i-1}(h_{r-t})
 \overline{F_{i-1}^{i-2}}(h_s)\overline{F_i^{i-1}}(h_t)\\
 &\quad=
 F_{i-1}^{i}(h_{r-s})\overline{F_i^{i-1}}(h_s)
 F_{i-1}^{i}(h_{r-t})\overline{F_i^{i-1}}(h_t)
 =\Li{s}{i} \cdot \Li{t}{i}.
\end{align*}
Expanding $\Ri{s}{i-1}$ and $\Ri{s}{i-1}$ completes the proof of \eqref{eq:L-choice-transfer} since $\Li{s}{i}$ and $\Li{t}{i}$ are commutative idempotents. Explicitly, for every
$x\in\Hom_{\Foam^\varepsilon}(\Gamma_{i-1},\Gamma_{i-1})$, we have:
\[
 F_{i-1}^{i}(h_{r-s})
 \bigl(\Ri{s}{i-1}x\Ri{s}{i-1}\bigr)
 \overline{F_i^{i-1}}(h_s)
 =
 \Li{s}{i}
 \bigl(F_{i-1}^{i}(h_{r-s})x\overline{F_i^{i-1}}(h_s)\bigr)
 \Li{s}{i}.
\]
\end{proof}

\subsection{Iterating the corner isomorphisms}\label{sec:iterated-transfer}

For an idempotent $\mathfrak{E}_{i,S}$ and $s\in S$, Equation~\eqref{eq:exact-word-transfer} expresses $\mathfrak{E}_{i,S}$ using $\mathfrak{E}_{i-1,S-\{s\}}$. Repeating this once for each element of $S$ expresses $\mathfrak{E}_{i,S}$ using $\mathfrak{E}_{i-|S|,\emptyset}$. The following proposition presents the details of this idea. 

\begin{proposition}\label{prop:iterated-transfer}
Let $0\leq i\leq l$ and $S\subseteq\{0,\ldots,r\}$. Put
$k=|S|\leq i$, choose an ordering $s_1,\ldots,s_k$ of the elements of
$S$, and set
\[
 j:=i-k,
 \qquad
 d(S):=\sum_{v=1}^{k}(2s_v-r).
\]
Assign $s_v$ to the decoration of the maps between $\Gamma_{j+v-1}$ and $\Gamma_{j+v}$. Define:
\begin{align}
 \overline{F_{i,S}}
 &:={\overline{F_{j+1}^{j}}}(h_{s_1})
    {\overline{F_{j+2}^{j+1}}}(h_{s_2})\cdots
    {\overline{F_i^{i-1}}}(h_{s_k}),
 \label{eq:iterated-lowering}\\
 F_{i,S}
 &:=F_{i-1}^{i}(h_{r-s_k})\cdots
    F_{j+1}^{j+2}(h_{r-s_2})
    F_j^{j+1}(h_{r-s_1}).
 \label{eq:iterated-raising}
\end{align}
For $S=\emptyset$, we understand that $d(S)=0$ and $F_{i,\emptyset}=\overline{F_{i,\emptyset}}= \idfoam_{\Gamma_i}$. Then
\begin{align}
 \mathfrak{E}_{i,S}=F_{i,S}\,\mathfrak{E}_{j,\emptyset}\,\overline{F_{i,S}}, \label{eq:iterated-exact-word}\\
  \mathfrak{E}_{j,\emptyset}=\overline{F_{i,S}}\,\mathfrak{E}_{i,S}\,F_{i,S}. \label{eq:iterated-exact-word-lower-transfer}
\end{align}
\end{proposition}

\begin{proof}
We will only prove \eqref{eq:iterated-exact-word}. The other equality follows similarly. Begin by inducting on $k$. The base case $k=0$ is clear. For $k > 0$, set $S'=S-\{s_k\}$. Equation~\eqref{eq:exact-word-transfer} gives:
\[
 \mathfrak{E}_{i,S}=F_{i-1}^{i}(h_{r-s_k})\mathfrak{E}_{i-1,S'} \overline{F_i^{i-1}}(h_{s_k}).
\]
Using the induction hypothesis we can substitute $\mathfrak{E}_{i-1,S'}$. The adjacent maps concatenate as in Equations~\eqref{eq:iterated-lowering} and \eqref{eq:iterated-raising}. This proves Equation~\eqref{eq:iterated-exact-word}.
\end{proof}

Note that reordering the elements of $S$ may change the sign of $\overline{F_{i,S}}$, while $\mathfrak{E}_{i,S}$ only depends on $i$ and $S$. 

\begin{theorem}\label{thm:upper-factor-through}
Using notation from Proposition \ref{prop:iterated-transfer}; for every $0\leq i\leq l$, the following equation is a decomposition into mutually orthogonal idempotents:
\begin{equation}\label{eq:upper-factor-through}
 \idfoam_{\Gamma_i}=
 \sum_{\substack{S\subseteq\{0,\ldots,r\}\\
         \max(0,i+r+1-l)\leq |S|\leq\min(r+1,i)}}
 F_{i,S}\mathfrak{E}_{i-|S|,\emptyset}\overline{F_{i,S}}.
\end{equation}
\end{theorem}

\begin{proof}
This is a direct consequence of Proposition~\ref{prop:iterated-transfer} and Proposition~\ref{prop:exact-word-projectors}.
\end{proof}

\begin{remark}
    The decomposition of $\idfoam_{\Gamma_i}$ in Theorem \ref{thm:upper-factor-through} and the preceding statements has been constructed so that all of the $\Li{s}{i}$'s are factored out of $\mathfrak{E}_{i,S}$, transferring $\mathfrak{E}_{i,S}$ (downwards) to $\mathfrak{E}_{i- |S|, \emptyset} \in \Hom_{\Foam^\varepsilon}(\Gamma_{i- |S|}, \Gamma_{i- |S|})$. Alternatively, we could factor all of the $\Ri{s}{i}$'s from $\mathfrak{E}_{i,S}$, transferring them (upwards) to $\mathfrak{E}_{i + r +1 - |S|, \{0,\ldots,r\}} \in \Hom_{\Foam^\varepsilon}(\Gamma_{i + r +1 - |S|}, \Gamma_{i + r +1 - |S|})$. This would lead to an analogue of Corollary \ref{cor:equiv to kekule in kar}.
\end{remark}

\subsection{Decategorification}\label{sec:incidence}

In this section we show that Theorem \ref{thm:upper-factor-through}  categorifies equations from Equation \eqref{eq:upper-factor-through}, how this yields decompositions in the Karoubi envelope, and how they relate to the Kekulé relations.

% First, we introduce some notation for the (split) Grothendieck ring of a graded Karoubi envelope. Let $K=\Q(q)$ and let
% \[
%  \mathcal{G}:=K\otimes_{\Z[q,q^{-1}]}K_0\bigl(\Kar(\Foam^\varepsilon)\bigr),
% \]
% elements of $\mathcal{G}$ will be denoted with square brackets $[\_]$, and $[W\{d\}]=q^d[W]$. 

\begin{theorem}\label{thm:incidence-convolution}
For $0\leq i\leq l$, $d(S)=\sum_{s_v \in S}(2s_v-r)$, there is an isomorphism 
\begin{equation}\label{eq:residual-web-decomposition}
 \mathcal E(\Gamma_i)\cong
 \bigoplus_{u=0}^{l-r-1}
 \ \bigoplus_{\substack{S\subseteq\{0,\ldots,r\}\\|S|=i-u}}
 q^{-d(S)}\mathcal E(\mathfrak{E}_{u,\emptyset}).
\end{equation}
Equivalently,
\begin{equation}\label{eq:incidence-convolution}
 \mathcal E(\Gamma_i) \cong
 \osum_{u=0}^{l-r-1}
 \begin{bmatrix}r+1\\i-u\end{bmatrix}_q \mathcal E(\mathfrak{E}_{u,\emptyset}),
 \qquad 0\leq i\leq l.
\end{equation}
Equation \eqref{eq:residual-web-decomposition} is categorified by Theorem \ref{thm:upper-factor-through}. 
\end{theorem}

\begin{proof}
The result can be deduced from Theorem \ref{thm:upper-factor-through}. To do so, note that Proposition \ref{prop:iterated-transfer} shows
\begin{align*}
     \mathfrak{E}_{i,S}= \mathfrak{E}_{i,S}F_{u,S}\mathfrak{E}_{u,\emptyset}\overline{F_{u,S}}\mathfrak{E}_{i,S}, \\\
  \mathfrak{E}_{u,\emptyset}=\mathfrak{E}_{u,\emptyset}\overline{F_{u,S}}\mathfrak{E}_{i,S}F_{u,S}\mathfrak{E}_{u,\emptyset}. 
\end{align*}
Therefore $ \mathfrak{E}_{u,\emptyset} \overline{F_{u,S}}  \mathfrak{E}_{i,S}$ and $\mathfrak{E}_{i,S} F_{u,S} \mathfrak{E}_{u,\emptyset}$ are mutually inverse in $\Kar(\Foam^\varepsilon).$ A degree computation shows that:
\begin{align}\label{eq:isomorphism for fixed S+i}
    \mathcal E(\mathfrak{E}_{i,S}) \cong q^{-d(S)} \mathcal E(\mathfrak{E}_{u ,\emptyset}).
\end{align}
This proves Equation \eqref{eq:residual-web-decomposition}.

Equation~\eqref{eq:incidence-convolution} follows from the quantum binomial theorem
\begin{equation}\label{eq:balanced-q-binomial}
 \prod_{s=0}^{r}\left(1+q^{2s-r}z\right)  =\sum_{k=0}^{r+1}  \begin{bmatrix}r+1\\k\end{bmatrix}_qz^k,
\end{equation}
since the coefficient of $z^k$ on the left is $\sum_{|S|=k}q^{d(S)}$.
\end{proof}

Next, we invert Equation \eqref{eq:incidence-convolution} to express $\mathcal E(\mathfrak{E}_{u,\emptyset})$ in terms of the objects $\mathcal E(\Gamma_i)$. Then, since $\mathfrak{E}_{u,\emptyset} = 0$ for $u \geq l-r$, we obtain relations between $\mathcal E(\Gamma_i)$'s and show that these are indeed the transformed Kekulé relations of Definition \ref{def:true kekule}.

\begin{proposition}\label{lem:karoubi categorification}
For every $0\leq j\leq l$, there is an isomorphism
\begin{equation}\label{eq:inverse-incidence}
 \mathcal E(\mathfrak{E}_{j,\emptyset}) \oplus \osum_{\substack{0\leq k\leq j\\ k\text{ odd}}} \begin{bmatrix}r+k\\k\end{bmatrix}_q \mathcal E(\Gamma_{j-k})
 \cong \osum_{\substack{0\leq k\leq j\\ k\text{ even}}}  \begin{bmatrix}r+k\\k\end{bmatrix}_q \mathcal E(\Gamma_{j-k}).
\end{equation}
Here $\mathcal E(\mathfrak{E}_{j,\emptyset})=0$ whenever $j\geq l-r$.
\end{proposition}

\begin{proof}
Set
\begin{align*}
 &c_a:=\begin{bmatrix}r+1\\a\end{bmatrix}_q, \quad
 d_k:=\begin{bmatrix}r+k\\k\end{bmatrix}_q, \\
 &C(z) := \sum_{a=0}^{r+1}c_a z^a, \quad
 D_{\mathrm{ev}}(z) := \sum_{\substack{k\geq0\\k\text{ even}}}d_k z^k, \quad
 D_{\mathrm{odd}}(z) := \sum_{\substack{k\geq0\\k\text{ odd}}} d_k z^k.
\end{align*}
We can invert $C(z)$ using the quantum binomial theorem (see Equation~\eqref{eq:balanced-q-binomial}) and the fact that $h_i(q^{-r},q^{2-r},\cdots,q^r)=\begin{bmatrix}
    r+i\\
    r
  \end{bmatrix}_q$, to obtain:
\begin{equation*}\label{eq:C-inverse}
 C(z)^{-1}=\sum_{k\geq0}(-1)^k \begin{bmatrix}r+k\\k\end{bmatrix}_qz^k = D_{\mathrm{ev}}(z)-D_{\mathrm{odd}}(z).
\end{equation*}
For a power series $p$, denote by $p(x)|_{x^i}$ the coefficient of $p(x)$ at $x^i$. Taking the coefficients of $z^{j-u}$ in $C(z)D_{\mathrm{ev}}(z)=1+C(z)D_{\mathrm{odd}}(z),$ and writing $\delta_{ju}$ for the Kronecker delta function gives:
\begin{align*}
\delta_{ju} &= \left. C(z)\bigl(D_{\mathrm{ev}}(z)- D_{\mathrm{odd}}(z)\bigr) \right|_{z^{j-u}} \\
&= \left.C(z)D_{\mathrm{ev}}(z)\right|_{z^{j-u}} - \left.C(z)D_{\mathrm{odd}}(z)\right|_{z^{j-u}} \\
&= \sum_{\substack{k=0\\k\text{ even}}}^{j} d_kc_{j-k-u} - \sum_{\substack{k=0\\k\text{ odd}}}^{j} d_kc_{j-k-u},
\end{align*}
which allows us to compute using Equation~\eqref{eq:incidence-convolution}:
\begin{align*}
\osum_{\substack{0\leq k\leq j\\k\text{ even}}} d_k\mathcal E(\Gamma_{j-k}) 
 &\cong \osum_{u=0}^{l-r-1} \left( \sum_{\substack{0\leq k\leq j\\k\text{ even}}} d_kc_{j-k-u} \right)\mathcal E(\mathfrak{E}_{u,\emptyset}) \\
&\cong \osum_{u=0}^{l-r-1} \left( \delta_{ju} +  \sum_{\substack{0\leq k\leq j\\k\text{ odd}}} d_kc_{j-k-u} \right)\mathcal E(\mathfrak{E}_{u,\emptyset}) \\
&\cong \mathcal E(\mathfrak{E}_{j,\emptyset}) \oplus \osum_{\substack{0\leq k\leq j\\k\text{ odd}}} d_k\mathcal E(\Gamma_{j-k}).
\end{align*}
Finally, Proposition \ref{prop:exact-word-projectors} gives $\mathfrak{E}_{j,\emptyset}=0$ when $j\geq l-r$.
\end{proof}

Next we establish that Theorem \ref{thm:upper-factor-through} categorifies precisely the Kekulé relations. Analyzing the proof shows that Proposition \ref{lem:karoubi categorification} presents the same web relations as Theorem \ref{thm:hocat categorificaiton} on the nose (after passing from the ambient categories in those statements to $\Web$).

\begin{corollary}\label{cor:equiv to kekule in kar}
    For $l - r \leq j \leq l$, the set of Equations \eqref{eq:inverse-incidence}:
    \begin{align*}
   \osum_{\substack{0\leq k\leq j\\ k\text{ odd}}} \begin{bmatrix}r+k\\k\end{bmatrix}_q \mathcal E(\Gamma_{j-k})
 \cong \osum_{\substack{0\leq k\leq j\\ k\text{ even}}}  \begin{bmatrix}r+k\\k\end{bmatrix}_q \mathcal E(\Gamma_{j-k}),
\end{align*}
recover the transformed Kekulé relations from Definition \ref{def:true kekule}.
\end{corollary}
\begin{proof}
    This follows by noting that $\mathcal E(\mathfrak{E}_{j,\emptyset})=0$ for the permissed choices of $j$ and considering $\mathcal E(\Gamma_{j-k})$ as $ \langle \Gamma_{j-k} \rangle$ and noting that $\sum_{k=0}^{j}(-1)^k 
    \begin{bmatrix} r+k\\k \end{bmatrix}_q \langle\Gamma_{j-k}\rangle$  are the transformed Kekulé relations.
\end{proof}
\begin{remark}
Let $S=\{0,1,\cdots,r\}$. The methods used in proving Proposition~\ref{thm:incidence-convolution} can be used to construct an isomorphism in the opposite direction, demonstrating 
    \[\mathcal E(\mathfrak{E}_{j,S}) \oplus \osum_{\substack{0\leq k\leq l-j\\ k\text{ odd}}} \begin{bmatrix}r+k\\k\end{bmatrix}_q \mathcal E(\Gamma_{j+k})
 \cong \osum_{\substack{0\leq k\leq l-j\\ k\text{ even}}}  \begin{bmatrix}r+k\\k\end{bmatrix}_q \mathcal E(\Gamma_{j+k}).\]
This equation, together with Proposition~\ref{thm:incidence-convolution}, and the isomorphism $\mathcal E(\mathfrak{E}_{j,\emptyset})=\mathcal E(\mathfrak{E}_{j+r+1,S})$ implies that 
 \[\osum_{\substack{0\leq k\leq j\\ k\text{ odd}}} \begin{bmatrix}r+k\\k\end{bmatrix}_q \mathcal E(\Gamma_{j-k})
 \oplus 
 \osum_{\substack{0\leq k\leq l-j-r-1\\ k\text{ even}}} \begin{bmatrix}r+k\\k\end{bmatrix}_q \mathcal E(\Gamma_{j+k+r+1})
 \] 
 is isomorphic to 
 \[\osum_{\substack{0\leq k\leq j\\ k\text{ even}}} \begin{bmatrix}r+k\\k\end{bmatrix}_q \mathcal E(\Gamma_{j-k})
 \oplus 
 \osum_{\substack{0\leq k\leq l-j-r-1\\ k\text{ odd}}} \begin{bmatrix}r+k\\k\end{bmatrix}_q \mathcal E(\Gamma_{j+k+r+1})
 \] 
 This categorifies the analogues of Kekulé relations which do not use the webs $\Gamma_{j+1},\cdots,\Gamma_{j+r}$. These equations are not included in Equation~\ref{KekuleRelation} for any choice of parameters therein. 
\end{remark}

\subsection{A construction of $\mathfrak{E}_{i,\emptyset}$ in the homotopy category}
In this section we will construct a chain complex whose homology (here by homology we mean chain homotopy equivalent to a complex with zero differential) recovers the results from Proposition \ref{lem:karoubi categorification}. This shows that the morphisms $\mathfrak{E}_{i,\emptyset}$ arise in the homotopy category of $\Foam^\varepsilon$. 

Even though the construction works for any choice of the parameters, we restrict ourselves to following particular case in order to make the construction explicit. 

\begin{example}
    Fix $m=3$ and $a=(0,0,0)$, $b=(1,1,1)$ as in Example~\ref{EgSmoothing} and Section~\ref{sec:munu}. It follows that $r=1$ and $l=4$. The construction starts with the bidirectional long exact sequences, $C_0$ and $C_1$, from Theorem~\ref{thm:hocat categorificaiton} and their shifts. We can arrange those into a two-dimensional anti-commutative diagram as displayed in Figure \ref{fig:2D complex}.
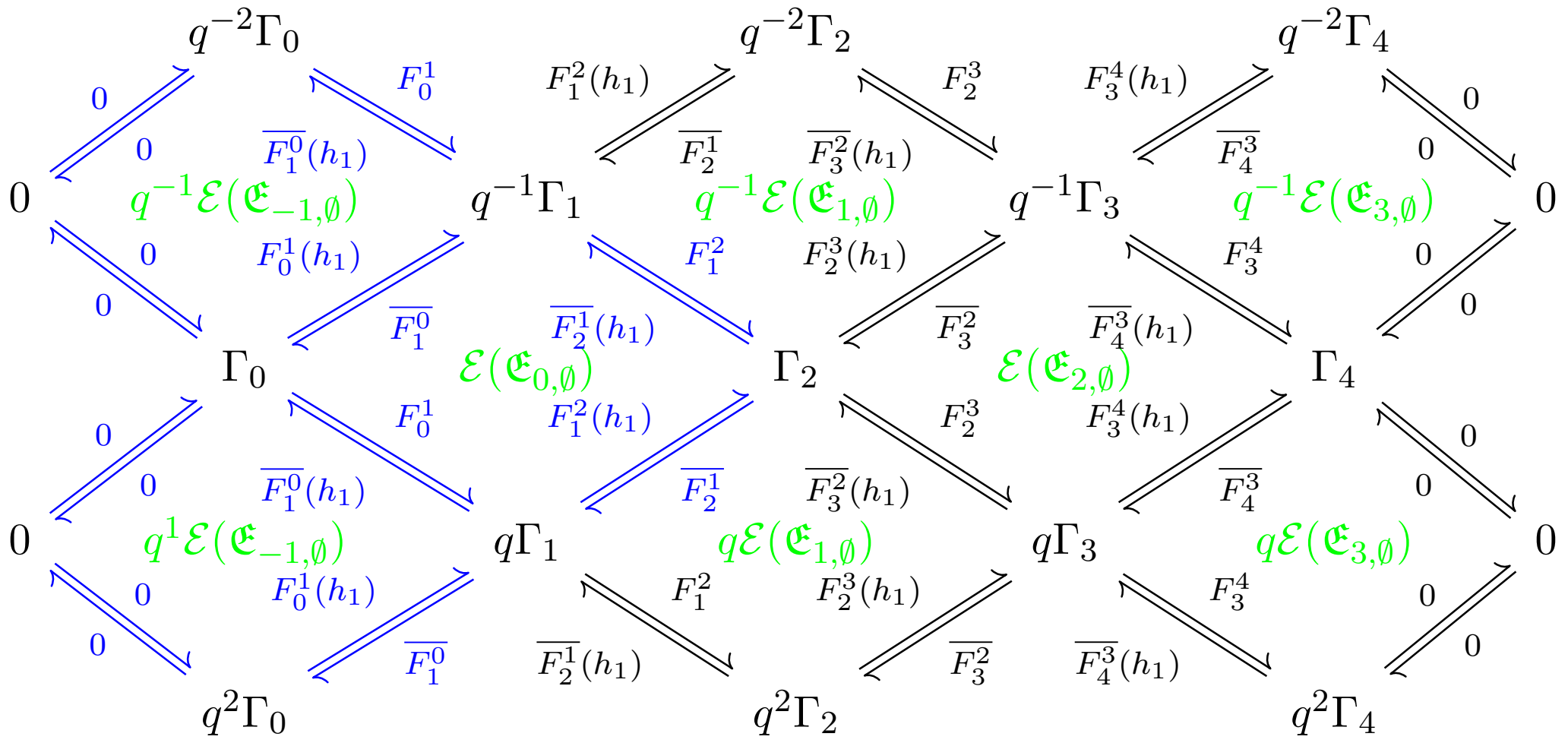
\begin{figure}[htbp]
  \centering
  \begin{tikzcd}[column sep=1.2em, row sep=1.8em, cells={nodes={font=\small}}]
    & q^{-2}\Gamma_0 \arrow[rd, "F_0^1", harpoon, shift left=0.2ex, blue] \arrow[ld, "0", harpoon, shift left=0.2ex, blue]
    &
    & q^{-2}\Gamma_2 \arrow[rd, "F_2^3", harpoon, shift left=0.2ex] \arrow[ld, "\overline{F_2^1}", harpoon, shift left=0.2ex]
    &
    & q^{-2}\Gamma_4 \arrow[rd, "0", harpoon, shift left=0.2ex] \arrow[ld, "\overline{F_4^3}", harpoon, shift left=0.2ex]
    & \\
  0 \arrow[ru, "0", harpoon, shift left=0.2ex, blue] \arrow[rd, "0", harpoon, shift left=0.2ex, blue]
    & {\color{green}q^{-1}\mathcal{E}(\mathfrak{E}_{-1,\emptyset})}
    & q^{-1}\Gamma_1 \arrow[ru, "F_1^2(h_1)", harpoon, shift left=0.2ex] \arrow[rd, "F_1^2", harpoon, shift left=0.2ex, blue] \arrow[lu, "\overline{F_1^0}(h_1)", harpoon, shift left=0.2ex, blue] \arrow[ld, "\overline{F_1^0}", harpoon, shift left=0.2ex, blue]
    & {\color{green}q^{-1}\mathcal{E}(\mathfrak{E}_{1,\emptyset})}
    & q^{-1}\Gamma_3 \arrow[ru, "F_3^4(h_1)", harpoon, shift left=0.2ex] \arrow[rd, "F_3^4", harpoon, shift left=0.2ex] \arrow[lu, "\overline{F_3^2}(h_1)", harpoon, shift left=0.2ex] \arrow[ld, "\overline{F_3^2}", harpoon, shift left=0.2ex]
    & {\color{green}q^{-1}\mathcal{E}(\mathfrak{E}_{3,\emptyset})}
    & 0 \arrow[lu, "0", harpoon, shift left=0.2ex] \arrow[ld, "0", harpoon, shift left=0.2ex] \\
    & \Gamma_0 \arrow[ru, "F_0^1(h_1)", harpoon, shift left=0.2ex, blue] \arrow[rd, "F_0^1", harpoon, shift left=0.2ex, blue] \arrow[lu, "0", harpoon, shift left=0.2ex, blue] \arrow[ld, "0", harpoon, shift left=0.2ex, blue]
    & {\color{green}\mathcal{E}(\mathfrak{E}_{0,\emptyset})}
    & \Gamma_2 \arrow[ru, "F_2^3(h_1)", harpoon, shift left=0.2ex] \arrow[rd, "F_2^3", harpoon, shift left=0.2ex] \arrow[lu, "\overline{F_2^1}(h_1)", harpoon, shift left=0.2ex, blue] \arrow[ld, "\overline{F_2^1}", harpoon, shift left=0.2ex, blue]
    & {\color{green}\mathcal{E}(\mathfrak{E}_{2,\emptyset})}
    & \Gamma_4 \arrow[ru, "0", harpoon, shift left=0.2ex] \arrow[rd, "0", harpoon, shift left=0.2ex] \arrow[lu, "\overline{F_4^3}(h_1)", harpoon, shift left=0.2ex] \arrow[ld, "\overline{F_4^3}", harpoon, shift left=0.2ex]
    & \\
  0 \arrow[ru, "0", harpoon, shift left=0.2ex, blue] \arrow[rd, "0", harpoon, shift left=0.2ex, blue]
    & {\color{green}q^{1}\mathcal{E}(\mathfrak{E}_{-1,\emptyset})}
    & q\Gamma_1 \arrow[ru, "F_1^2(h_1)", harpoon, shift left=0.2ex, blue] \arrow[rd, "F_1^2", harpoon, shift left=0.2ex] \arrow[lu, "\overline{F_1^0}(h_1)", harpoon, shift left=0.2ex, blue] \arrow[ld, "\overline{F_1^0}", harpoon, shift left=0.2ex, blue]
    & {\color{green}q\mathcal{E}(\mathfrak{E}_{1,\emptyset})}
    & q\Gamma_3 \arrow[ru, "F_3^4(h_1)", harpoon, shift left=0.2ex] \arrow[rd, "F_3^4", harpoon, shift left=0.2ex] \arrow[lu, "\overline{F_3^2}(h_1)", harpoon, shift left=0.2ex] \arrow[ld, "\overline{F_3^2}", harpoon, shift left=0.2ex]
    & {\color{green}q\mathcal{E}(\mathfrak{E}_{3,\emptyset})}
    & 0 \arrow[lu, "0", harpoon, shift left=0.2ex] \arrow[ld, "0", harpoon, shift left=0.2ex] \\
    & q^2\Gamma_0 \arrow[ru, "F_0^1(h_1)", harpoon, shift left=0.2ex, blue] \arrow[lu, "0", harpoon, shift left=0.2ex, blue]
    &
    & q^2\Gamma_2 \arrow[ru, "F_2^3(h_1)", harpoon, shift left=0.2ex] \arrow[lu, "\overline{F_2^1}(h_1)", harpoon, shift left=0.2ex]
    &
    & q^2\Gamma_4 \arrow[ru, "0", harpoon, shift left=0.2ex] \arrow[lu, "\overline{F_4^3}(h_1)", harpoon, shift left=0.2ex]
    &
  \end{tikzcd}
  \caption{Two-dimensional anti-commutative complex recovering $\mathfrak{E}_{i,\emptyset}$.}
  \label{fig:2D complex}
\end{figure}

The displayed morphisms $\mathfrak{E}_{j,\emptyset}$ drawn in green are positioned so that each (shifted) web $\Gamma_i$ decomposes in the Karoubi envelope as a sum of the surrounding $\mathfrak{E}_{j,\emptyset}$ according to Theorem \ref{thm:incidence-convolution}.

Alternatively, the positioning of the green $\mathfrak{E}_{j,\emptyset}$ can be justified as follows. Pick any of the (shifted) webs $\Gamma_i$ in the complex from Figure \ref{fig:2D complex}. Starting at $\Gamma_i$ traverse a square (along precisely $4 = 2(r+1)$ differentials) and multiply the maps along the edges. Proposition~\ref{prop:iterated-transfer} tells us which of the $\mathfrak{E}_{j,\emptyset}$ this product recovers, and that is the morphism drawn in green in Figure \ref{fig:2D complex}.

We are now ready to explain how $\mathcal E(\mathfrak{E}_{j,\emptyset})$ can be recovered in a homotopy category. Let $0\leq j\leq l = 4$, construct a chain complex $C^j$ by truncating the two-dimensional complex from Figure~\ref{fig:2D complex} to the chain groups $\Gamma_{j-k}$, $k\geq 0$, with appropriate shifts. It can be shown that the homology of the total chain complex of $C^j$ recovers Equation \eqref{eq:inverse-incidence}. That is, every $\mathcal E(\mathfrak{E}_{j,\emptyset})$ can be realized in the homotopy category of $\Foam^\varepsilon$. Moreover, for $l - r  = 3 \leq j \leq l = 4$,  the total complex of $C^j$ is null-homotopic, which leads to another categorification of the Kekulé relations as in Proposition \ref{lem:karoubi categorification}.

For example, consider the truncation, corresponding to $j = 2$, with chain groups $\Gamma_2$, $q^{-1}\Gamma_1$, $q\Gamma_1$, $q^{-2}\Gamma_0$, $\Gamma_0$, $q^2\Gamma_0$ and four $0$'s (two of which are not displayed in Figure \ref{fig:2D complex}), whose differentials are colored blue in Figure \ref{fig:2D complex}. By the above discussion, it is straightforward to check that the homology of the total complex of this truncation recovers $\mathcal{E}(\mathfrak{E}_{2,\emptyset})$ as a part of Equation \eqref{eq:inverse-incidence}
\begin{align*}
\mathcal{E}(\mathfrak{E}_{2,\emptyset})\oplus[2]_q\mathcal{E}(\Gamma_1)=\mathcal{E}(\Gamma_2)\oplus[3]_q\mathcal{E}(\Gamma_0).
\end{align*}
\end{example}
This construction can be generalized to any choice of $m,a$ and $b$. The resulting analogue of the two-dimensional complex in Figure \ref{fig:2D complex} would be $(r+1)$-dimensional. More specifically, for every $0\leq j\leq l$, we can construct a truncated $(r+1)$-complex from Equation~\eqref{LESsF}, whose chain groups are $\Gamma_{j-k}$ with $k\geq 0$ and proper shifts that finally goes to $\Gamma_j$ and differentials same with Equation~\eqref{LESsF}. We can prove that this $(r+1)$-complex contracts to $\mathcal E(\mathfrak{E}_{j,\emptyset})$ in the homotopy category of $\Kar(\Foam^\varepsilon).$ And this is to say, each $\mathcal E(\mathfrak{E}_{j,\emptyset})$ can be realize as a object in homotopy category of $\Foam^\varepsilon$.

\section{Categorification of the $6$-valent vertex}\label{sec:munu}
In this section we resolve \cite[Conjecture 4.29]{grlj20266valentvertexmathfrakglnweb} posed by M. Khovanov, H. Wu, M. Zhang, and the first author. The relation categorified in that paper is not quite the same as the Kekulé relations that we have already categorified, see the next definition or \cite[Equation 2]{grlj20266valentvertexmathfrakglnweb}.

\def\SixVertex
{\begin{tric}
\draw [decoration={markings,mark=at position 0.6 with {\arrow{>}}},postaction={decorate}] (0,0)--(2,0);
\draw [decoration={markings,mark=at position 0.6 with {\arrow{>}}},postaction={decorate}](-2,0)--(0,0);
\draw [decoration={markings,mark=at position 0.6 with {\arrow{>}}},postaction={decorate}] (1,1.7)--(0,0);
\draw [decoration={markings,mark=at position 0.6 with {\arrow{>}}},postaction={decorate}] (0,0)--(-1,-1.7);
\draw [decoration={markings,mark=at position 0.6 with {\arrow{>}}},postaction={decorate}] (1,-1.7)--(0,0);
\draw [decoration={markings,mark=at position 0.6 with {\arrow{>}}},postaction={decorate}] (0,0)--(-1,1.7);
\end{tric}
}

\begin{definition}\label{def:six valent}
Fix the parameters from Definition \ref{def:kekule} as $m = 3, a = (0,0,0),$ $b = (1,1,1),$ and $N \geq 4$. Therefore $l = N -1,$ and $r = N-4$. \emph{The $6$-valent vertex $\Wsix$} is a morphism in $\Web^\varepsilon$, given by:
    \begin{align} \label{GLN6Vertex}
    \Wsix &= \ \SixVertex \ :=  \Gamma_1  - [N-3]_q \Gamma_0 . 
\end{align}
\end{definition}

This relation was categorified in the aforementioned work for $N = 4,5$, where $\Gamma_0, \Gamma_1$ were denoted as $\Wone, \Wzero$ respectively. The following statements resolve the conjecture for all $N \geq 4$, using the results from previous sections.

\begin{lemma}\label{lem:six vertex categorification}
    For $N \ge 4$, the following is a decomposition into mutually orthogonal idempotents:
    \begin{align}\label{eq:six vertex categorification}
        \idfoam_{\Gamma_1} &= \sum_{s  = 0 }^{N - 4}  \Li{s}{1} + \prod_{s = 0}^{N-4}  \Ri{s}{1}.
    \end{align}
\end{lemma}
\begin{proof}
The orthogonality and idempotency follow from Lemma \ref{lem:Li,Ri prop}. Using Lemma \ref{lem:prod of Ri Li vanishing} and Equation \eqref{eq:main L,R equation}, compute:
    \begin{align*}
        \idfoam_{\Gamma_1}  &= \prod_{s = 0}^{N-4} \idfoam_{\Gamma_1} \\
        &= \prod_{s = 0}^{N-4} \left( \Li{s}{1} + \Ri{s}{1}\right) \\
        &= \sum_{s  = 0 }^{N - 4} \Li{s}{1} \prod_{\substack{s' \in \{0,\ldots,N-4\}\\ s' \neq s}} \Ri{s'}{1} + \prod_{s = 0}^{N-4} \Ri{s}{1} \\
        &= \sum_{s  = 0 }^{N - 4} \Li{s}{1} \prod_{\substack{s' \in \{0,\ldots,N-4\}\\ s' \neq s}} \left( \idfoam_{\Gamma_1} - \Li{s'}{1}\right) + \prod_{s = 0}^{N-4} \Ri{s}{1} \\
        &= \sum_{s  = 0 }^{N - 4} \Li{s}{1} + \prod_{s = 0}^{N-4} \Ri{s}{1}
    \end{align*}
\end{proof}

\begin{proposition}\label{lem:munu as Ri}
The following foam equality holds:
\begin{align*}
    \prod_{s = 0}^{N-4} \Ri{s}{1} = \mu \cdot \nu. 
\end{align*}
See Figure \ref{fig:munu} for a movie diagram of the foam $\mu \cdot \nu$. 
\begin{figure}[htbp]
  \centering
  \newcommand{\ExMovieStyleScale}{0.85}
  \scalebox{\ExMovieStyleScale}{%
  \begin{tikzpicture}[
    scale=0.5,
    line/.style={
      line width=0.9pt,
      line cap=round,
      line join=round
    },
    flow/.style={
      line,
      postaction={decorate},
      decoration={
        markings,
        mark=at position 0.53 with
        {\arrow{Stealth[length=2.2mm,width=1.7mm]}}
      }
    },
    every node/.style={
      font=\normalsize,
      inner sep=1.5pt
    }
  ]

%==================================================
% Left diagram
%==================================================
\begin{scope}

% Circumradius of the regular hexagon
\def\R{1.65}
\pgfmathsetmacro{\H}{sqrt(3)*\R/2}

% Vertices of the regular hexagon
\coordinate (LT) at (-0.5*\R, \H);
\coordinate (RT) at ( 0.5*\R, \H);
\coordinate (RM) at ( \R,0);
\coordinate (RB) at ( 0.5*\R,-\H);
\coordinate (LB) at (-0.5*\R,-\H);
\coordinate (LM) at (-\R,0);

% Outer endpoints
\coordinate (LTout) at (-1*\R, 2*\H);
\coordinate (RTout) at ( 1*\R, 2*\H);
\coordinate (RMout) at ( 2*\R,0);
\coordinate (RBout) at ( 1*\R,-2*\H);
\coordinate (LBout) at (-1*\R,-2*\H);
\coordinate (LMout) at (-2*\R,0);

% Directed edges of the hexagon
\draw[flow] (RT) -- (LT)
  node[pos=.52,above=5pt] {$2$};

\draw[flow] (LT) -- (LM)
  node[pos=.50,left=4pt] {$1$};

\draw[flow] (LM) -- (LB)
  node[pos=.54,left=4pt] {$2$};

\draw[flow] (LB) -- (RB)
  node[pos=.50,below=5pt] {$1$};

\draw[flow] (RB) -- (RM)
  node[pos=.50,right=4pt] {$2$};

\draw[flow] (RM) -- (RT)
  node[pos=.50,right=4pt] {$1$};

% Six outward-extending branches
\draw[flow] (LT) -- (LTout)
  node[pos=.93,left=3pt] {$1$};

\draw[flow] (RTout) -- (RT)
  node[pos=.07,right=3pt] {$1$};

\draw[flow] (LMout) -- (LM)
  node[pos=.08,below=4pt] {$1$};

\draw[flow] (RM) -- (RMout)
  node[pos=.92,below=4pt] {$1$};

\draw[flow] (LB) -- (LBout)
  node[pos=.93,left=3pt] {$1$};

\draw[flow] (RBout) -- (RB)
  node[pos=.07,right=3pt] {$1$};

\end{scope}

\begin{scope}[xshift=7cm]

% Circumradius of the regular hexagon
\def\R{1.65}
\pgfmathsetmacro{\H}{sqrt(3)*\R/2}

\node at (0,0) {$3$};

% Vertices of the regular hexagon
\coordinate (LT) at (-0.5*\R, \H);
\coordinate (RT) at ( 0.5*\R, \H);
\coordinate (RM) at ( \R,0);
\coordinate (RB) at ( 0.5*\R,-\H);
\coordinate (LB) at (-0.5*\R,-\H);
\coordinate (LM) at (-\R,0);

% Outer endpoints
\coordinate (LTout) at (-1*\R, 2*\H);
\coordinate (RTout) at ( 1*\R, 2*\H);
\coordinate (RMout) at ( 2*\R,0);
\coordinate (RBout) at ( 1*\R,-2*\H);
\coordinate (LBout) at (-1*\R,-2*\H);
\coordinate (LMout) at (-2*\R,0);

% Directed edges of the hexagon
% \draw[flow] (RT) -- (LT);
% \draw[flow] (LT) -- (LM);
% \draw[flow] (LM) -- (LB);
% \draw[flow] (LB) -- (RB);
% \draw[flow] (RB) -- (RM);
% \draw[flow] (RM) -- (RT);
\draw[flow] (LT) -- (RT);
\draw[flow] (LM) -- (LT);
\draw[flow] (LB) -- (LM);
\draw[flow] (RB) -- (LB);
\draw[flow] (RM) -- (RB);
\draw[flow] (RT) -- (RM);

% Six outward-extending branches
\draw[line] (LT) -- (LTout);
\draw[line] (RTout) -- (RT);
\draw[line] (LMout) -- (LM);
\draw[line] (RM) -- (RMout);
\draw[line] (LB) -- (LBout);
\draw[line] (RBout) -- (RB);
\end{scope}

\begin{scope}[xshift=14cm]

% Circumradius of the regular hexagon
\def\R{1.65}
\pgfmathsetmacro{\H}{sqrt(3)*\R/2}

% Vertices of the regular hexagon
\coordinate (LT) at (-0.5*\R, \H);
\coordinate (RT) at ( 0.5*\R, \H);
\coordinate (RM) at ( \R,0);
\coordinate (RB) at ( 0.5*\R,-\H);
\coordinate (LB) at (-0.5*\R,-\H);
\coordinate (LM) at (-\R,0);

% Outer endpoints
\coordinate (LTout) at (-1*\R, 2*\H);
\coordinate (RTout) at ( 1*\R, 2*\H);
\coordinate (RMout) at ( 2*\R,0);
\coordinate (RBout) at ( 1*\R,-2*\H);
\coordinate (LBout) at (-1*\R,-2*\H);
\coordinate (LMout) at (-2*\R,0);

% Directed edges of the hexagon
\draw[flow] (LT) -- (RT)
  node[pos=.52,above=5pt] {$1$};

\draw[flow] (LM) -- (LT)
  node[pos=.50,left=4pt] {$2$};

\draw[flow] (LB) -- (LM)
  node[pos=.54,left=4pt] {$1$};

\draw[flow] (RB) -- (LB)
  node[pos=.50,below=5pt] {$2$};

\draw[flow] (RM) -- (RB)
  node[pos=.50,right=4pt] {$1$};

\draw[flow] (RT) -- (RM)
  node[pos=.50,right=4pt] {$2$};

% Six outward-extending branches
\draw[flow] (LT) -- (LTout)
  node[pos=.93,left=3pt] {$1$};

\draw[flow] (RTout) -- (RT)
  node[pos=.07,right=3pt] {$1$};

\draw[flow] (LMout) -- (LM)
  node[pos=.08,below=4pt] {$1$};

\draw[flow] (RM) -- (RMout)
  node[pos=.92,below=4pt] {$1$};

\draw[flow] (LB) -- (LBout)
  node[pos=.93,left=3pt] {$1$};

\draw[flow] (RBout) -- (RB)
  node[pos=.07,right=3pt] {$1$};

\end{scope}

\begin{scope}[xshift=21cm]

% Circumradius of the regular hexagon
\def\R{1.65}
\pgfmathsetmacro{\H}{sqrt(3)*\R/2}

\node at (0,0) {$3$};

% Vertices of the regular hexagon
\coordinate (LT) at (-0.5*\R, \H);
\coordinate (RT) at ( 0.5*\R, \H);
\coordinate (RM) at ( \R,0);
\coordinate (RB) at ( 0.5*\R,-\H);
\coordinate (LB) at (-0.5*\R,-\H);
\coordinate (LM) at (-\R,0);

% Outer endpoints
\coordinate (LTout) at (-1*\R, 2*\H);
\coordinate (RTout) at ( 1*\R, 2*\H);
\coordinate (RMout) at ( 2*\R,0);
\coordinate (RBout) at ( 1*\R,-2*\H);
\coordinate (LBout) at (-1*\R,-2*\H);
\coordinate (LMout) at (-2*\R,0);

% Directed edges of the hexagon
\draw[flow] (RT) -- (LT);
\draw[flow] (LT) -- (LM);
\draw[flow] (LM) -- (LB);
\draw[flow] (LB) -- (RB);
\draw[flow] (RB) -- (RM);
\draw[flow] (RM) -- (RT);
% \draw[flow] (LT) -- (RT);
% \draw[flow] (LM) -- (LT);
% \draw[flow] (LB) -- (LM);
% \draw[flow] (RB) -- (LB);
% \draw[flow] (RM) -- (RB);
% \draw[flow] (RT) -- (RM);

% Six outward-extending branches
\draw[line] (LT) -- (LTout);
\draw[line] (RTout) -- (RT);
\draw[line] (LMout) -- (LM);
\draw[line] (RM) -- (RMout);
\draw[line] (LB) -- (LBout);
\draw[line] (RBout) -- (RB);
\end{scope}

%==================================================
% Right diagram
%==================================================
\begin{scope}[xshift=28cm]

% Circumradius of the regular hexagon
\def\R{1.65}
\pgfmathsetmacro{\H}{sqrt(3)*\R/2}

% Vertices of the regular hexagon
\coordinate (LT) at (-0.5*\R, \H);
\coordinate (RT) at ( 0.5*\R, \H);
\coordinate (RM) at ( \R,0);
\coordinate (RB) at ( 0.5*\R,-\H);
\coordinate (LB) at (-0.5*\R,-\H);
\coordinate (LM) at (-\R,0);

% Outer endpoints
\coordinate (LTout) at (-1*\R, 2*\H);
\coordinate (RTout) at ( 1*\R, 2*\H);
\coordinate (RMout) at ( 2*\R,0);
\coordinate (RBout) at ( 1*\R,-2*\H);
\coordinate (LBout) at (-1*\R,-2*\H);
\coordinate (LMout) at (-2*\R,0);

% Directed edges of the hexagon
\draw[flow] (RT) -- (LT)
  node[pos=.52,above=5pt] {$2$};

\draw[flow] (LT) -- (LM)
  node[pos=.50,left=4pt] {$1$};

\draw[flow] (LM) -- (LB)
  node[pos=.54,left=4pt] {$2$};

\draw[flow] (LB) -- (RB)
  node[pos=.50,below=5pt] {$1$};

\draw[flow] (RB) -- (RM)
  node[pos=.50,right=4pt] {$2$};

\draw[flow] (RM) -- (RT)
  node[pos=.50,right=4pt] {$1$};

% Six outward-extending branches
\draw[flow] (LT) -- (LTout)
  node[pos=.93,left=3pt] {$1$};

\draw[flow] (RTout) -- (RT)
  node[pos=.07,right=3pt] {$1$};

\draw[flow] (LMout) -- (LM)
  node[pos=.08,below=4pt] {$1$};

\draw[flow] (RM) -- (RMout)
  node[pos=.92,below=4pt] {$1$};

\draw[flow] (LB) -- (LBout)
  node[pos=.93,left=3pt] {$1$};

\draw[flow] (RBout) -- (RB)
  node[pos=.07,right=3pt] {$1$};

\end{scope}
\end{tikzpicture}}
\caption{A movie depiction of the foam $\mu \cdot \nu$.}
\label{fig:munu}
\end{figure}
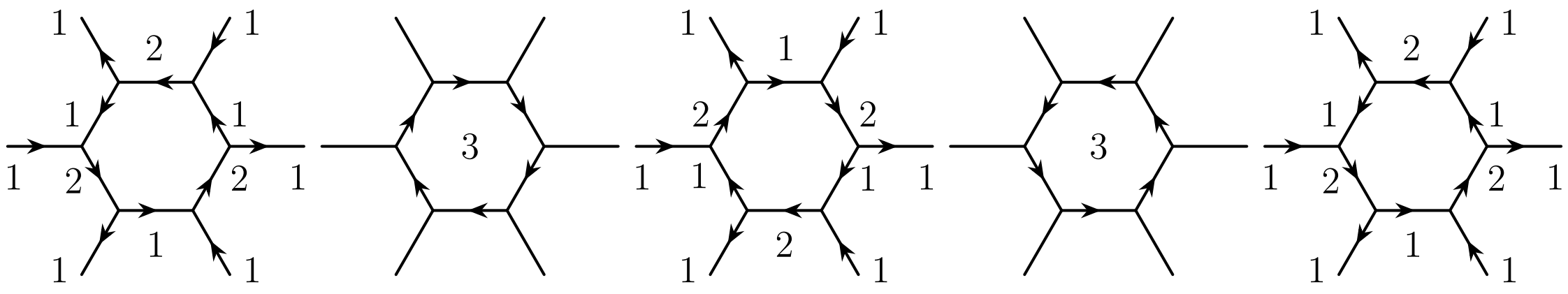
\end{proposition}
\begin{proof}
We split the proof of this proposition into three lemmas. Lemma \ref{lem:prod Ri N-3 N-3} will establish that $\prod_{s = 0}^{N-4} \Ri{s}{1} = (-1)^{\frac{1}{2} (N +1)N} F_{N-2}^1 F_{1}^{N-2}$. Lemma \ref{lem:D2 = id on gammaN-2}, will show that $F_{N-2}^1 F_{1}^{N-2} = (-1)^{\frac{1}{2} (N +1)N} F_{N-2}^1 D^2 F_{1}^{N-2}$ and lastly, Lemma \ref{lem:munu = n-3 D D n-3} demonstrates that $\mu \cdot \nu = F_{N-2}^1 D^2 F_{1}^{N-2}$, completing the proof. The latter two results follow from explicit foam evaluation.
\end{proof}

\begin{remark}
    It was shown in \cite[Theorem 4.2]{grlj20266valentvertexmathfrakglnweb} that $\mathcal E(\Wsix)$ is indecomposable when $\mu \cdot \nu$ is a projection. Hence Proposition \ref{lem:munu as Ri} shows that $\mathcal E(\Wsix)$ is indecomposable in $\Kar(\Foam^\varepsilon)$.
\end{remark}

\begin{lemma}\label{lem:prod Ri N-3 N-3}
    The following foam equality holds:
    \begin{align*}
        \prod_{s = 0}^{N-4} \Ri{s}{1} = (-1)^{\frac{1}{2} (N +1)N} F_{N-2}^1 F_{1}^{N-2}.
    \end{align*}
\end{lemma}
\begin{proof}
    This follows from iteratively applying Corollaries \ref{cor:iterated composition of slice morphisms} and \ref{cor:composition of slice morphisms}. Expanding the expression on the left (we understand $\prod_{i = a}^b c_i = c_b c_{b-1}\cdots c_a$):
    \begin{align*}
        \prod_{s = 0}^{N-4} \Ri{s}{1} &=\prod_{s = 0}^{N-4} \overline{F_{2}^1} (h_s)F_{1}^{2} (h_{N - 4-s}) \\
         &= \prod_{s = 0}^{N-4} \overline{F_{N -2 -s}^{N- 3 -s}} (h_s)  \prod_{s = 0}^{N-4} F_{1 + s}^{2 + s} (h_{N - 4-s}) \\
        % &=  \overline{F_{1}^{0}}(h_{N-4}) \cdot \overline{F_{2}^{1}}(h_{N-3}) \cdots\overline{F_{N -2}^{N- 3}} (h_0)  F_{N-2}^{N-3}(h_0) \cdot F_{N-1}^{N-2}(h_1)  \cdots F_{1 }^{2} (h_{N - 4})\\
        &=    \left(\prod_{s = 0}^{N-4}\sum_{v = 0}^s(-1)^v E_v(x_{\overline{B}}) F_{N - 2 -s}^{N- 3 -s} (h_{s-v}) \right)  F_1^{N-2} \\
        &=    \left(\prod_{s = 0}^{N-4}F_{N - 2 -s}^{N- 3-s} (h_{s}) \right)  F_1^{N-2} \\
        &=  (-1)^{\frac{1}{2} (N-4)(N-3)} F_{N-2}^1 F_{1}^{N-2},
    \end{align*}
    where the rearrangement and cancellation of terms is justified by the aforementioned corollaries and Lemma \ref{lem:Li,Ri prop}. The result follows by noting that $(-1)^{\frac{1}{2} (N-4)(N-3)} = (-1)^{\frac{1}{2} (N +1)N}$.
\end{proof}
\begin{remark}\label{rmk:Ris generalizeation}
    The same proof shows that for $t\neq 0$, $\prod_{s = 0}^{t} \Ri{s}{i} = (-1)^{\frac{1}{2} (t +1)t} F_{i + t + 1}^i F_{i}^{i + t + 1}$.
\end{remark} 

\begin{lemma}\label{lem:D2 = id on gammaN-2}
    The following foam equality holds:
    \begin{align*}
        F_{N-2}^1 F_{1}^{N-2} = (-1)^{\frac{1}{2} (N +1)N} F_{N-2}^1 D^2 F_{1}^{N-2}.
    \end{align*}
\end{lemma}
\begin{proof}
    This follows from foam evaluation. Note that Remark \ref{rmk:D4 = id} says that $F_{N-2}^1 F_{1}^{N-2} = F_{N-2}^1 D^4 F_{1}^{N-2}$. It remains to show is that $D^2$ acts as $ (-1)^{\frac{1}{2} (N +1)N} \idfoam_{\Gamma_1'}$ on $\Gamma_1'$. Here $\Gamma_1'$ denotes the middle web in Figure \ref{fig:munu}, that is the codomain of $DF_{1}^{N-2}$.

    We will compare $D^2$ with $\idfoam_{\Gamma_1'}$. There are three distinct families of colorings of $\Gamma_1'$ as displayed in Figure \ref{fig:colorings of Gamma1'}. 
\begin{figure}[htbp]
\centering
\begin{tikzpicture}[
  scale=0.5,
  web edge/.style={
    line width=0.9pt,
    line cap=round,
    line join=round
  },
  web arrow/.style args={#1}{
    postaction={decorate},
    decoration={
      markings,
      mark=at position 0.53 with {
        \arrow[draw=#1,fill=#1]{
          Stealth[length=2.2mm,width=1.7mm]
        }
      }
    }
  }
]

\foreach \panel in {0,1,2} {
  \begin{scope}[shift={({7.5*\panel},0)}]

  % Hexagon vertices
  \def\R{1.65}
  \pgfmathsetmacro{\H}{sqrt(3)*\R/2}

  \coordinate (LT) at (-0.5*\R, \H);
  \coordinate (RT) at ( 0.5*\R, \H);
  \coordinate (RM) at ( \R,0);
  \coordinate (RB) at ( 0.5*\R,-\H);
  \coordinate (LB) at (-0.5*\R,-\H);
  \coordinate (LM) at (-\R,0);

  % Outer endpoints
  \coordinate (LTout) at (-\R, 2*\H);
  \coordinate (RTout) at ( \R, 2*\H);
  \coordinate (RMout) at ( 2*\R,0);
  \coordinate (RBout) at ( \R,-2*\H);
  \coordinate (LBout) at (-\R,-2*\H);
  \coordinate (LMout) at (-2*\R,0);

\ifcase\panel
  % First coloring
  \def\singleedges{%
    LT/RT/green,%
    LB/LM/green,%
    RM/RB/green,%
    LT/LTout/blue,%
    RTout/RT/blue,%
    LMout/LM/blue,%
    RM/RMout/blue,%
    LB/LBout/blue,%
    RBout/RB/blue%
  }
  \def\doubleedges{%
    LM/LT/blue/green,%
    RB/LB/blue/green,%
    RT/RM/blue/green%
  }
\or
  % Second coloring
  \def\singleedges{%
    LT/RT/red,%
    LB/LM/blue,%
    RM/RB/green,%
    LT/LTout/blue,%
    RTout/RT/green,%
    LMout/LM/red,%
    RM/RMout/red,%
    LB/LBout/green,%
    RBout/RB/blue%
  }
  \def\doubleedges{%
    LM/LT/blue/red,%
    RB/LB/blue/green,%
    RT/RM/green/red%
  }
\or
  % Third coloring
  \def\singleedges{%
    LT/RT/green,%
    LB/LM/green,%
    RM/RB/red,%
    LT/LTout/blue,%
    RTout/RT/red,%
    LMout/LM/blue,%
    RM/RMout/green,%
    LB/LBout/red,%
    RBout/RB/green%
  }
  \def\doubleedges{%
    LM/LT/blue/green,%
    RB/LB/red/green,%
    RT/RM/green/red%
  }
\fi

  % Single edges with matching colored arrowheads
  \foreach \from/\to/\edgecolor in \singleedges {
    \draw[
      web edge,
      color=\edgecolor,
      web arrow={\edgecolor}
    ] (\from) -- (\to);
  }

  % Double edges with one black arrowhead
  \foreach \from/\to/\firstcolor/\secondcolor in \doubleedges {
    \draw[web edge,color=\firstcolor]
      ($(\from)!0.06cm!90:(\to)$) --
      ($(\to)!0.06cm!-90:(\from)$);

    \draw[web edge,color=\secondcolor]
      ($(\from)!0.06cm!-90:(\to)$) --
      ($(\to)!0.06cm!90:(\from)$);

    \path[web arrow={black}] (\from) -- (\to);
  }

  \end{scope}
}

\end{tikzpicture}
\caption{The three different families of colorings of $\Gamma_1'$.}
\label{fig:colorings of Gamma1'}
\end{figure}
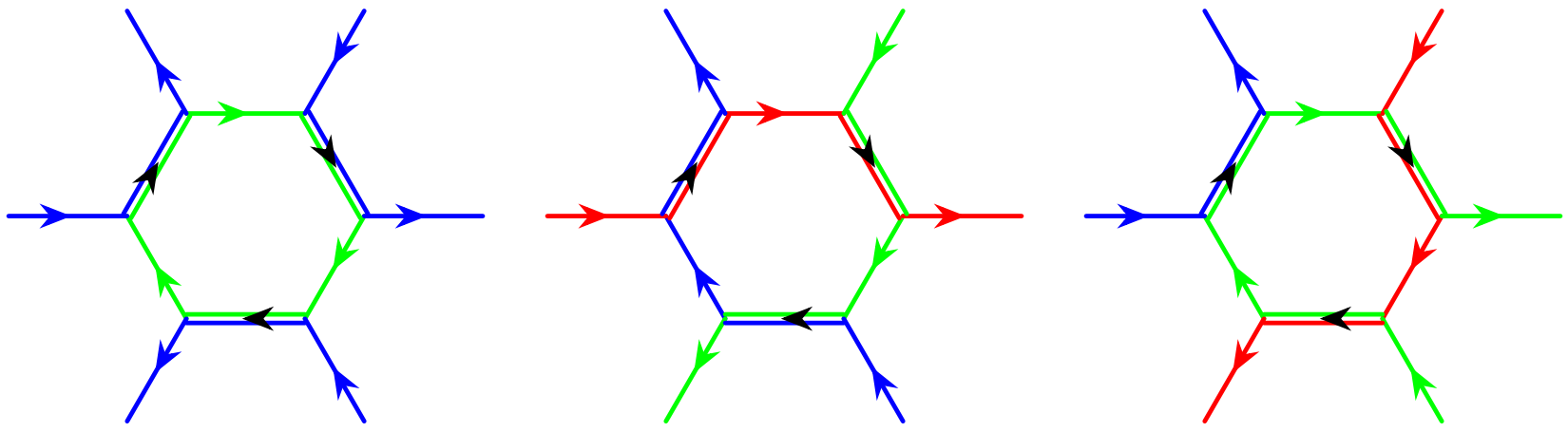

For readability, write $G = D^2 -  \idfoam_{\Gamma_1'}$. We will first analyze the leftmost coloring of Figure \ref{fig:colorings of Gamma1'}, call it $c_1$ and extend it to a coloring of foam $G$. All the bichrome surfaces of $(G, c_1)$ have trivial Euler characteristic. Denote by $o\in O$, $a\in A$, $b_k\in B_k$, and $b_l\in B_l$ arbitrary elements of these sets (see Definition \ref{def:smoothing}). Then
  \begin{itemize}
    \item The monochrome surfaces have Euler characteristics $\chi(G_{o}, c_1)=2, \chi(G_a, c_1) = 2, \chi(G_{b_i}, c_1) = -2(i-1)$, where $b_i$ denotes either $b_k$ or $b_l$.
    \item The positive circles can be counted as $\theta^+_{ao}(G, c_1) = 2\mathbbm{1}\{o > a\}$, $\theta^+_{ab_i}(G, c_1) = \pm k (1 - \mathbbm{1}\{b_i > a\})$, $\theta^+_{b_i o}(G, c_1) = \pm i  \mathbbm{1}\{b_i > o\}$, where $b_i$ denotes either $b_k$ or $b_l$. Lastly, $\theta^+_{b_k b_l}(G, c_1) = k \mathbbm{1}\{b_k > b_l\} +  l \mathbbm{1}\{b_l > b_k\}$.
  \end{itemize}
  Sum these contributions:
  \begin{align*}
    \sum_{i=1}^N {i \chi(F_i(G,c_1))/2}  &\equiv a + \sum_{o \in O} o + \sum_{k = 1}^N\sum_{b \in B_k} (k+1) b \\
    \sum_{1\leq i < j \leq N} \theta^+_{ij}(G,c_1) & \equiv \sum_{k = 1}^N\sum_{b \in B_k} k b \\
    s(G,c_1) &\equiv  a + \sum_{o \in O} o  + \sum_{k = 1}^N\sum_{b \in B_k}  b  = \sum_{k = 1}^N k = \frac{N(N+1)}{2}
  \end{align*}

Next, we analyze the middle coloring of Figure \ref{fig:colorings of Gamma1'}, call it $c_2$. Denote the green pigment on the middle web of Figure \ref{fig:colorings of Gamma1'} by $f$, the red pigment by $g$ and the blue pigment by $e$, denote by $o$ an arbitrary element of $O$. Note that $e,f,g$ must all be distinct colors, without loss of generality, assume $e<f<g$. All the bichrome surfaces of $(G, c_2)$ have trivial Euler characteristic. Then
  \begin{itemize}
    \item The monochrome surfaces have Euler characteristics $0$, except for $\chi(G_{o}, c_2)=2$.
    \item The positive circles can be counted as $\theta^+_{ef}(G, c_2) = \theta^+_{gf}(G, c_2) = \theta^+_{eg}(G, c_2) = 0$, $\theta^+_{bo}(G, c_2) = \mathbbm{1}\{b > o\}$, where $b$ denotes an element of $\{e,f,g\}$.
  \end{itemize}
  Add the contributions:
  \begin{align*}
    \sum_{i=1}^N {i \chi(F_i(G,c_2))/2}  &\equiv \sum_{o \in O} o  = \frac{N(N+1)}{2} - e - f - g\\
    \sum_{1\leq i < j \leq N} \theta^+_{ij}(G,c_2) & \equiv \sum_{k \in \{e,f,g\}}^N\sum_{o \in O} \mathbbm{1}\{b > o\} \equiv e -1 + f -2 + g -3 \\
    s(G,c_2) &\equiv \frac{N(N+1)}{2} - e - f - g +  e -1 + f -2 + g -3 \equiv\frac{N(N+1)}{2}.
  \end{align*}

We omit the evaluation for the rightmost coloring of Figure \ref{fig:colorings of Gamma1'} as it is analogous to the previous one.
\end{proof}

\begin{remark}
    As in Remark \ref{rmk:Ris generalizeation}, Lemma~\ref{lem:D2 = id on gammaN-2} can be generalized to show 
    \begin{align*}
        F_{N-1-i}^i F_{i}^{N-1-i} = (-1)^{\frac{1}{2} (N +1)N+i-1} F_{N-1-i}^i D^2 F_{i}^{N-1-i},\quad i=1,2,\cdots,\left\lfloor \frac{N-1}{2}\right\rfloor.
    \end{align*}
    The argument from the proof of Lemma \ref{lem:munu = n-3 D D n-3} bellow then generalizes to show that 
    \begin{align*}
        \xi_i^1 \circ \xi_i^2=\prod_{s = 0}^{N-2-2i} \Ri{s}{i} = (-1)^{\frac{1}{2} (t +1)t} F_{i + t + 1}^i F_{i}^{i + t + 1},
    \end{align*}
    where $\xi_i^1 \circ \xi_i^2$ is depicted in Figure \ref{fig:xi} and should be understood as the generalization of $\mu \cdot \nu = \xi_1^1 \circ \xi_1^2$.
    \begin{figure}[htbp]
  \centering
  \newcommand{\ExMovieStyleScale}{0.85}
  \scalebox{\ExMovieStyleScale}{%
  \begin{tikzpicture}[
    scale=0.5,
    line/.style={
      line width=0.9pt,
      line cap=round,
      line join=round
    },
    flow/.style={
      line,
      postaction={decorate},
      decoration={
        markings,
        mark=at position 0.53 with
        {\arrow{Stealth[length=2.2mm,width=1.7mm]}}
      }
    },
    every node/.style={
      font=\normalsize,
      inner sep=1.5pt
    }
  ]

%==================================================
% Left diagram
%==================================================
\begin{scope}

% Circumradius of the regular hexagon
\def\R{1.65}
\pgfmathsetmacro{\H}{sqrt(3)*\R/2}

% Vertices of the regular hexagon
\coordinate (LT) at (-0.5*\R, \H);
\coordinate (RT) at ( 0.5*\R, \H);
\coordinate (RM) at ( \R,0);
\coordinate (RB) at ( 0.5*\R,-\H);
\coordinate (LB) at (-0.5*\R,-\H);
\coordinate (LM) at (-\R,0);

% Outer endpoints
\coordinate (LTout) at (-1*\R, 2*\H);
\coordinate (RTout) at ( 1*\R, 2*\H);
\coordinate (RMout) at ( 2*\R,0);
\coordinate (RBout) at ( 1*\R,-2*\H);
\coordinate (LBout) at (-1*\R,-2*\H);
\coordinate (LMout) at (-2*\R,0);

% Directed edges of the hexagon
\draw[flow] (RT) -- (LT)
  node[pos=.52,above=5pt] {$i+1$};

\draw[flow] (LT) -- (LM)
  node[pos=.50,left=4pt] {$i$};

\draw[flow] (LM) -- (LB)
  node[pos=.54,left=4pt] {$i+1$};

\draw[flow] (LB) -- (RB)
  node[pos=.50,below=5pt] {$i$};

\draw[flow] (RB) -- (RM)
  node[pos=.50,right=4pt] {$i+1$};

\draw[flow] (RM) -- (RT)
  node[pos=.50,right=4pt] {$i$};

% Six outward-extending branches
\draw[flow] (LT) -- (LTout)
  node[pos=.93,left=3pt] {$1$};

\draw[flow] (RTout) -- (RT)
  node[pos=.07,right=3pt] {$1$};

\draw[flow] (LMout) -- (LM)
  node[pos=.08,above=4pt] {$1$};

\draw[flow] (RM) -- (RMout)
  node[pos=.92,above=4pt] {$1$};

\draw[flow] (LB) -- (LBout)
  node[pos=.93,left=3pt] {$1$};

\draw[flow] (RBout) -- (RB)
  node[pos=.07,right=3pt] {$1$};

\end{scope}

\begin{scope}[xshift=7cm]

% Circumradius of the regular hexagon
\def\R{1.65}
\pgfmathsetmacro{\H}{sqrt(3)*\R/2}

\node at (0,0) {$2i+1$};

% Vertices of the regular hexagon
\coordinate (LT) at (-0.5*\R, \H);
\coordinate (RT) at ( 0.5*\R, \H);
\coordinate (RM) at ( \R,0);
\coordinate (RB) at ( 0.5*\R,-\H);
\coordinate (LB) at (-0.5*\R,-\H);
\coordinate (LM) at (-\R,0);

% Outer endpoints
\coordinate (LTout) at (-1*\R, 2*\H);
\coordinate (RTout) at ( 1*\R, 2*\H);
\coordinate (RMout) at ( 2*\R,0);
\coordinate (RBout) at ( 1*\R,-2*\H);
\coordinate (LBout) at (-1*\R,-2*\H);
\coordinate (LMout) at (-2*\R,0);

% Directed edges of the hexagon
% \draw[flow] (RT) -- (LT);
% \draw[flow] (LT) -- (LM);
% \draw[flow] (LM) -- (LB);
% \draw[flow] (LB) -- (RB);
% \draw[flow] (RB) -- (RM);
% \draw[flow] (RM) -- (RT);
\draw[flow] (LT) -- (RT);
\draw[flow] (LM) -- (LT);
\draw[flow] (LB) -- (LM);
\draw[flow] (RB) -- (LB);
\draw[flow] (RM) -- (RB);
\draw[flow] (RT) -- (RM);

% Six outward-extending branches
\draw[line] (LT) -- (LTout);
\draw[line] (RTout) -- (RT);
\draw[line] (LMout) -- (LM);
\draw[line] (RM) -- (RMout);
\draw[line] (LB) -- (LBout);
\draw[line] (RBout) -- (RB);
\end{scope}

\begin{scope}[xshift=14cm]

% Circumradius of the regular hexagon
\def\R{1.65}
\pgfmathsetmacro{\H}{sqrt(3)*\R/2}

% Vertices of the regular hexagon
\coordinate (LT) at (-0.5*\R, \H);
\coordinate (RT) at ( 0.5*\R, \H);
\coordinate (RM) at ( \R,0);
\coordinate (RB) at ( 0.5*\R,-\H);
\coordinate (LB) at (-0.5*\R,-\H);
\coordinate (LM) at (-\R,0);

% Outer endpoints
\coordinate (LTout) at (-1*\R, 2*\H);
\coordinate (RTout) at ( 1*\R, 2*\H);
\coordinate (RMout) at ( 2*\R,0);
\coordinate (RBout) at ( 1*\R,-2*\H);
\coordinate (LBout) at (-1*\R,-2*\H);
\coordinate (LMout) at (-2*\R,0);

% Directed edges of the hexagon
\draw[flow] (LT) -- (RT)
  node[pos=.52,above=5pt] {$i$};

\draw[flow] (LM) -- (LT)
  node[pos=.50,left=4pt] {$i+1$};

\draw[flow] (LB) -- (LM)
  node[pos=.54,left=4pt] {$i$};

\draw[flow] (RB) -- (LB)
  node[pos=.50,below=5pt] {$i+1$};

\draw[flow] (RM) -- (RB)
  node[pos=.50,right=4pt] {$i$};

\draw[flow] (RT) -- (RM)
  node[pos=.50,right=4pt] {$i+1$};

% Six outward-extending branches
\draw[flow] (LT) -- (LTout)
  node[pos=.93,left=3pt] {$1$};

\draw[flow] (RTout) -- (RT)
  node[pos=.07,right=3pt] {$1$};

\draw[flow] (LMout) -- (LM)
  node[pos=.08,below=4pt] {$1$};

\draw[flow] (RM) -- (RMout)
  node[pos=.92,below=4pt] {$1$};

\draw[flow] (LB) -- (LBout)
  node[pos=.93,left=3pt] {$1$};

\draw[flow] (RBout) -- (RB)
  node[pos=.07,right=3pt] {$1$};

\end{scope}

\begin{scope}[xshift=21cm]

% Circumradius of the regular hexagon
\def\R{1.65}
\pgfmathsetmacro{\H}{sqrt(3)*\R/2}

\node at (0,0) {$2i+1$};

% Vertices of the regular hexagon
\coordinate (LT) at (-0.5*\R, \H);
\coordinate (RT) at ( 0.5*\R, \H);
\coordinate (RM) at ( \R,0);
\coordinate (RB) at ( 0.5*\R,-\H);
\coordinate (LB) at (-0.5*\R,-\H);
\coordinate (LM) at (-\R,0);

% Outer endpoints
\coordinate (LTout) at (-1*\R, 2*\H);
\coordinate (RTout) at ( 1*\R, 2*\H);
\coordinate (RMout) at ( 2*\R,0);
\coordinate (RBout) at ( 1*\R,-2*\H);
\coordinate (LBout) at (-1*\R,-2*\H);
\coordinate (LMout) at (-2*\R,0);

% Directed edges of the hexagon
\draw[flow] (RT) -- (LT);
\draw[flow] (LT) -- (LM);
\draw[flow] (LM) -- (LB);
\draw[flow] (LB) -- (RB);
\draw[flow] (RB) -- (RM);
\draw[flow] (RM) -- (RT);
% \draw[flow] (LT) -- (RT);
% \draw[flow] (LM) -- (LT);
% \draw[flow] (LB) -- (LM);
% \draw[flow] (RB) -- (LB);
% \draw[flow] (RM) -- (RB);
% \draw[flow] (RT) -- (RM);

% Six outward-extending branches
\draw[line] (LT) -- (LTout);
\draw[line] (RTout) -- (RT);
\draw[line] (LMout) -- (LM);
\draw[line] (RM) -- (RMout);
\draw[line] (LB) -- (LBout);
\draw[line] (RBout) -- (RB);
\end{scope}

%==================================================
% Right diagram
%==================================================
\begin{scope}[xshift=28cm]

% Circumradius of the regular hexagon
\def\R{1.65}
\pgfmathsetmacro{\H}{sqrt(3)*\R/2}

% Vertices of the regular hexagon
\coordinate (LT) at (-0.5*\R, \H);
\coordinate (RT) at ( 0.5*\R, \H);
\coordinate (RM) at ( \R,0);
\coordinate (RB) at ( 0.5*\R,-\H);
\coordinate (LB) at (-0.5*\R,-\H);
\coordinate (LM) at (-\R,0);

% Outer endpoints
\coordinate (LTout) at (-1*\R, 2*\H);
\coordinate (RTout) at ( 1*\R, 2*\H);
\coordinate (RMout) at ( 2*\R,0);
\coordinate (RBout) at ( 1*\R,-2*\H);
\coordinate (LBout) at (-1*\R,-2*\H);
\coordinate (LMout) at (-2*\R,0);

% Directed edges of the hexagon
\draw[flow] (RT) -- (LT)
  node[pos=.52,above=5pt] {$i+1$};

\draw[flow] (LT) -- (LM)
  node[pos=.50,left=4pt] {$i$};

\draw[flow] (LM) -- (LB)
  node[pos=.54,left=4pt] {$i+1$};

\draw[flow] (LB) -- (RB)
  node[pos=.50,below=5pt] {$i$};

\draw[flow] (RB) -- (RM)
  node[pos=.50,right=4pt] {$i+1$};

\draw[flow] (RM) -- (RT)
  node[pos=.50,right=4pt] {$i$};

% Six outward-extending branches
\draw[flow] (LT) -- (LTout)
  node[pos=.93,left=3pt] {$1$};

\draw[flow] (RTout) -- (RT)
  node[pos=.07,right=3pt] {$1$};

\draw[flow] (LMout) -- (LM)
  node[pos=.08,above=4pt] {$1$};

\draw[flow] (RM) -- (RMout)
  node[pos=.92,above=4pt] {$1$};

\draw[flow] (LB) -- (LBout)
  node[pos=.93,left=3pt] {$1$};

\draw[flow] (RBout) -- (RB)
  node[pos=.07,right=3pt] {$1$};

\end{scope}
\end{tikzpicture}}
\caption{A movie depiction of the foam $\xi_i^1 \circ \xi_i^2$.}
\label{fig:xi}
\end{figure}
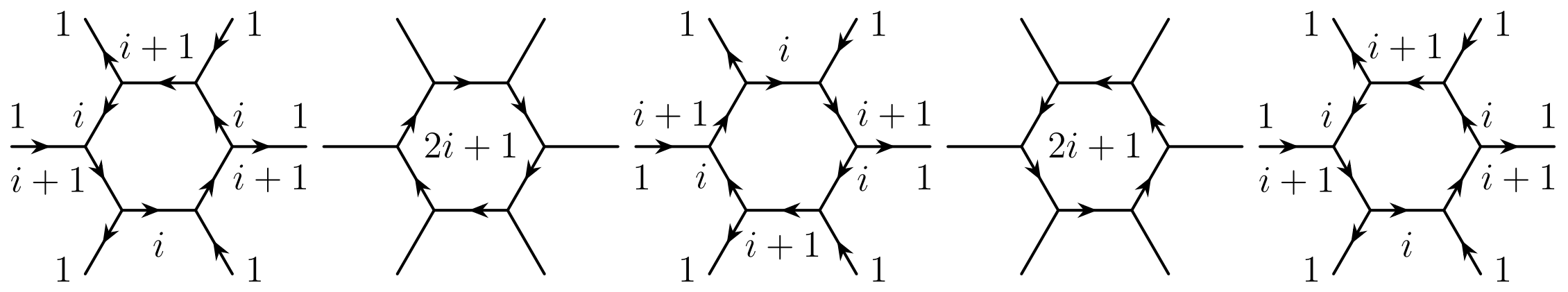
\end{remark}

\begin{lemma}\label{lem:munu = n-3 D D n-3}
    The following foam equality holds:
    \begin{align*}
         \mu \cdot \nu = F_{N-2}^1 D^2 F_{1}^{N-2}.
    \end{align*}
\end{lemma}
\begin{proof}
  The same sequence of rotations and reflections carries $\nu$ to $\mu$ and $D F_1^{N-2}$ to $F_{N-2}^1 D$. Under the natural correspondence of coloringss, the bichrome surfaces of $\nu$ and $D F_1^{N-2}$ have the same Euler characteristics. The result then follows by the argument used in the proof of \cite[Lemma 3.21]{RW2020}. 
\end{proof}

Combining the preceding results gives the following theorem.

\begin{theorem}[Conjecture 4.29 in \cite{grlj20266valentvertexmathfrakglnweb}]
For $N \geq 4$, there exists a decomposition of $\Gamma_1$ in $\Kar(\Foam^{\varepsilon})$ given as: 
\[
 \mathcal{E}(\Gamma_1) \cong q^{4-N}\mathcal{E}(\Gamma_0) \oplus q^{6-N}\mathcal{E}(\Gamma_0)  \oplus \dots \oplus q^{N-4}\mathcal{E}(\Gamma_0) \oplus \mathcal{E}(\mu \cdot \nu).
\]
\end{theorem}
\begin{proof}
The result follows from Proposition \ref{lem:munu as Ri} if we show that Equation \eqref{eq:six vertex categorification} can be expressed in the Karoubi envelope as 
\begin{align}\label{eq:decategorification without munu}
     \mathcal{E}(\Gamma_1) \cong q^{4-N}\mathcal{E}(\Gamma_0) \oplus q^{6-N}\mathcal{E}(\Gamma_0)  \oplus \dots \oplus q^{N-4}\mathcal{E}(\Gamma_0) \oplus \mathcal{E}\left(\prod_{s = 0}^{N-4}  \Ri{s}{1} \right).
\end{align}

Compute the following degrees using Lemma \ref{lem:deg compuation} (here $|B| = 3$):
\begin{align*}
    \foamdeg(\overline{F_{1}^{0}}(h_s)) &= 2s - N +1 + 3 = 2s - N + 4\\
    \foamdeg(F_{0}^1(h_{N - 4-s} )) & = 2(N - 4-s) - N +1 + 3 = - 2s + N - 4
\end{align*}

Recall from Definition \ref{def:Li,Ri} that:
\begin{align*}
    \Li{s}{1}&= F_{0}^1(h_{N - 4-s} )\overline{F_{1}^{0}}(h_s), \\
    \Ri{s}{0} &=\overline{F_{1}^{0}}(h_s)F_{0}^1(h_{N - 4-s} ).
\end{align*}
Lemma \ref{lem:Li,Ri prop} and Lemma \ref{lem:prod of Ri Li vanishing} show that $\Ri{s}{0} + \Li{s}{0} = \idfoam_{\Gamma_0}$ and $\Li{s}{0} = 0$. Hence,
\begin{align*}
    \Ri{s}{0} = \idfoam_{\Gamma_0}.
\end{align*}
That is $F_{0}^1(h_{N - 4-s} )$ and $\overline{F_{1}^{0}}(h_s)$ are inverses in $\Kar(\Foam^\varepsilon),$ demonstrating:
\begin{align*}
    \mathcal E(\Li{s}{1}) \cong q^{N - 4 - 2s}\mathcal E (\Gamma_0).
\end{align*}
An application of Lemma~\ref{lem:six vertex categorification} completes the proof.
\end{proof}

\section{Indecomposability of $2m$-valent vertices}\label{sec:indecomposabilis 2m vertex}
Throughout this section we work over $\C(q)$, fix $m$, and assume
\begin{align*}
    N\geq 2,\qquad 0\leq j<i,\qquad i+j<N,\qquad a =(0,\ldots,0), \qquad b = (1, \ldots, 1). 
\end{align*}

We will prove that the $2m$-valent vertices (see Definition \ref{def:2m valent vertex}) are indecomposable.

\subsection{Auxiliary computations in $\dot{\mathcal U}_q^{N}(\GL_2)$}

In order to perform the web evaluations required to determine indecomposability, we employ Skew-Howe duality from \cite{CKM2014}. This section carries out the computations in $\dot{\mathcal U}_q^{N}(\GL_2)$ that will be utilized in the next section. We summarize the definition of $\dot{\mathcal U}_q^{N}(\GL_2)$ and provide some equalities within.

The quantum group $\dot{\mathcal U}_q^{N}(\GL_2)$ from \cite[Section 4.1]{CKM2014} has objects given by pairs $(a,b)$ with $0\leq a,b\leq N$, objects outside these bounds are zero. We write $1_{(a,b)}$ for the identity of $(a,b)$. The generators $E,F$ act as
\begin{align*}
    E1_{(a,b)} &: (a,b)\longrightarrow(a+1,b-1),\\
    F1_{(a,b)} &: (a,b)\longrightarrow(a-1,b+1).
\end{align*}
Their divided powers are
\begin{align*}
    E^{(k)}=\frac{E^k}{[k]_q!},\qquad
    F^{(k)}=\frac{F^k}{[k]_q!},\qquad
    [k]_q!=\prod_{a=1}^k[a]_q.
\end{align*}

\begin{lemma}\label{lem:two-m triangular}
For $0\leq k\leq j$,
\begin{align}\label{eq:two-m triangular}
    FE E^kF^k1_{(i,N-j)}
    ={}&E^{k+1}F^{k+1}1_{(i,N-j)}\notag\\
    &+[k+1]_q[N-i-j+k]_q E^kF^k1_{(i,N-j)}.
\end{align}
Moreover, $F^{j+1}1_{(i,N-j)}=0$. In particular,
\begin{align}\label{eq:two-m Newton}
    E^kF^k1_{(i,N-j)}
    &=\prod_{s=0}^{k-1}\bigl(FE1_{(i,N-j)}-[s+1]_q[N-i-j+s]_q 1_{(i,N-j)}\bigr),\notag\\
    0&=\prod_{s=0}^j\bigl(FE1_{(i,N-j)}-[s+1]_q[N-i-j+s]_q1_{(i,N-j)}\bigr).
\end{align}
\end{lemma}

\begin{proof}
Relation \cite[Equation (4.1)]{CKM2014} implies
\begin{align*}
    FE^{k+1}1_{(a,b)}
    =E^{k+1}F1_{(a,b)}-[k+1]_q[a-b+k]_q E^k1_{(a,b)}.
\end{align*}
$F^k$ changes weights as $F^k 1_{(i,N-j)} = 1_{(i-k,N-j+k)} F^k 1_{(i,N-j)}$. Substituting into the preceding identity proves Equation~\eqref{eq:two-m triangular}. Iterating Equation~\eqref{eq:two-m triangular} proves Equation~\eqref{eq:two-m Newton}, since $(i-j-1,N+1)$ is the zero object in $\dot{\mathcal U}_q^{N}(\GL_2)$. 
\end{proof}

\begin{proposition}\label{prop:two-m powers}
For every $m\geq 0$,
\begin{align*}
    &(FE)^m1_{(i,N-j)}
    =\sum_{k=0}^{\min(m,j)}
    \left(
    \sum_{\substack{a_0,\ldots,a_k\geq 0\\a_0+\cdots+a_k=m-k}}
    \prod_{s=0}^k\bigl([s+1]_q[N-i-j+s]_q\bigr)^{a_s}
    \right)E^kF^k1_{(i,N-j)}  \\
    &=\sum_{k=0}^{\min(m,j)}
h_{m-k}\left(
    [N-i-j]_q,\,
    [2]_q[N-i-j+1]_q,\,
    \ldots,\,
    [k+1]_q[N-i-j+k]_q
\right)E^kF^k1_{(i,N-j)}.
\end{align*} 
\end{proposition}

\begin{proof}
The result follows by induction using Equation~\eqref{eq:two-m triangular}, and the fact that terms with $k=j+1$ (which appear in the induction step) vanish by Lemma \ref{lem:two-m triangular}.
\end{proof}

\subsection{Web pairings and indecomposability.}

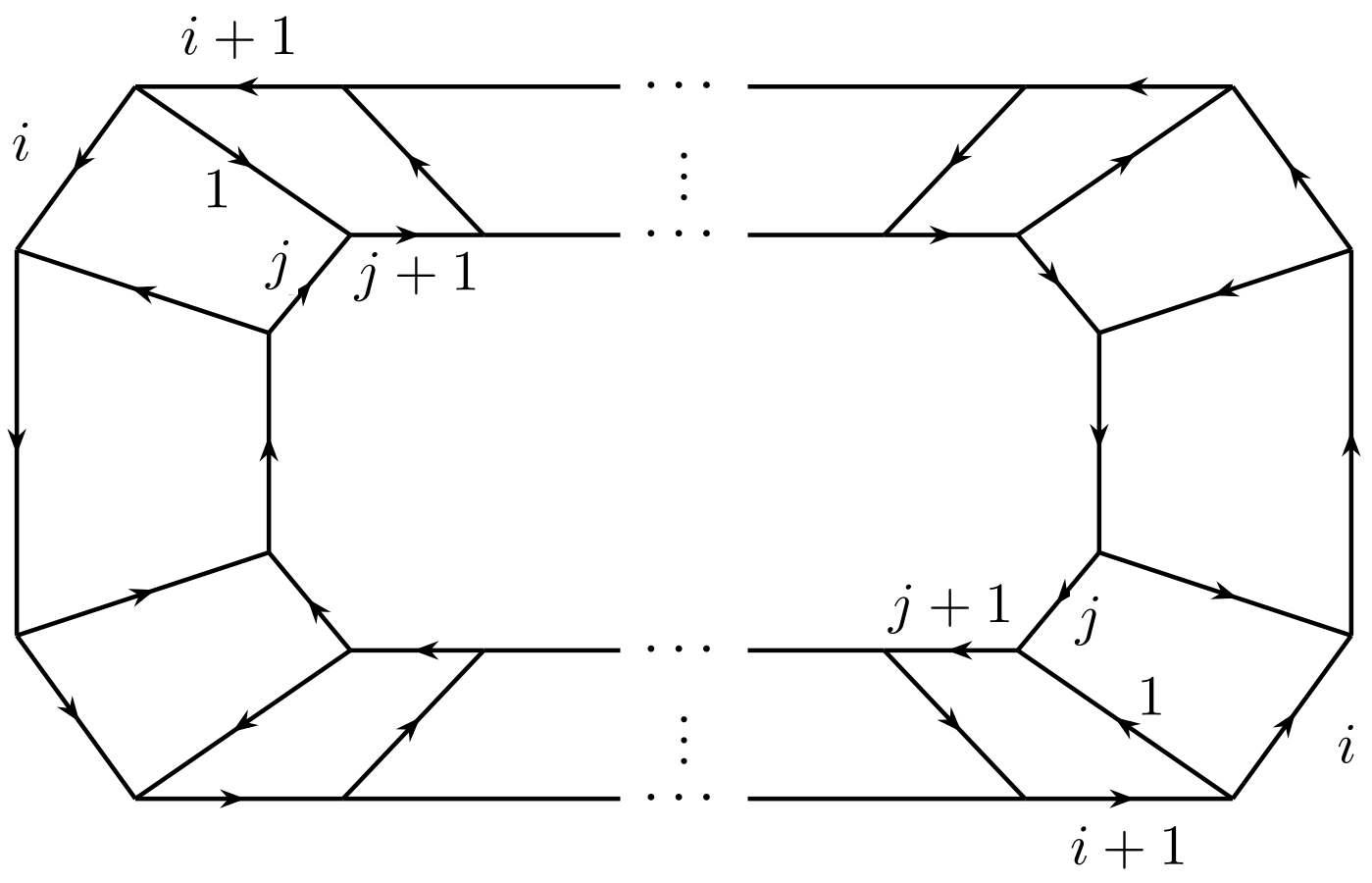
\begin{figure}[htbp!]
\centering
\tikzset{
  outer cycle/.style={
    draw=black,
    line width=.85pt,
    postaction={decorate},
    decoration={markings,mark=at position .53 with
      {\arrow{Stealth[length=1.7mm,width=1.25mm]}}}
  },
  inner cycle/.style={
    draw=black,
    %blue!65!black,
    line width=.85pt,
    postaction={decorate},
    decoration={markings,mark=at position .53 with
      {\arrow{Stealth[length=1.7mm,width=1.25mm]}}}
  },
  radial edge/.style={
    draw=black,
    line width=.85pt,
    postaction={decorate},
    decoration={markings,mark=at position .55 with
      {\arrow{Stealth[length=1.7mm,width=1.25mm]}}}
  },
  cycle label/.style={font=\small,inner sep=1pt,fill=white}
}
\begin{tikzpicture}[x=1cm,y=1cm]
  % Visible vertices, in counterclockwise order.  The two omitted
  % stretches at the top and bottom contain the remaining vertices.
  \foreach \k/\x/\y in {
     0/-2.3/ 2.4, 1/-3.7/ 2.4, 2/-4.5/ 1.3,
     3/-4.5/-1.3, 4/-3.7/-2.4, 5/-2.3/-2.4,
     6/ 2.3/-2.4, 7/ 3.7/-2.4, 8/ 4.5/-1.3,
     9/ 4.5/ 1.3,10/ 3.7/ 2.4,11/ 2.3/ 2.4}
    \coordinate (O\k) at (\x,\y);
  \foreach \k/\x/\y in {
     0/-1.35/ 1.4, 1/-2.25/ 1.4, 2/-2.8/ .74,
     3/-2.8/-.74, 4/-2.25/-1.4, 5/-1.35/-1.4,
     6/ 1.35/-1.4, 7/ 2.25/-1.4, 8/ 2.8/-.74,
     9/ 2.8/ .74,10/ 2.25/ 1.4,11/ 1.35/ 1.4}
    \coordinate (I\k) at (\x,\y);

  % Gamma_i is counterclockwise.  Reflection reverses the arrows on
  % the inner copy, representing overline{Gamma_j}.
  \foreach \k in {0,1,2,3,4,6,7,8,9,10}{
    \pgfmathtruncatemacro{\kp}{\k+1}
    \draw[outer cycle] (O\k) -- (O\kp);
    \draw[inner cycle] (I\kp) -- (I\k);
  }

  % At consecutive vertices the 1-labeled radial edges alternate
  % direction, satisfying the flow rule for i/i+1 and j/j+1.
  \foreach \k in {0,2,4,6,8,10}
    \draw[radial edge] (I\k) -- (O\k);
  \foreach \k in {1,3,5,7,9,11}
    \draw[radial edge] (O\k) -- (I\k);

  % Broken portions of both cycles; each column of vertical dots is
  % an omitted family of radial 1-edges.
  \foreach \y/\a/\b in {2.4/2.3/.43,-2.4/2.3/.43}{
    \draw[line width=.85pt] (-\a,\y) -- (-\b,\y);
    \draw[line width=.85pt] (\b,\y) -- (\a,\y);
    \node[fill=white,inner sep=0pt] at (0,\y) {$\cdots$};
  }
  \foreach \y/\a/\b in {1.4/1.35/.43,-1.4/1.35/.43}{
    \draw[draw=black,line width=.85pt] (-\a,\y) -- (-\b,\y);
    \draw[draw=black,line width=.85pt] (\b,\y) -- (\a,\y);
    \node[fill=white,inner sep=0pt] at (0,\y) {$\cdots$};
  }
  \node[fill=white,inner sep=0pt] at (0, 1.9) {$\vdots$};
  \node[fill=white,inner sep=0pt] at (0,-1.9) {$\vdots$};

  % Representative labels; all other labels follow by alternation.
  \node[cycle label] at (-3.0, 2.73) {$i+1$};
  \node[cycle label] at (-4.47, 2.03) {$i$};
  \node[cycle label] at ( 3.0,-2.73) {$i+1$};
  \node[cycle label] at ( 4.47,-2.03) {$i$};

  \node[cycle label,text=black] at (-1.80,1.12) {$j+1$};
  \node[cycle label,text=black] at (-2.73,1.20) {$j$};
  \node[cycle label,text=black] at ( 1.80,-1.12) {$j+1$};
  \node[cycle label,text=black] at ( 2.73,-1.20) {$j$};

  \node[cycle label] at (-3.15,1.71) {$1$};
  \node[cycle label] at ( 3.15,-1.71) {$1$};
\end{tikzpicture}
\caption{The web $\Gamma_i^m\overline{\Gamma_j^m}$, there are $2m$-many edges connecting the two polygons.}
\label{fig: web GammaiGammaj}
\end{figure}

In this section we evaluate $\langle \Gamma_i^m\overline{\Gamma_j^m} \rangle \in \C(q)$, where the superscript $m$ denotes the dependence on $m$ (see Figure \ref{fig: web GammaiGammaj}). Expanding the evaluation in $q$, will allow us to deduce indecombosability of $2m$-valent vertices.

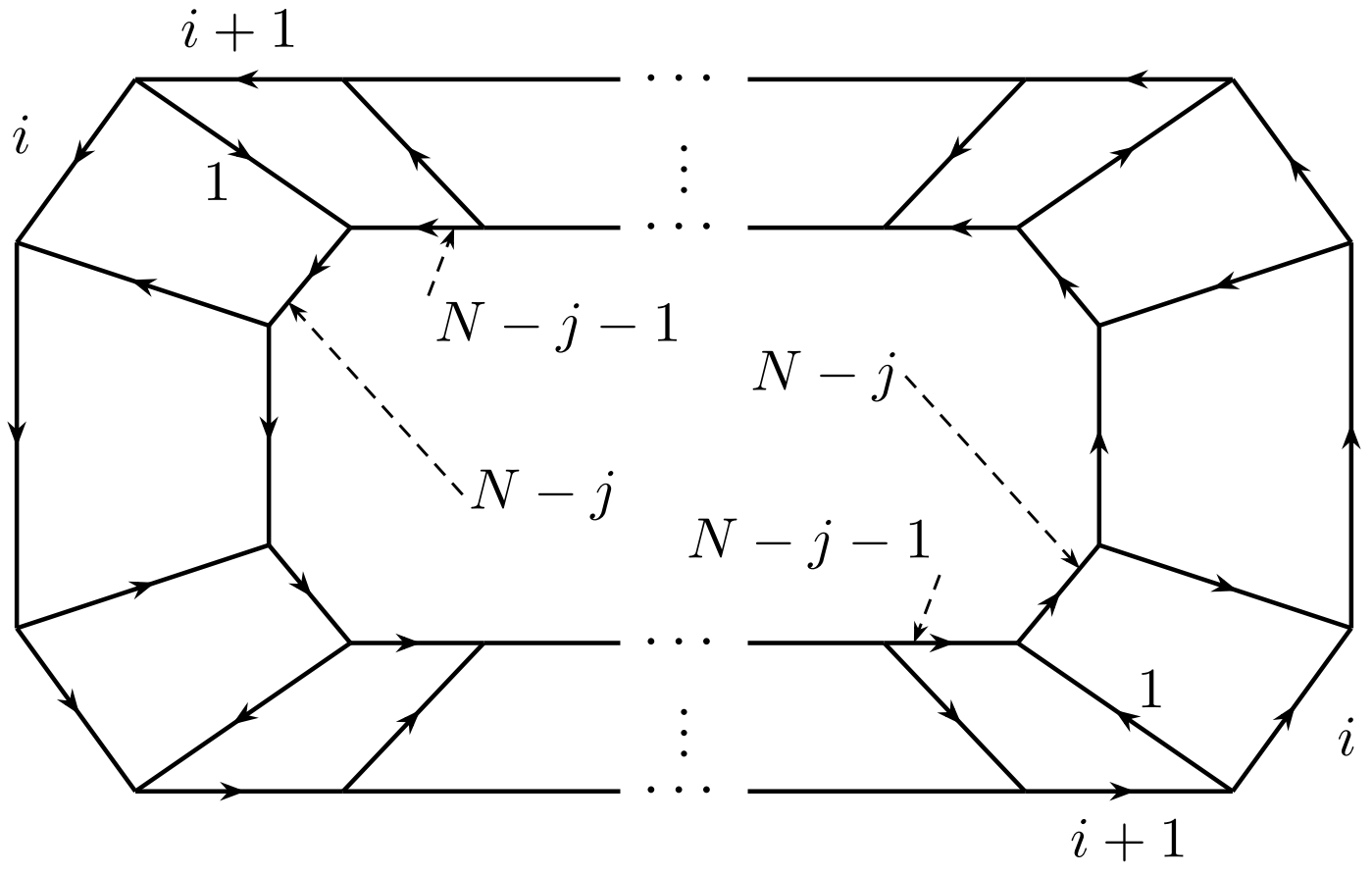
\begin{figure}[htbp!]
\centering
% Schematic closure of the two Kekule webs from main(34).tex when
% a=(0,...,0) and b=(1,...,1).  The dots suppress repeated pairs
% of trivalent vertices; every radial edge has label 1.
\tikzset{
  outer cycle/.style={
    draw=black,
    line width=.85pt,
    postaction={decorate},
    decoration={markings,mark=at position .53 with
      {\arrow{Stealth[length=1.7mm,width=1.25mm]}}}
  },
  inner cycle/.style={
    draw=black,
    line width=.85pt,
    postaction={decorate},
    decoration={markings,mark=at position .53 with
      {\arrow{Stealth[length=1.7mm,width=1.25mm]}}}
  },
  radial edge/.style={
    draw=black,
    line width=.85pt,
    postaction={decorate},
    decoration={markings,mark=at position .55 with
      {\arrow{Stealth[length=1.7mm,width=1.25mm]}}}
  },
  label pointer/.style={
    draw=black,
    line width=.55pt,
    densely dashed,
    -{Stealth[length=1.4mm,width=1mm]}
  },
  cycle label/.style={font=\small,inner sep=1pt,fill=white}
}
\begin{tikzpicture}[x=1cm,y=1cm]
  % Visible vertices, in counterclockwise order.  The two omitted
  % stretches at the top and bottom contain the remaining vertices.
  \foreach \k/\x/\y in {
     0/-2.3/ 2.4, 1/-3.7/ 2.4, 2/-4.5/ 1.3,
     3/-4.5/-1.3, 4/-3.7/-2.4, 5/-2.3/-2.4,
     6/ 2.3/-2.4, 7/ 3.7/-2.4, 8/ 4.5/-1.3,
     9/ 4.5/ 1.3,10/ 3.7/ 2.4,11/ 2.3/ 2.4}
    \coordinate (O\k) at (\x,\y);
  \foreach \k/\x/\y in {
     0/-1.35/ 1.4, 1/-2.25/ 1.4, 2/-2.8/ .74,
     3/-2.8/-.74, 4/-2.25/-1.4, 5/-1.35/-1.4,
     6/ 1.35/-1.4, 7/ 2.25/-1.4, 8/ 2.8/-.74,
     9/ 2.8/ .74,10/ 2.25/ 1.4,11/ 1.35/ 1.4}
    \coordinate (I\k) at (\x,\y);

  % The outer polygon is counterclockwise.  The inner polygon is the
  % dual of the reflected web: reverse its arrows and complement its
  % internal labels from ell to N-ell.
  \foreach \k in {0,1,2,3,4,6,7,8,9,10}{
    \pgfmathtruncatemacro{\kp}{\k+1}
    \draw[outer cycle] (O\k) -- (O\kp);
    \draw[inner cycle] (I\k) -- (I\kp);
  }

  % At consecutive vertices the 1-labeled radial edges alternate
  % direction, satisfying the flow rule after dualizing the inner cycle.
  \foreach \k in {0,2,4,6,8,10}
    \draw[radial edge] (I\k) -- (O\k);
  \foreach \k in {1,3,5,7,9,11}
    \draw[radial edge] (O\k) -- (I\k);

  % Broken portions of both cycles; each column of vertical dots is
  % an omitted family of radial 1-edges.
  \foreach \y/\a/\b in {2.4/2.3/.43,-2.4/2.3/.43}{
    \draw[line width=.85pt] (-\a,\y) -- (-\b,\y);
    \draw[line width=.85pt] (\b,\y) -- (\a,\y);
    \node[fill=white,inner sep=0pt] at (0,\y) {$\cdots$};
  }
  \foreach \y/\a/\b in {1.4/1.35/.43,-1.4/1.35/.43}{
    \draw[draw=black,line width=.85pt] (-\a,\y) -- (-\b,\y);
    \draw[draw=black,line width=.85pt] (\b,\y) -- (\a,\y);
    \node[fill=white,inner sep=0pt] at (0,\y) {$\cdots$};
  }
  \node[fill=white,inner sep=0pt] at (0, 1.9) {$\vdots$};
  \node[fill=white,inner sep=0pt] at (0,-1.9) {$\vdots$};

  % Representative labels; all other labels follow by alternation.
  \node[cycle label] at (-3.0, 2.73) {$i+1$};
  \node[cycle label] at (-4.47, 2.03) {$i$};
  \node[cycle label] at ( 3.0,-2.73) {$i+1$};
  \node[cycle label] at ( 4.47,-2.03) {$i$};

  % Each dashed pointer ends at the inner edge it labels.
  \node[cycle label] (tl-long) at (-.85, .73) {$N-j-1$};
  \node[cycle label] (tl-slant) at (-.95,-.40) {$N-j$};
  \node[cycle label] (br-long) at ( .85,-.73) {$N-j-1$};
  \node[cycle label] (br-slant) at ( .95, .40) {$N-j$};
  \draw[label pointer] (tl-long.north west) -- (-1.55, 1.40);
  \draw[label pointer] (tl-slant.west) -- (-2.67, .90);
  \draw[label pointer] (br-long.south east) -- (1.55,-1.40);
  \draw[label pointer] (br-slant.east) -- (2.67,-.90);

  \node[cycle label] at (-3.15,1.71) {$1$};
  \node[cycle label] at ( 3.15,-1.71) {$1$};
\end{tikzpicture}
\caption{The web $\Gamma_i^m\overline{\Gamma_j^m}'$.}
\label{fig: web GammaiGammaj '}
\end{figure}

First step in completing this evaluation is ladderizing  $\Gamma_i^m\overline{\Gamma_j^m}$. To do so, we use the dual web $\overline{\Gamma_j^m}'$ of $\overline{\Gamma_j^m}$ (see Definition \ref{def:duality morphism}), this gives us the web displayed in Figure \ref{fig: web GammaiGammaj '}.

\begin{lemma}
    It holds that $\langle \Gamma_i^m\overline{\Gamma_j^m} \rangle = \langle \Gamma_i^m\overline{\Gamma_j^m}' \rangle$.
\end{lemma}
\begin{proof}
    This is a well known fact, see \cite{CKM2014, Mor2007}.
\end{proof}

Before applying the results from the previous section we require an auxiliary computation which also computes $\Gamma_i^2\overline{\Gamma_j^2}'$.

\begin{proposition}\label{prop:two-m basis closure}
For $0\leq k\leq j$, the following closed web evaluates as
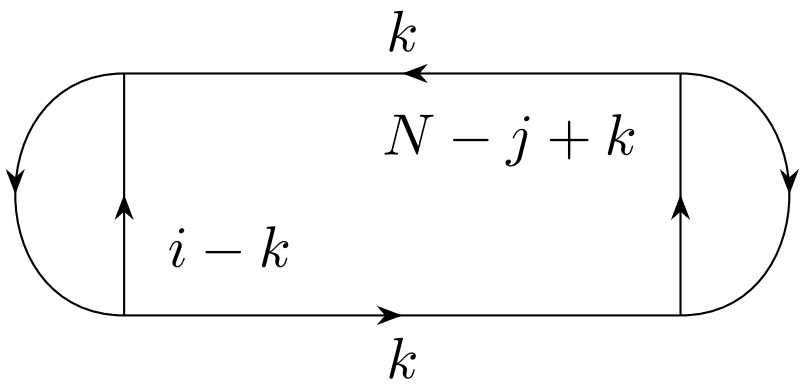
\begin{equation}\label{eq:two-m basis closure}
    \left\langle
   \begin{tikzpicture}[
    scale=0.8,
    baseline={([yshift=-.5ex]current bounding box.center)},
    every node/.style={font=\small},
    decoration={markings,mark=at position 0.5 with {\arrow{Stealth}}}
]
    \draw[postaction={decorate}] (-2.3,-1)--(-2.3,1);
    \draw[postaction={decorate}] ( 2.3,-1)--( 2.3,1);
    \draw[postaction={decorate}] ( 2.3,1)--(-2.3,1);
    \draw[postaction={decorate}] (-2.3,-1)--( 2.3,-1);
    \draw[postaction={decorate}]
        (-2.3,1)..controls (-3.5,1) and (-3.5,-1)..(-2.3,-1);
    \draw[postaction={decorate}]
        (2.3,1)..controls (3.5,1) and (3.5,-1)..(2.3,-1);

    \node[above] at (0,1) {$k$};
    \node[below] at (0,-1) {$k$};
    \node[right=4pt] at (-2.3,-0.45) {$i-k$};
    \node[left=4pt]  at ( 2.3, 0.45) {$N-j+k$};
    \node[left=4pt]  at (-3.2,0) {$i$};
    \node[right=4pt] at ( 3.2,0) {$N-j$};
\end{tikzpicture}
    \right\rangle
    =
    \begin{bmatrix}N\\k\end{bmatrix}_q
    \begin{bmatrix}N-k\\i-k\end{bmatrix}_q
    \begin{bmatrix}N-k\\j-k\end{bmatrix}_q.
\end{equation}
The pictured web is the planar closure of the ladder web corresponding to $E^{(k)}F^{(k)}1_{(i,N-j)}$. 
\end{proposition}

\begin{proof}
This is a simple web computation done by applying the (diagon and circle) relations from Figure~\ref{relations}.
\end{proof}

Combining Proposition \ref{prop:two-m basis closure} and the results from the previous section yields the following.

\begin{theorem}\label{prop:two-m closed powers}
For every $m\geq 0$,
\begin{align}\label{eq:computing the pairing}
 \bigl\langle \Gamma_i^m &\overline{\Gamma_j^m}\bigr\rangle
    =\bigl\langle \Gamma_i^m\overline{\Gamma_j^m}'\bigr\rangle  = \sum_{k=0}^{\min(m,j)}([k]_q!)^2
    \begin{bmatrix}N\\k\end{bmatrix}_q
    \begin{bmatrix}N-k\\i-k\end{bmatrix}_q
    \begin{bmatrix}N-k\\j-k\end{bmatrix}_q\notag\\
    &\cdot
    h_{m-k}\left(
    [N-i-j]_q,\,
    [2]_q[N-i-j+1]_q,\,
    \ldots,\,
    [k+1]_q[N-i-j+k]_q
\right).
\end{align}
\end{theorem}

\begin{proof}
This is a direct consequence of Propositions~\ref{prop:two-m powers} and  \ref{prop:two-m basis closure}, and the definition of divided powers in $\dot{\mathcal U}_q^{N}(\GL_2)$. 
\end{proof}

The results from this section suffice to determine indecomposability of $2m$-valent vertices. First we recall a condition which determines indecomposability of a web by analyzing its self pairing.

\begin{lemma}[Lemma 2.43 in \cite{grlj20266valentvertexmathfrakglnweb}] \label{WebPairIndecompo}
    A \(\GL_N\) web \(W\) with boundary \(\varepsilon\) is indecomposable if
    \[
        q^{\boundarydeg(\varepsilon)}\langle \overline{W}W \rangle \in 1+q\mathbb{Z}[q].
    \]
    Equivalently, \(\langle \overline{W}W \rangle\) has coefficient \(1\) in degree \(-\boundarydeg(\varepsilon)\) and has no terms in smaller degrees.
\end{lemma}

The following generalizes Definition \ref{def:six valent} by defining the $2m$-valent vertices for any choice of $m$. It should be compared to Proposition \ref{lem:karoubi categorification}, in particular, $\mathcal{W}^j_{2m}$ corresponds to $\mathcal E(\mathfrak{E}_{j,\emptyset})$ in the Karoubi envelope. Note that for $m > 3$, there are multiple non equivalent $2m$-valent vertices. Recall that for the choice of sequences $a$ and $b$ in this section, it holds that $l = N-1$ and $r = N - 1 -m$.

\begin{definition}\label{def:2m valent vertex}
    For $j\in \{1, \ldots, m-2\}$, define \emph{$2m$-valent vertices}  $ \mathcal{W}^j_{2m} \in \Web$  as
    \begin{align*}
        \mathcal{W}^j_{2m} := \sum_{k = 0}^j (-1)^k\begin{bmatrix}
            r + k \\
            k
        \end{bmatrix}_q\Gamma_{j-k}^m.
    \end{align*}
\end{definition}

If we set $m= 3$ and $j = 1$, this recovers the definition of $\Wsix$. It was shown in \cite[Proposition 3.8]{grlj20266valentvertexmathfrakglnweb} that $\Wsix$ is indecomposable. The same holds for all other vertices from Definition \ref{def:2m valent vertex}.

\begin{corollary}\label{prop:two-m extreme coefficients}
Let $m \geq 3$ and $j\in \{1, \ldots, m-2\}$. Then the vertex $\mathcal{W}^j_{2m}$ is  indecomposable.
\end{corollary}

\begin{proof}
By Lemma \ref{WebPairIndecompo}, it suffices to prove that
\begin{align*}
    P(q) := q^{m(N-1)} \left\langle \mathcal{W}^j_{2m}\overline{\mathcal{W}^j_{2m}} \right\rangle \in 1 + q\Z[q].
\end{align*}

By Definition \ref{def:2m valent vertex} and Theorem \ref{prop:two-m closed powers}, we have:
\begin{align}\label{eq:two-m direct-expansion}
    P&(q) = q^{m(N-1)} \sum_{k = 0}^j \sum_{c,d = 0}^{j-k} (-1)^{c+d} ([k]_q!)^2
    \begin{bmatrix}N\\k\end{bmatrix}_q  \begin{bmatrix}N-m+c-1\\c\end{bmatrix}_q
    \begin{bmatrix}N-m+d-1\\d\end{bmatrix}_q \\
    & \cdot
    \begin{bmatrix}N-k\\j-c-k\end{bmatrix}_q
    \begin{bmatrix}N-k\\j-d-k\end{bmatrix}_q h_{m-k}\left([N-2j+c+d]_q,\cdots,
    [k+1]_q[N-2j+c+d+k]_q\right) \notag.
\end{align}
For each $k$, the term indexed by $(c,d)$ comes from pairing $ \begin{bmatrix}r+c\\c\end{bmatrix}_q\Gamma_{j-c}^m$ with $ \begin{bmatrix}r+d\\d\end{bmatrix}_q\overline{\Gamma_{j-d}^m}$ and there are terms with negative degrees of $q$. We will show that as we sum over $0 \leq c,d \leq j-k$ terms of negative degree cancel.

For each $0 \leq t \leq k$, write:
\begin{align*}
    \lambda_t &= [t+1]_q[N-2j+c+d+t]_q 
    \\
    &= \frac{[t+1]_q}{1-q^2}\left(q^{1-N-t+2j-c-d} - q^{1+N+t-2j+c+d}\right).
\end{align*}
Using the above, we can expand $h_{m-k}$ as:
\begin{align*}
    h_{m-k}(\lambda_0,\cdots,\lambda_k)
    &= \frac{1}{(1-q^2)^{m-k}}\sum_{\ell = 0}^{m-k} (-1)^\ell
     q^{(m-k-\ell)(1-N+2j-c-d)+\ell(1+N-2j+c+d)} C_{k,\ell}(q), \\
    C_{k,\ell}(q) &:= \sum_{\substack{n_0+\cdots+n_k=m-k\\n_0,\cdots,n_k\geq0}} \sum_{\substack{e_0+\cdots+e_k=\ell\\0\leq e_t\leq n_t}} \prod_{t=0}^k \binom{n_t}{e_t} q^{t(2e_t-n_t)}[t+1]_q^{n_t}.
\end{align*}
In particular, $C_{k,\ell}(q)$ is independent of $c,d$ and has no terms below degree $-2k(m-k-\ell)$. The dependence of this expansion on $c,d$ is contained in $q^{(m-k-2\ell)(2j-c-d)}$. Thus we can write $q^{(m-k-2\ell)(2j-c-d)}$ as $q^{(m-k-2\ell)(j-c)}q^{(m-k-2\ell)(j-d)}$ and the sum over $c,d$ in Equation \eqref{eq:two-m direct-expansion} splits into two sums, one over $c$ and the other one over $d$. It suffices to show that negative terms in one of those sums cancel. Set
\begin{align*}
    S_{k,\ell}(q) := \sum_{c = 0}^{j-k} (-1)^c
    \begin{bmatrix}N-m+c-1\\c\end{bmatrix}_q
    \begin{bmatrix}N-k\\j-k-c\end{bmatrix}_q q^{(m-k-2\ell)(j-c)},
\end{align*}
so that 
\begin{align*}
&\sum_{c,d=0}^{j-k}(-1)^{c+d}
\begin{bmatrix}N-m+c-1\\c\end{bmatrix}_q
\begin{bmatrix}N-m+d-1\\d\end{bmatrix}_q
\begin{bmatrix}N-k\\j-k-c\end{bmatrix}_q
\begin{bmatrix}N-k\\j-k-d\end{bmatrix}_q
q^{(m-k-2\ell)(2j-c-d)}\\
&=\left(
\sum_{c=0}^{j-k}(-1)^c
\begin{bmatrix}N-m+c-1\\c\end{bmatrix}_q
\begin{bmatrix}N-k\\j-k-c\end{bmatrix}_q
q^{(m-k-2\ell)(j-c)}
\right)^2\\
&=S_{k,\ell}(q)^2.
\end{align*}
Then Equation \eqref{eq:two-m direct-expansion} becomes:
\begin{align*}
    P(q) &= \sum_{k = 0}^j \sum_{\ell = 0}^{m-k} P_{k,\ell}(q), \\
    P_{k,\ell}(q) &:= \frac{(-1)^\ell q^{k(N-1)+2N\ell}}{(1-q^2)^{m-k}}
    ([k]_q!)^2 \begin{bmatrix}N\\k\end{bmatrix}_q C_{k,\ell}(q) S_{k,\ell}(q)^2.
\end{align*}
Next, we focus on the lowest coefficient of $S_{k,\ell}$, for ease of reading, set
\begin{align*}
    \alpha_t := N+m-2k-1-2\ell-2t, \qquad \beta_u := N-m-1-2u.
\end{align*}
The inverse quantum binomial theorem and the above conventions give the following expansions of terms in $S_{k,\ell}$:
\begin{align*}
    \begin{bmatrix}N-m+c-1\\c\end{bmatrix}_q
    &= \sum_{\substack{\nu \in \mathbb{Z}_{\ge 0}^{N-m}\\|\nu| = c}}
    q^{\sum_{u = 0}^{N-m-1}\nu_u\beta_u}.
\end{align*}
The rest of the terms can be expanded using the quantum binomial theorem and the definition of $\alpha_t$ as 
\begin{align*}
    \begin{bmatrix}N-k\\j-k-c\end{bmatrix}_q
    &= \sum_{\substack{T \subseteq \{0,\cdots,N-k-1\}\\|T| = j-k-c}}
    q^{\sum_{t \in T} N - k - 1 - 2t}, \\
    q^{(m-k-2\ell)(j-k-c)}\begin{bmatrix}N-k\\j-k-c\end{bmatrix}_q
    &= \sum_{\substack{T \subseteq \{0,\cdots,N-k-1\}\\|T| = j-k-c}}
    q^{\sum_{t \in T}\alpha_t}.
\end{align*}

Since $j-c = k+(j-k-c)$, the $c$-th summand of $S_{k,\ell}$ is therefore the sum of the terms
\begin{align*}
    &(-1)^c q^{k(m-k-2\ell)+\sum_{t \in T}\alpha_t+\sum_{u = 0}^{N-m-1}\nu_u\beta_u}, \\
    &\qquad |T| = j-k-c, \quad |\nu| = c.
\end{align*}

Each pair $(T,\nu)$ specifies one monomial in $P(q)$, $T$ specifics which exponents $\alpha_t$ occur, and $\nu_u$ specifies the number of exponents $\beta_u$. The sign of this monomial is $(-1)^c$, where $c = |\nu|$.

Note that some of $\alpha_t$ equal $\beta_u$. Specifically,
\begin{align*}
    \alpha_{m-k-\ell+u}
    &= N+m-2k-1-2\ell-2(m-k-\ell+u) \\
    &= N-m-1-2u = \beta_u.
\end{align*}

Write $t_u = m-k-\ell+u$, with $0 \leq u \leq N-m-1$. Whenever $t_u \in T$, removing $t_u$ from $T$ and increasing $\nu_u$ by $1$ leaves the power of $q$ unchanged and replaces $c$ by $c+1$. The sign changes, so the two terms cancel. 
% More precisely, choose the smallest $u$ for which $t_u \in T$ or $\nu_u > 0$, and pair
% \begin{align*}
%     (T,\nu) \longleftrightarrow \left(T \setminus \{t_u\},\nu+e_u\right)
% \end{align*}
% if $t_u \in T$, where $e_u$ is the $u$-th unit vector. If $t_u \notin T$ and $\nu_u > 0$, make the reverse change. The same $u$ is selected after the change, so applying the rule twice restores the original term.

Therefore, only the terms corresponding to $\nu = 0$ and $T$ disjoint from $\{m-k-\ell,\cdots,N-k-\ell-1\}$ remain in $S_{k,\ell}$. These terms have $c = 0$, so $|T| = j-k$ and their sign is positive. Hence,
\begin{align}\label{eq:two-m direct-cancellation}
    S_{k,\ell}(q) = q^{k(m-k-2\ell)}
    \sum_{\substack{T \subseteq \{0,\cdots,N-k-1\}\\|T| = j-k\\
        T \cap \{m-k-\ell,\cdots,N-k-\ell-1\} = \varnothing}}
    q^{\sum_{t \in T}\alpha_t}.
\end{align}

Every subset $T$ in Equation \eqref{eq:two-m direct-cancellation} contains $j-k$ indices. To minimize $\sum_{t \in T}\alpha_t$, we therefore choose the $j-k$ smallest remaining exponents. Their sum is:
\begin{align*}
    d_{k,\ell} = \begin{cases}
        (j-k-\ell)(N-m+j-k-\ell)-\ell(N-m+\ell), & 0 \leq \ell \leq j-k, \\
        (j-k)(j-k-N+m-2\ell), & j-k \leq \ell \leq m-k.
    \end{cases}
\end{align*}
This choice is unique, which shows that $S_{k,\ell}(q)^2$ has lowest degree $2k(m-k-2\ell)+2d_{k,\ell}$ and lowest coefficient $1$.
The lowest degrees of $([k]_q!)^2$ and $\begin{bmatrix}N\\k\end{bmatrix}_q$ are $-k(k-1)$ and $-k(N-k)$, respectively, while $(1-q^2)^{-(m-k)}$ has constant term $1$ and no negative powers. Adding these degrees and using the bound on $C_{k,\ell}$ gives the following lower bound for the degree of $P_{k,\ell}(q)$:
\begin{align*}
    &k(N-1)+2N\ell-k(k-1)-k(N-k)-2k(m-k-\ell)+2k(m-k-2\ell)+2d_{k,\ell} \\
    &\qquad = 2\ell(N-k)+2d_{k,\ell}.
\end{align*}
Substituting the two cases for $d_{k,\ell}$ gives a lower bound on the degree of $P_{k,\ell}(q)$:
\begin{align}\label{eq:two-m direct-degree-bound}
    \begin{cases}
        2\left((j-k-\ell)(N-m+j-k)+\ell(m-j)\right), & 0 \leq \ell \leq j-k, \\
        2\left((m-k-\ell)(j-k)+(\ell-j+k)(N-j)\right), & j-k \leq \ell \leq m-k.
    \end{cases}
\end{align}
Both cases are nonnegative because $N-m \geq 1$ and $m-j > 0$. They vanish only when $k = j$ and $\ell = 0$. We show that there is a choice of $k, \ell$ and $T$ that attains this lower bound.

For $k = j$ and $\ell = 0$, the only subset in Equation \eqref{eq:two-m direct-cancellation} is $T = \varnothing$, so $S_{j,0}(q) = q^{j(m-j)}$. The lowest coefficient of $C_{j,0}$ is $1$. The factorial, quantum binomial, and $(1-q^2)^{-(m-j)}$ factors also have lowest coefficient $1$, showing that the constant term of $P_{j,0}$ is $1$. Therefore:
\begin{align*}
    P(q)= q^{m(N-1)}\left\langle \mathcal{W}^j_{2m}\overline{\mathcal{W}^j_{2m}} \right\rangle \in 1 + q\Z[q].
\end{align*}
Applying Lemma \ref{WebPairIndecompo} proves that $\mathcal{W}^j_{2m}$ is indecomposable.
\end{proof}

\begin{remark}
    Already for $m=4$, only a single additional $8$-valent vertex does not suffice in constructing a web basis for the boundary sequence $\varepsilon = 1^+ 1^-1^+ 1^-1^+ 1^-1^+ 1^-$. That is due to the dimension of said web space being $24$ (for $N \geq 4$). There are $14$ planar matchings on $\varepsilon$ and $8$ ways of arranging the $6$-valent vertex and a single cup or cap, leaving one dimension unaccounted for. For this small case, it can be checked that the missing basis element is $\mathcal{W}_4^2$. This sort of a naive approach fails for $m > 4$. 
\end{remark}

We have always assumed that $r\geq0$ in order to recover the Kekulé relations. For the case $r<0$, even though the Kekulé relations are degenerated, we can ask whether $\Gamma^m_j$'s are decomposable. By using the same strategy as in the proof of Corollary~\ref{prop:two-m extreme coefficients}, we can show that $\Gamma^m_j$ is indecomposable if $r=N-m-1<0$, or equivalently, $N\leq m$, obtaining the following corollary.

\begin{corollary}
    Let $m\geq0$, $a=(0,\cdots,0)$ and $b=(1,\cdots,1)$. Then $\Gamma^m_j$ is decomposable if and only if $N>m$. 
\end{corollary}

%\section{Further directions}\label{sec:further direction}

%\newpage

\end{document}